\documentclass[a4paper,10pt]{article}
\usepackage[T2A]{fontenc}
\usepackage[utf8]{inputenc}

\usepackage[margin=2.5cm]{geometry}

\usepackage{amsmath, amssymb, amsthm}
\usepackage{hyperref}
\usepackage{xcolor}
\usepackage{color, colortbl}
\usepackage{tikz}

\usepackage[shortlabels]{enumitem}

\usepackage{chngcntr}

\makeatletter
\newcommand{\customitem}[1]{%
\item[#1]\protected@edef\@currentlabel{#1}%
}
\makeatother

\usepackage{hhline}
\usepackage{diagbox}
\usepackage{multirow}

\definecolor{myred}{rgb}{0.81, 0.06, 0.13}
\definecolor{myblue}{rgb}{0,0.5,0.6}

\newcommand{\Mod}[1]{\ (\mathrm{mod}\ #1)}

\def\lrc{\lambda_{rc}}
\def\lrr{\lambda_{rr}}
\def\lcc{\lambda_{cc}}
\def\lrrc{\lambda_{rrc}}

\def\F{{\mathbb F}}

\def\Z{{\mathbb Z}}

\DeclareMathOperator{\PG}{PG}
\DeclareMathOperator{\AG}{AG}

\DeclareMathOperator{\Aut}{Aut}

\newcommand{\cfl}{\cellcolor{black!30}}

\theoremstyle{plain}
\newtheorem{theorem}{Theorem}[section]
\newtheorem{corollary}[theorem]{Corollary}
\newtheorem{proposition}[theorem]{Proposition}
\newtheorem{lemma}[theorem]{Lemma}

\newtheorem{question}[theorem]{Question}
\newtheorem{conjecture}[theorem]{Conjecture}
\newtheorem{construction}[theorem]{Construction}
\theoremstyle{definition}
\newtheorem{definition}[theorem]{Definition}
\newtheorem{example}[theorem]{Example}
\newtheorem{problem}[theorem]{Problem}
\newtheorem{observation}[theorem]{Observation}
\theoremstyle{remark}
\newtheorem{remark}[theorem]{Remark}

\begin{document}

\title{Resolvable Triple Arrays}
\author{Alexey Gordeev, Lars-Daniel {\"O}hman\\
\small{Department of Mathematics and Mathematical Statistics, Ume\r{a} University, Sweden}}
\date{}

\maketitle

\begin{abstract}
We present a new construction of triple arrays by combining a symmetric 2-design with a resolution of another 2-design.
This is the first general method capable of producing non-extremal triple arrays.
We call the triple arrays which can be obtained in this way \emph{resolvable}.
We employ the construction to produce the first examples of $(21 \times 15, 63)$-triple arrays, and enumerate all resolvable $(7 \times 15, 35)$-triple arrays, of which there was previously only a single known example.
An infinite subfamily of Paley triple arrays turns out to be resolvable.

We also introduce a new intermediate object, \emph{unordered} triple arrays, that are to triple arrays what symmetric 2-designs are to Youden rectangles, and propose a strengthening of Agrawal's long-standing conjecture on the existence of extremal triple arrays.
For small parameters, we completely enumerate all unordered triple arrays, and use this data to corroborate the new conjecture.
We construct several infinite families of resolvable unordered triple arrays, and, in particular, show that all $((q + 1) \times q^2, q(q + 1))$-triple arrays are resolvable and are in correspondence with finite affine planes of order $q$.
\end{abstract}

\section{Introduction}\label{sec:intro}

An $(r \times c, v)$-\textit{triple array} is an $r \times c$ array on $v$ symbols that is \textit{binary} (no repeated symbols in any row or column), \textit{equireplicate} (each symbol occurs the same number of times), and satisfies the following three intersection conditions:

\begin{enumerate}
    \customitem{(RC)} Any row and column have a constant number of symbols in common,
	\customitem{(RR)} Any pair of distinct rows has a constant number of symbols in common,
	\customitem{(CC)} Any pair of distinct columns has a constant number of symbols in common.
\end{enumerate}

Triple arrays were introduced by Agrawal~\cite{agrawalMethodsConstructionDesigns1966}, though sporadic examples had appeared earlier in the literature.
For a historical overview, see Section 13 of the survey by Bailey~\cite{baileyRelationsPartitions2017}.
Triple arrays with $v = \max(r, c)$ are better known as \textit{Youden rectangles}, or, when $r = c$, as \textit{Latin squares}, and triple arrays with $v = rc$ have no repeated symbols and exist trivially for all $r$, $c$.
Regarding what other parameter combinations are possible for triple arrays, the following theorem was shown independently using linear algebraic techniques by Bailey and Heidtmann \cite[unpublished]{bailey1994extremal}, Bagchi~\cite[Corollary~1]{bagchiTwowayDesigns1998} and McSorley, Phillips, Wallis and Yucas~\cite[Theorem~3.2]{mcsorleyDoubleArraysTriple2005a}.
More recently, the present authors and Markstr{\"o}m gave an alternative purely combinatorial proof using double counting in~\cite[Corollary 7.2(a)]{gordeevTripleArrays2026}.

\begin{theorem}
Any $(r \times c, v)$-triple array with $v > \max(r, c)$ has $v \geq r + c - 1$.
\end{theorem}

The case $v = r + c - 1$ is the most well-studied, and such triple arrays are usually called \textit{extremal}.
Agrawal~\cite{agrawalMethodsConstructionDesigns1966} proposed a construction of extremal triple arrays based on a \textit{symmetric 2-design} and gave examples for a range of parameter sets, but it has yet to be proven that the construction always works.
Preece, Wallis and Yucas~\cite{preecePaleyTripleArrays2005} used Hadamard matrices to construct $(q \times (q + 1), 2q)$-triple arrays for all odd prime powers $q \geq 5$, which they named \textit{Paley triple arrays}.
Nilson and the second author~\cite{nilsonTripleArraysYouden2015} introduced a method of constructing triple arrays from \textit{Youden rectangles} and, in particular, showed that all Paley triple arrays can be constructed this way.
Nilson and Cameron~\cite{nilsonTripleArraysDifference2017} further investigated this method for Youden rectangles developed from \textit{difference sets} and constructed $((2u^2 - u) \times (2u^2 + u), 4u^2 - 1)$-triple arrays for positive integers $u$ with the square-free part dividing 6.
These arrays, together with Paley triple arrays, form the only two known infinite families.
Notably, all these methods can only produce extremal triple arrays.

In contrast, much less is known about \textit{non-extremal} triple arrays, i.e. those with $r + c - 1 < v < rc$.
For a long while, no non-extremal examples were known.
At the British Combinatorial Conference in Aberystwyth in 1973, using different terminology, Preece asked for a $(7 \times 15, 35)$-triple array (also listed by Preece as an open question in~\cite[Table~1]{preeceNonorthogonalGraecolatinDesigns1976}), which is the smallest possible non-extremal parameter set.
An example was subsequently found in the 2000s by McSorley, Phillips, Wallis and Yucas~\cite[Section 10]{mcsorleyDoubleArraysTriple2005a}, but they did not specify the details of how this design was found, and provided no analysis of its structure.
Such an analysis was, however, later given by Yucas~\cite{yucas_structure_7x15}, who also commented that the design was originally found by computer search.
Until now, this remained the only example of a non-extremal triple array in the literature.

Yucas showed how to deconstruct the aforementioned $(7 \times 15, 35)$-triple array using two key ingredients: a \textit{resolution} of the 2-design corresponding to the \textit{projective geometry} $\PG(3, 2)$, and the symmetric 2-design known as the \textit{Fano plane}.
The resolution he used is one of the seven original \textit{Kirkman parades} (see Mulder~\cite{mulder1917kirkman} and Cole~\cite{cole_kirkman_parades}), and Yucas further claimed that the same construction can be carried out using other parades, which, as we will see in Section~\ref{subsec:7x15}, is not entirely true.
He did not explicitly mention, however, that the same idea can also in principle be applied to other compatible pairs of 2-designs, one \textit{resolvable} and one symmetric, potentially resulting in triple arrays with various parameters, both non-extremal and extremal.

In the present paper, we iron out the details of a general construction method based on the ideas of Yucas~\cite{yucas_structure_7x15}, and employ it in conjunction with computational methods to produce further examples of non-extremal triple arrays, notably the first $(21\times 15,63)$-triple arrays.
The triple arrays produced by this method possess an interesting structure in addition to the basic properties of triple arrays. We call these triple arrays \emph{resolvable}, and enumerate completely all resolvable $(7 \times 15, 35)$-triple arrays.
It turns out that the method is also applicable for some extremal parameter sets, in particular, we show that an infinite subfamily of Paley triple arrays is resolvable.

The special case of triple arrays where $c=v$ corresponds to Youden rectangles. Forgetting the order of elements in the columns of a Youden rectangle, and treating these column sets as the blocks of a block design results in a symmetric 2-design. Conversely, as proven by Smith and Hartley~\cite{smithConstructionYoudenSquares1948}, all Youden rectangles can be constructed by ordering the elements of a symmetric 2-design. For general triple arrays, we analogously introduce the notion of an \emph{unordered} triple array, as a block design with two types of blocks, for the rows and for the columns, respectively. The corresponding ordering problem turns out to be considerably harder than in the Youden rectangle case.

Agrawal's construction mentioned above does produce unordered triple arrays from symmetric 2-designs, but it has not yet been proven that the construction results in triple arrays. The traditional formulation of the associated conjecture is that if there exists a symmetric 2-design, then there exists a triple array on the corresponding parameters. We generalize this by conjecturing that any extremal unordered triple array (except for one small counterexample) can in fact be ordered to give a triple array. Our enumerative work provides some evidence for the truth of the generalized conjecture. In the non-extremal case, the situation is different, as we find several $(7 \times 15, 35)$-unordered triple arrays that cannot be ordered.

Using finite projective geometries and Hadamard matrices, we construct three infinite series of \emph{resolvable unordered triple arrays}, including a non-extremal family.
Furthermore, we show that all $((q + 1) \times q^2, q(q + 1))$-triple arrays are resolvable, and that the corresponding unordered triple arrays are in one-to-one correspondence with affine planes of order $q$.
We present two equivalent reformulations of the problem of constructing a triple array on these parameters, one in terms of derangements, and another in terms of multipartite hypergraphs.
The first reformulation leads to the first combinatorial explanation for the non-existence of a $(3 \times 4, 6)$-triple array.

The rest of the paper is organized as follows. In Section~\ref{sec:prelim}, we present more formally notions from basic design theory and finite geometry, and in Section~\ref{sec:UTA} we introduce unordered triple arrays. In Section~\ref{sec:agrawal}, we present Agrawal's construction and relate it to unordered triple arrays.
Section~\ref{sec:constr} contains a new construction of unordered triple arrays, introduces the notion of resolvable triple arrays, and gives some concrete examples of families of arrays resulting from the construction. In Section~\ref{sec:paley}, we show that some Paley triple arrays are in fact resolvable. Section~\ref{sec:enum} presents our enumerative results. In Section~\ref{sec:AG->TA}, we investigate more closely triple arrays on parameters $((q+1) \times q^2, q(q+1))$. Section~\ref{sec:concl} concludes with a range of open problems.

\section{Preliminaries}
\label{sec:prelim}

\subsection{2-designs}
\label{subsec:designs}

We start with the definition and basic properties of \textit{2-designs}, also often referred to as \textit{balanced incomplete block designs (BIBD)}.
For further background on designs as well as for proofs of the claims that follow, see \emph{e.g.} the textbook of Beth, Jungnickel and Lenz~\cite{bethDesignTheory1999}.

A \textit{block design} $\mathcal{D} = (P, \mathcal{B})$ on a \textit{point} set $P$ is a collection $\mathcal{B}$ of subsets of $P$ called \textit{blocks}.
It may have \textit{identical blocks}, i.e. one set of points may appear as a block multiple times.
In other words, in general, $\mathcal{B}$ is a multiset.
The points of a block design are sometimes called \textit{treatments}, \textit{symbols}, or simply \textit{elements}.
To avoid confusion with the \emph{symbols} used in triple arrays, we will generally use the term \textit{points} when talking about block designs.

The \textit{dual} of a block design $\mathcal{D} = (P, \mathcal{B})$ is a block design $\mathcal{D}^* := (\mathcal{B}, P^*)$ with the roles of points and blocks reversed: the points of $\mathcal{D}^*$ are the blocks of $\mathcal{D}$, and, for each point $p \in P$ of $\mathcal{D}$, the dual $\mathcal{D}^*$ has a block $p^* \in P^*$ with $B \in p^*$ if and only if $p \in B$.
When talking about dual designs, we will often make no distinction between points $P$ of $\mathcal{D}$ and blocks $P^*$ of $\mathcal{D}^*$ and treat them as one and the same.

\begin{definition}
A \textit{2-$(v, k, \lambda)$-design} (\textit{2-design} for short) is a block design on a set of $v$ points in which each block has size $k$ and each pair of points occurs together in $\lambda$ blocks.
\end{definition}

The 2 in the term 2-designs refers to the \emph{pairs} mentioned in the definition.
It can be shown by double counting arguments that, setting
\begin{equation}\label{eq:2des-params}
r := \frac{\lambda(v - 1)}{k - 1},\quad b := \frac{vr}{k},    
\end{equation}
each point of a 2-$(v, k, \lambda)$ design belongs to $r$ blocks, and the total number of blocks is $b$.
We will refer to 2-$(v, k, \lambda)$ designs as \textit{2-$(v, b, r, k, \lambda)$ designs} when it is important to highlight all parameters. 
A more involved dependency between the parameters is the following.

\begin{theorem}[Fisher's inequality]\label{thm:fisher}
For any 2-$(v, b, r, k, \lambda)$ design, it holds that $b \geq v$.
\end{theorem}

The 2-designs that attain equality in Theorem~\ref{thm:fisher}, i.e. those with the same number of blocks and points, are called \textit{symmetric}.
Note that a symmetric 2-design not only has $b = v$, but also $r = k$.
The name \textit{symmetric} comes from the fact that the dual of a symmetric 2-design is also a symmetric 2-design.
In contrast, Theorem~\ref{thm:fisher} implies that the dual of a non-symmetric 2-design cannot be a 2-design.
This duality observation can be formulated as the following theorem, giving an alternative characterization of symmetric 2-designs.

\begin{theorem}\label{thm:sym-dual}
In a symmetric 2-$(v, k, \lambda)$ design, any two distinct blocks share $\lambda$ common points.
Conversely, any block design of $v$ blocks, each of size $k$, on a set of $v$ points, with two distinct blocks sharing $\lambda$ common points, is a 2-$(v, k, \lambda)$ design.
\end{theorem}

For a block design $\mathcal{D} = (P, \mathcal{B})$ and a positive integer $m$, denote by $m\mathcal{B}$ the collection of $mb$ blocks in which each block $B \in \mathcal{B}$ is repeated $m$ times.
The \textit{$m$-multiple of $\mathcal{D}$} is the block design $m\mathcal{D} := (P, m\mathcal{B})$.
Note that the $m$-multiple of a 2-$(v, b, r, k, \lambda)$ design is a 2-$(v, mb, mr, k, m\lambda)$ design.
The \textit{complement} of a block design $\mathcal{D} = (P, \mathcal{B})$ is a block design $\overline{\mathcal{D}} := (P, \overline{\mathcal{B}})$
, where $\overline{\mathcal{B}}$ is a multiset containing a block $P \setminus B$ for each block $B$ of $\mathcal{B}$.
The complement of a 2-$(v, b, r, k, \lambda)$ design is a 2-$(v, b, b - r, v - k, \lambda + b - 2r)$ design.
In particular, the complement of a symmetric 2-$(v, k, \lambda)$ design is a symmetric 2-$(v, v - k, \lambda + v - 2k)$ design.

Looking further at the internal structure of 2-designs, in particular at the prospect of partitioning the collection of blocks in a  `balanced' way, the following definition is of use.

\begin{definition}
A \textit{parallel class} in a 2-design is a set of blocks that partition the point set.
A \textit{resolution} of a 2-design is a partition of the collection of blocks into parallel classes.
A 2-design that admits a resolution is called \textit{resolvable}.
\end{definition}

Note that a parallel class in a 2-$(v, k, \lambda)$ design consists of $\frac{v}{k}$ blocks, and a resolution has $\frac{bk}{v} = r$ parallel classes.
Due to triple arrays being the main focus of the present paper, the letters $v$ and $r$ will from now on usually be reserved for, respectively, the number of symbols and rows in an array, and not for the parameters of a 2-design.

\subsection{Finite geometries}
\label{subsec:fin-geom-def}

For integer $n \geq 2$ and a prime power $q$, let $\PG(n, q)$ denote the standard $n$-dimensional \textit{projective space} over the finite field with $q$ elements. The space $\PG(n,q)$ can also be seen as a 2-design, with points and blocks of the design corresponding to points and lines of the geometry.
For integer $1 \leq i < n$, let $\PG_i(n, q)$ denote the 2-design with points and blocks corresponding to points and $i$-dimensional subspaces of $\PG(n, q)$.
Similarly, let $\AG(n, q)$ denote the standard $n$-dimensional \textit{affine space} over the finite field with $q$ elements, and $\AG_i(n, q)$ the 2-design with points and blocks corresponding to points and $i$-dimensional affine subspaces of $\AG(n, q)$.

For $n=2$ and a prime power $q$, $\PG(2, q)$ and $\AG(2, q)$ are called the standard projective and affine \emph{planes} of order $q$, or the \textit{Galois} projective and affine planes over the finite field with $q$ elements.
In greater generality, any symmetric 2-$(n^2 + n + 1, n + 1, 1)$ design is usually also called a (non-standard) \textit{finite projective plane of order $n$}, and any 2-$(n^2, n, 1)$ design is called a (non-standard) \textit{finite affine plane of order $n$}.
When talking about these designs, we will use the terms `blocks' and `lines' interchangeably.

\subsection{Triple arrays and component designs}
\label{subsec:TAs}

In order to situate triple arrays in a more general context, we state the following definition.
An $r \times c$ \textit{row-column design} on $v$ symbols is a two-dimensional array with $r$ rows and $c$ columns, each cell of which is filled with one of $v$ symbols.
It is \textit{binary} if no symbol appears more than once in any row or column. 
It is \textit{equireplicate} with \textit{replication number} $e$ if every symbol occurs $e$ times in the array. With this terminology at hand, we may define triple arrays as follows.

\begin{definition}\label{def:TA}
An \textit{$(r \times c, v)$-triple array} (TA) is a binary equireplicate $r \times c$ row-column design on $v$ symbols in which, for some integers $\lrc$, $\lrr$, $\lcc$,
\begin{enumerate}
\customitem{(RC)}\label{TA:rc} any row and column have $\lrc$ common symbols,
\customitem{(RR)}\label{TA:rr} any two distinct rows have $\lrr$ common symbols,
\customitem{(CC)}\label{TA:cc} any two distinct columns have $\lcc$ common symbols.
\end{enumerate}
\end{definition}

The number of distinct symbols $v$ used in a triple array is clearly at most the number of cells $rc$.
On the other hand, since a triple array is binary, $v$ is at least as large as the numbers $r$ and $c$ of rows and columns.
Triple arrays with $v = r = c$ are \textit{Latin squares}, and triple arrays with $v = r$ or $v = c$ are \textit{Youden rectangles} (also often, somewhat confusingly, called \textit{Youden squares}).
Triple arrays with $v = rc$ have a unique symbol in each cell and thus exist trivially for all $r, c$.
For the purposes of this paper, we generally consider only \textit{non-trivial} triple arrays, i.e. those with $\max(r, c) < v < rc$.

The main connection between 2-designs and triple arrays is that each $(r \times c, v)$-triple array $T$ has two associated 2-designs called its \textit{component designs}: the \emph{row design} $RD_T$ and the \emph{column design} $CD_T$, constructed as follows.

Let $V$ be the symbol set of $T$, and let $R_i$ and $C_j$ be the sets of symbols appearing in row $i$ and column $j$ of $T$, respectively.
Then $RD_T$ is the dual of the block design $(V, \{R_i\}_{i=1}^r)$. In other words, the points of $RD_T$ are $R_1, \dots, R_r$, and $RD_T$ has a block $RD_T(v) := \{R_i\,:\,v \in R_i\}$ for each $v \in V$.
Definition~\ref{def:TA}, and condition~\ref{TA:rr} in particular, imply that $RD_T$ is a 2-$(r, v, c, e, \lrr)$ design.
Similarly, $CD_T$ is the dual of the block design $(V, \{C_j\}_{j=1}^c)$, has a block $CD_T(v) := \{C_j\,:\,v \in C_j\}$ for each $v \in V$, and is a 2-$(c, v, r, e, \lcc)$ design.

Noting that the total number of cells in an $(r \times c, v)$-triple array is $rc$, its replication number must be $e = \frac{rc}{v}$.
The parameters $\lrc, \lrr, \lcc$ are also determined by $r, c, v$ via similar double counts.
For $\lrc$, see~\cite[Theorem 3.1]{mcsorleyDoubleArraysTriple2005a}, and for $\lrr$ and $\lcc$, apply~\eqref{eq:2des-params} to the component designs, to obtain the following identities.

\begin{equation}\label{eq:TA-params}
\lrc = e = \frac{rc}{v},\quad \lrr = \frac{c(e - 1)}{r - 1},\quad \lcc = \frac{r(e - 1)}{c - 1}.
\end{equation}

We see that for an $(r \times c, v)$-triple array to exist, the parameters $r, c, v$ must satisfy the divisibility conditions 
\[
rc \equiv 0 \Mod{v},\quad c(e - 1) \equiv 0 \Mod{r - 1},\quad r(e - 1) \equiv 0 \Mod{c - 1}.
\]
We will call parameter sets $(r \times c, v)$ satisfying these divisibility conditions \emph{admissible for triple arrays}.
Throughout the rest of the text, unless we explicitly state otherwise, $(r \times c, v)$ denotes a parameter set admissible for triple arrays, and the corresponding parameters $e, \lrc, \lrr, \lcc$ are given by~\eqref{eq:TA-params}.

\section{Unordered triple arrays}
\label{sec:UTA}

The component designs of a triple array are very useful in analyzing its structure, since they essentially represent properties~\ref{TA:rr} and~\ref{TA:cc} from Definition~\ref{def:TA}.
However, it is not as easy to recover the property~\ref{TA:rc} just by looking at the component designs.
In light of this, we introduce the following new definition.

\begin{definition}\label{def:UTA}
An $(r \times c, v)$-\textit{unordered triple array} (UTA) on a set $V$ of $v$ symbols is a collection of $c$-sets $R_1, \dots, R_r \subseteq V$ called \textit{row-sets}, and $r$-sets $C_1, \dots, C_c \subseteq V$ called \textit{column-sets} such that, for some integers $e$, $\lrc$, $\lrr$, $\lcc$, each symbol appears in $e$ row-sets and $e$ column-sets, and
\[
|R_i \cap C_j| = \lrc \text{ for all } i, j,\quad |R_i \cap R_s| = \lrr \text{ for all } i \neq s,\quad |C_j \cap C_t| = \lcc \text{ for all }j \neq t.
\]
\end{definition}

It turns out that many properties of triple arrays only depend on the contents of rows and columns and not on the exact placement of symbols in cells.
The same double counts as for triple arrays in the previous section apply for unordered triple arrays as well, and show that the parameters $e, \lrc, \lrr, \lcc$ of an unordered triple array must likewise satisfy conditions~\eqref{eq:TA-params}.
The notion of \textit{admissible} parameters thus naturally extends to unordered triple arrays.
Similarly, we will refer to $e$ as the \textit{replication number}, call an unordered triple array \textit{extremal} if $v = r + c - 1$, and \textit{non-extremal} if $r + c - 1 < v < rc$.

A consequence of Theorem~7.3 in~\cite{jagerEnumerationConstructionRowColumn2025} is that for non-trivial unordered triple arrays, that is, where $\max\{r,c\} < v < rc$, it holds that $r\neq c$, so the row-sets and the column-sets of a non-trivial unordered triple array actually have different sizes, and can therefore easily be distinguished from each other.
An unordered triple array $U$ has two associated 2-designs, the row and column designs $RD_U$ and $CD_U$, defined exactly as for triple arrays, that is, $RD_U$ is the dual of the block design $(V, \{R_i\}_{i = 1}^r)$ given by the row-sets, and $CD_U$ is the dual of the block design $(V, \{C_j\}_{j = 1}^c)$ given by the column-sets.

Any $(r \times c, v)$-triple array $T$ immediately gives rise to an $(r \times c, v)$-unordered triple array $U_T$ by treating rows and columns as unordered row-sets and column-sets.
We will call $U_T$ the \textit{underlying} unordered triple array of $T$.
Clearly, the component designs of $T$ and $U_T$ coincide.
An example of a $(4 \times 9, 12)$-triple array and a representation of its underlying unordered triple array are given in Figure~\ref{fig:4-9-12-TA}.

\begin{figure}[!ht]
\begin{center}
\begin{tabular}{c c c}
\begin{tabular}{c|*{9}{c}|}
\multicolumn{1}{c}{} & \multicolumn{1}{c}{\scriptsize $C_{1}$} & \multicolumn{1}{c}{\scriptsize $C_{2}$} & \multicolumn{1}{c}{\scriptsize $C_{3}$} & \multicolumn{1}{c}{\scriptsize $C_{4}$} & \multicolumn{1}{c}{\scriptsize $C_{5}$} & \multicolumn{1}{c}{\scriptsize $C_{6}$} & \multicolumn{1}{c}{\scriptsize $C_{7}$} & \multicolumn{1}{c}{\scriptsize $C_{8}$} & \multicolumn{1}{c}{\scriptsize $C_{9}$} \\
\hhline{~---------}
{\scriptsize $R_{ 1}$} & $ 5$ & $ 7$ & $ 6$ & $ 8$ & $11$ & $ 9$ & $ 4$ & $12$ & $10$ \\
{\scriptsize $R_{ 2}$} & $ 3$ & $ 8$ & $ 7$ & $ 5$ & $ 2$ & $ 1$ & $12$ & $ 6$ & $11$ \\
{\scriptsize $R_{ 3}$} & $ 4$ & $ 3$ & $10$ & $ 2$ & $ 7$ & $ 5$ & $ 1$ & $ 9$ & $12$ \\
{\scriptsize $R_{ 4}$} & $ 6$ & $ 9$ & $ 1$ & $10$ & $ 4$ & $11$ & $ 8$ & $ 2$ & $ 3$ \\
\hhline{~---------}
\end{tabular}
& &
\renewcommand{\arraystretch}{0.5}
\setlength\tabcolsep{1pt}
\begin{tabular}{c|*{4}{c}|*{9}{c}|}
\multicolumn{1}{c}{} & \multicolumn{1}{c}{\scriptsize $R_{1}$} & \multicolumn{1}{c}{\scriptsize $R_{2}$} & \multicolumn{1}{c}{\scriptsize $R_{3}$} & \multicolumn{1}{c}{\scriptsize $R_{4}$} & \multicolumn{1}{c}{\scriptsize $C_{1}$} & \multicolumn{1}{c}{\scriptsize $C_{2}$} & \multicolumn{1}{c}{\scriptsize $C_{3}$} & \multicolumn{1}{c}{\scriptsize $C_{4}$} & \multicolumn{1}{c}{\scriptsize $C_{5}$} & \multicolumn{1}{c}{\scriptsize $C_{6}$} & \multicolumn{1}{c}{\scriptsize $C_{7}$} & \multicolumn{1}{c}{\scriptsize $C_{8}$} & \multicolumn{1}{c}{\scriptsize $C_{9}$} \\
\hhline{~-------------}
{\scriptsize  1} & & \cfl & \cfl & \cfl & & & \cfl & & & \cfl & \cfl & & \\
{\scriptsize  2} & & \cfl & \cfl & \cfl & & & & \cfl & \cfl & & & \cfl & \\
{\scriptsize  3} & & \cfl & \cfl & \cfl & \cfl & \cfl & & & & & & & \cfl \\
{\scriptsize  4} & \cfl & & \cfl & \cfl & \cfl & & & & \cfl & & \cfl & & \\
{\scriptsize  5} & \cfl & \cfl & \cfl & & \cfl & & & \cfl & & \cfl & & & \\
{\scriptsize  6} & \cfl & \cfl & & \cfl & \cfl & & \cfl & & & & & \cfl & \\
{\scriptsize  7} & \cfl & \cfl & \cfl & & & \cfl & \cfl & & \cfl & & & & \\
{\scriptsize  8} & \cfl & \cfl & & \cfl & & \cfl & & \cfl & & & \cfl & & \\
{\scriptsize  9} & \cfl & & \cfl & \cfl & & \cfl & & & & \cfl & & \cfl & \\
{\scriptsize 10} & \cfl & & \cfl & \cfl & & & \cfl & \cfl & & & & & \cfl \\
{\scriptsize 11} & \cfl & \cfl & & \cfl & & & & & \cfl & \cfl & & & \cfl \\
{\scriptsize 12} & \cfl & \cfl & \cfl & & & & & & & & \cfl & \cfl & \cfl \\
\hhline{~-------------}
\end{tabular}
\end{tabular}
\end{center}
\caption{A $(4 \times 9, 12)$-triple array and a representation of its underlying unordered triple array: the cell at the intersection of row $x$ and column $R_i$ ($C_j$) is grey when symbol $x$ belongs to row-set $R_i$ (column-set $C_j$). Representations of unordered triple arrays in later figures follow the same convention.}
\label{fig:4-9-12-TA}
\end{figure}

One of the main problems in the study of triple arrays is the construction of a triple array on a given set of admissible parameters.
Since any triple array $T$ has an underlying unordered triple array $U_T$, one way to attack this problem is to split it into two
separate sequential problems:

\begin{problem}[UTA construction problem]\label{pbm:exUTA}
Given an admissible parameter set $(r \times c, v)$, find an $(r \times c, v)$-unordered triple array.
\end{problem}

\begin{problem}[UTA ordering problem]\label{pbm:UTAtoTA}
Given an $(r \times c, v)$-unordered triple array $U$, find an $(r \times c, v)$-triple array $T$ with $U_T = U$.
\end{problem}

\begin{figure}[!ht]
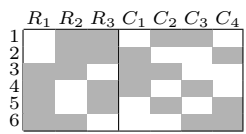

\begin{center}
\renewcommand{\arraystretch}{0.5}
\setlength\tabcolsep{1pt}
\begin{tabular}{c|*{3}{c}|*{4}{c}|}
\multicolumn{1}{c}{} & \multicolumn{1}{c}{\scriptsize $R_{1}$} & \multicolumn{1}{c}{\scriptsize $R_{2}$} & \multicolumn{1}{c}{\scriptsize $R_{3}$} & \multicolumn{1}{c}{\scriptsize $C_{1}$} & \multicolumn{1}{c}{\scriptsize $C_{2}$} & \multicolumn{1}{c}{\scriptsize $C_{3}$} & \multicolumn{1}{c}{\scriptsize $C_{4}$} \\
\hhline{~-------}
{\scriptsize  1} & & \cfl & \cfl & & \cfl & \cfl & \\
{\scriptsize  2} & & \cfl & \cfl & \cfl & & & \cfl \\
{\scriptsize  3} & \cfl & \cfl & & \cfl & \cfl & & \\
{\scriptsize  4} & \cfl & & \cfl & \cfl & & \cfl & \\
{\scriptsize  5} & \cfl & & \cfl & & \cfl & & \cfl \\
{\scriptsize  6} & \cfl & \cfl & & & & \cfl & \cfl \\
\hhline{~-------}
\end{tabular}
\end{center}
\caption{A $(3 \times 4, 6)$-unordered triple array.  Note that there are no $(3 \times 4, 6)$-triple arrays.}
\label{fig:UTA-no-sol}
\end{figure}

Both these problems turn out to be difficult and, in general, do not always have any solutions.
For instances of the \nameref{pbm:exUTA} with no solutions, see Examples~\ref{ex:7-15-21} and~\ref{ex:15-28-42}.
As for the \nameref{pbm:UTAtoTA}, it has long been known, via exhaustive checking, that there are no $(3 \times 4, 6)$-triple arrays, but
a $(3 \times 4, 6)$-unordered triple array is given in Figure~\ref{fig:UTA-no-sol}. We offer a novel combinatorial explanation of the non-existence of $(3 \times 4, 6)$-triple arrays in Corollary~\ref{cor:no346} of Section~\ref{sec:AG->TA}.
In Observation~\ref{obs:71535-no-sol} of Section~\ref{subsec:7x15}, we present more instances of the \nameref{pbm:UTAtoTA} with no solutions.

For an $(r \times c, v)$-unordered triple array $U$ with row-sets $R_1, \dots, R_r$ and columns-sets $C_1, \dots, C_c$, the \nameref{pbm:UTAtoTA} can be viewed as a special case of the following problem, by setting $A_{ij} := R_i \cap C_j$.

\begin{problem}\label{pbm:ADR}
Given a collection of finite sets $A_{ij}$, $1 \leq i \leq r$, $1 \leq j \leq c$, find an \textit{array of distinct representatives} $a_{ij} \in A_{ij}$ such that $a_{ij} \neq a_{sj}$ for $i \neq s$ and $a_{ij} \neq a_{it}$ for $j \neq t$.
\end{problem}

Fon-Der-Flaass~\cite{fon-der-flaassArraysDistinctRepresentatives1997} showed that the decision version of Problem~\ref{pbm:ADR} (i.e. the problem of deciding whether there is an array of distinct representatives $a_{ij} \in A_{ij}$ with the required properties for a given collection of finite sets $A_{ij}$) is NP-complete, even if the sets $A_{ij}$ are restricted to have form $A_{ij} = R_i \cap C_j$ for some sets $R_1, \dots, R_r$, $C_1, \dots, C_c$.
It should be noted, though, that he proved this by reducing any instance of the problem 3-SAT to an instance of Problem~\ref{pbm:ADR} in which there are sets $A_{ij}$ of sizes both 2 and 3.
This does not therefore imply NP-completeness of the decision version of the \nameref{pbm:UTAtoTA}, since in this instance we have $|A_{ij}| = |R_i \cap C_j| = \lrc$, i.e. a constant, and, moreover, there are further restrictions on the sets $R_i$, $C_j$.
In Section~\ref{sec:enum}, we explore the \nameref{pbm:UTAtoTA} computationally by treating it as a special case of another, more well-known NP-complete problem, the \textit{exact cover problem}.

To put the \nameref{pbm:exUTA} and the \nameref{pbm:UTAtoTA} into some further context, for the trivial parameters $v=r=c$, the \nameref{pbm:exUTA} is uniquely (up to renaming of symbols) solved by setting all column-sets and row-sets to $\{1,2,\ldots,v\}$, and the \nameref{pbm:UTAtoTA} is the problem of finding a $v\times v$ Latin square. Existence in this case is trivial, and research efforts have instead focused on other issues, such as how many different solutions there are. In the case $v=c>r$, the \nameref{pbm:exUTA} is equivalent to the problem of finding a symmetric 2-design, and the \nameref{pbm:UTAtoTA} is the problem of producing a Youden rectangle from a symmetric 2-design, which is always possible, as proven constructively by Smith and Hartley~\cite{smithConstructionYoudenSquares1948}. This goes to show that the \nameref{pbm:UTAtoTA} is explicitly solvable in at least some special cases.

\section{Agrawal's construction}
\label{sec:agrawal}

In this section, we consider the \nameref{pbm:exUTA} for extremal parameter sets.
Let $(r \times c, v)$ be an extremal admissible parameter set, i.e. $v = r + c - 1$.
Recalling~\eqref{eq:TA-params}, we have 
\[
e - 1 = \frac{rc}{r + c - 1} - 1 = \frac{(r - 1)(c - 1)}{v},
\]
from which we get
\begin{equation}\label{eq:extr-params}
\lcc = \frac{r(e - 1)}{c - 1} = \frac{r(r - 1)}{v} = \frac{r(v - c)}{v} = r - e,\quad \text{ and, similarly, }\quad \lrr = c - e.
\end{equation}

In 1966, Agrawal~\cite{agrawalMethodsConstructionDesigns1966} proposed a method for constructing `designs for two-way elimination of heterogeneity', now known as extremal triple arrays.
Although phrased differently, the construction proceeded in essence by first giving a recipe for constructing an unordered triple array from a symmetric 2-design, as follows. 

\begin{construction}[Agrawal's construction]\label{constr:Agrawal}
Start with a symmetric 2-$(r + c, r, \lcc)$ design $\mathcal{S}$ on the point set $V$, and fix one of its points $\sigma \in V$.
Let $C_1, \dots, C_c$ be the blocks of $\mathcal{S}$ not containing $\sigma$, and let $R_1, \dots, R_r$ be the complements of blocks $\overline{R_1}, \dots, \overline{R_r}$ of $\mathcal{S}$ containing $\sigma$, i.e. $R_i := V \setminus \overline{R_i}$.
\end{construction}

It is clear that given a symmetric 2-design, the construction can be carried out, and it is not hard to show that the resulting sets form an unordered triple array.

\begin{proposition}
Taking $R_i$ and $C_j$ from \nameref{constr:Agrawal} as row-sets and column-sets, respectively, gives an $(r \times c, r + c - 1)$-unordered triple array.        
\end{proposition}

\begin{proof}
Using Theorem~\ref{thm:sym-dual} and~\eqref{eq:extr-params}, we have, for all $i, j, s, t$ with $i \neq s$ and $j \neq t$,
\begin{align*}
|C_j \cap C_t| & = \lcc, \\    
|R_i \cap R_s| = |V| - |\overline{R_i}| - |\overline{R_s}| + |\overline{R_i} \cap \overline{R_s}| = (r + c) - 2r + \lcc & = \lrr, \\
|R_i \cap C_j| = |C_j| - |\overline{R_i} \cap C_j| = r - \lcc = e & = \lrc.\qedhere
\end{align*}

\end{proof}

\begin{example}\label{ex:agrawal_constr}
Applying \nameref{constr:Agrawal} to the symmetric 2-$(7, 3, 1)$ design with blocks
\[
\{0, 1, 2\}, \{0, 4, 5\}, \{0, 3, 6\}, \{2, 3, 4\}, \{1, 3, 5\}, \{1, 4, 6\}, \{2, 5, 6\},
\]
that is, the Fano plane, using $\sigma = 0$ as the fixed point, results in the $(3 \times 4, 6)$-unordered triple array from Figure~\ref{fig:UTA-no-sol}.
\end{example}

In our terminology, Agrawal thus only solved the problem of finding extremal $(r \times c, r + c - 1)$-\emph{unordered} triple arrays starting from any symmetric 2-$(r + c, r, \lcc)$ design.
In order to reach the end goal of constructing an $(r \times c, r + c - 1)$-triple array, it remains to solve the \nameref{pbm:UTAtoTA}. 
Agrawal did not provide a general way of doing this, but he remarked that he was able to find a solution in all the examples he considered, provided that $e = r - \lcc > 2$. The following conjecture, which still remains open, was therefore implicit in Agrawal's work.

\begin{conjecture}[Agrawal~\cite{agrawalMethodsConstructionDesigns1966}]\label{conj:Agrawal}
If there is a symmetric 2-$(r + c, r, \lcc)$ design with $r - \lcc > 2$, then there is an $(r \times c, r + c - 1)$-triple array.
\end{conjecture}

The only parameter set excluded by the condition $r - \lcc > 2$ in Conjecture~\ref{conj:Agrawal}, is $(3 \times 4,6)$, but as seen in Example~\ref{ex:agrawal_constr} the construction of an unordered triple array works in this case as well.
Conjecture~\ref{conj:Agrawal} is especially important due to the fact that \nameref{constr:Agrawal} of an unordered triple array can be reversed: let $R_1, \dots, R_r$ and $C_1, \dots, C_c$ be the row-sets and column-sets of an extremal $(r \times c, v)$-unordered triple array on a symbol set $V$.
Introduce a new symbol $\sigma$ and set $\overline{R_i} := V \setminus R_i \cup \{\sigma\}$.
Then the blocks $C_1, \dots, C_c, \overline{R_1}, \dots, \overline{R_r}$ form a symmetric 2-$(r + c, r, \lcc)$ design, which follows by repeating the arguments above in reverse. 
This implies the following theorem about extremal unordered triple arrays, first shown by Bailey and Heidtmann~\cite[unpublished]{bailey1994extremal} and later independently by McSorley, Phillips, Wallis and Yucas~\cite[Theorem 5.2]{mcsorleyDoubleArraysTriple2005a}.

\begin{theorem}\label{thm:Agrawal-univ}
Any $(r \times c, r + c - 1)$-unordered triple array can be obtained via \nameref{constr:Agrawal}.
In particular, if there is an $(r \times c, r + c - 1)$-(unordered) triple array, then there is a symmetric 2-$(r + c, r, \lcc)$ design.
\end{theorem}

\nameref{constr:Agrawal} and Theorem~\ref{thm:Agrawal-univ} essentially reduce the \nameref{pbm:exUTA} in the extremal case to the widely studied, though in general still very much open question of existence of symmetric 2-designs.
In particular, Theorem~\ref{thm:Agrawal-univ} implies non-existence of extremal (unordered) triple arrays when non-existence of the corresponding symmetric 2-designs is known.
The following two examples were given in~\cite[Section 7]{mcsorleyDoubleArraysTriple2005a}.

\begin{example}\label{ex:7-15-21}
Even though the parameter set $(7 \times 15, 21)$ is admissible for triple arrays, there are no symmetric 2-$(22, 7, 2)$ designs, and thus, by Theorem~\ref{thm:Agrawal-univ}, no $(7 \times 15, 21)$-(unordered) triple arrays.
\end{example}

\begin{example}\label{ex:15-28-42}
Similarly, the parameter set $(15 \times 28, 42)$ is admissible for triple arrays, and, moreover, there are known examples of would-be component designs, i.e. 2-$(15, 10, 18)$ and 2-$(28, 10, 5)$ designs.
However, there are no symmetric 2-$(43, 15, 5)$ designs and thus no $(15 \times 28, 42)$-(unordered) triple arrays.
\end{example}

Since any triple array has an underlying unordered triple array, Theorem~\ref{thm:Agrawal-univ} in turn implies that any extremal triple array can be obtained from \emph{some} unordered triple array obtained from Agrawal's construction. It may, however, be the case that some extremal unordered triple arrays cannot be ordered. 
The following strengthening of Conjecture~\ref{conj:Agrawal} therefore seems to better reflect Agrawal's observation that in the cases he tried, ``such rearrangement is always possible''.

\begin{conjecture}\label{conj:Agrawal-gen}
For any $(r \times c, r + c - 1)$-unordered triple array $U$ with $e > 2$, there is an $(r \times c, r + c - 1)$-triple array $T$ with $U_T = U$.
\end{conjecture}

\section{Triple arrays via resolutions}
\label{sec:constr}

This section contains our main contribution, namely the first general method of solving the \nameref{pbm:exUTA} for some non-extremal parameter sets.

\subsection{Resolvable triple arrays}
\label{subsec:main-constr}

\begin{figure}[!ht]
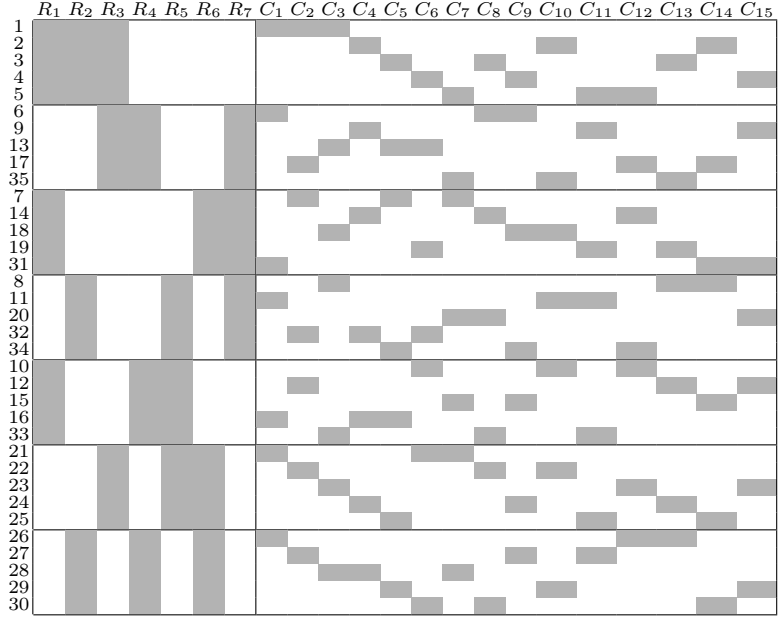

\begin{center}
\begin{tabular}{c}
\begin{tabular}{c|*{15}{c}|}
\multicolumn{1}{c}{} & \multicolumn{1}{c}{\scriptsize $C_{1}$} & \multicolumn{1}{c}{\scriptsize $C_{2}$} & \multicolumn{1}{c}{\scriptsize $C_{3}$} & \multicolumn{1}{c}{\scriptsize $C_{4}$} & \multicolumn{1}{c}{\scriptsize $C_{5}$} & \multicolumn{1}{c}{\scriptsize $C_{6}$} & \multicolumn{1}{c}{\scriptsize $C_{7}$} & \multicolumn{1}{c}{\scriptsize $C_{8}$} & \multicolumn{1}{c}{\scriptsize $C_{9}$} & \multicolumn{1}{c}{\scriptsize $C_{10}$} & \multicolumn{1}{c}{\scriptsize $C_{11}$} & \multicolumn{1}{c}{\scriptsize $C_{12}$} & \multicolumn{1}{c}{\scriptsize $C_{13}$} & \multicolumn{1}{c}{\scriptsize $C_{14}$} & \multicolumn{1}{c}{\scriptsize $C_{15}$} \\
\hhline{~---------------}
{\scriptsize $R_{ 1}$} & $31$ & $ 1$ & $18$ & $16$ & $ 7$ & $10$ & $ 5$ & $ 3$ & $ 4$ & $ 2$ & $33$ & $14$ & $19$ & $15$ & $12$ \\
{\scriptsize $R_{ 2}$} & $26$ & $32$ & $ 1$ & $ 2$ & $29$ & $30$ & $28$ & $20$ & $27$ & $11$ & $ 5$ & $34$ & $ 3$ & $ 8$ & $ 4$ \\
{\scriptsize $R_{ 3}$} & $ 1$ & $17$ & $13$ & $ 9$ & $ 3$ & $ 4$ & $21$ & $22$ & $ 6$ & $35$ & $25$ & $ 5$ & $24$ & $ 2$ & $23$ \\
{\scriptsize $R_{ 4}$} & $ 6$ & $27$ & $33$ & $28$ & $16$ & $13$ & $35$ & $30$ & $15$ & $10$ & $ 9$ & $26$ & $12$ & $17$ & $29$ \\
{\scriptsize $R_{ 5}$} & $16$ & $12$ & $23$ & $32$ & $34$ & $21$ & $15$ & $33$ & $24$ & $22$ & $11$ & $10$ & $ 8$ & $25$ & $20$ \\
{\scriptsize $R_{ 6}$} & $21$ & $22$ & $28$ & $24$ & $25$ & $19$ & $ 7$ & $14$ & $18$ & $29$ & $27$ & $23$ & $26$ & $30$ & $31$ \\
{\scriptsize $R_{ 7}$} & $11$ & $ 7$ & $ 8$ & $14$ & $13$ & $32$ & $20$ & $ 6$ & $34$ & $18$ & $19$ & $17$ & $35$ & $31$ & $ 9$ \\
\hhline{~---------------}
\end{tabular}\\
\\
\renewcommand{\arraystretch}{0.5}
\setlength\tabcolsep{1pt}
\begin{tabular}{c|*{7}{c}|*{15}{c}|}
\multicolumn{1}{c}{} & \multicolumn{1}{c}{\scriptsize $R_{1}$} & \multicolumn{1}{c}{\scriptsize $R_{2}$} & \multicolumn{1}{c}{\scriptsize $R_{3}$} & \multicolumn{1}{c}{\scriptsize $R_{4}$} & \multicolumn{1}{c}{\scriptsize $R_{5}$} & \multicolumn{1}{c}{\scriptsize $R_{6}$} & \multicolumn{1}{c}{\scriptsize $R_{7}$} & \multicolumn{1}{c}{\scriptsize $C_{1}$} & \multicolumn{1}{c}{\scriptsize $C_{2}$} & \multicolumn{1}{c}{\scriptsize $C_{3}$} & \multicolumn{1}{c}{\scriptsize $C_{4}$} & \multicolumn{1}{c}{\scriptsize $C_{5}$} & \multicolumn{1}{c}{\scriptsize $C_{6}$} & \multicolumn{1}{c}{\scriptsize $C_{7}$} & \multicolumn{1}{c}{\scriptsize $C_{8}$} & \multicolumn{1}{c}{\scriptsize $C_{9}$} & \multicolumn{1}{c}{\scriptsize $C_{10}$} & \multicolumn{1}{c}{\scriptsize $C_{11}$} & \multicolumn{1}{c}{\scriptsize $C_{12}$} & \multicolumn{1}{c}{\scriptsize $C_{13}$} & \multicolumn{1}{c}{\scriptsize $C_{14}$} & \multicolumn{1}{c}{\scriptsize $C_{15}$} \\
\hhline{~----------------------}
{\scriptsize  1} & \cfl & \cfl & \cfl & & & & & \cfl & \cfl & \cfl & & & & & & & & & & & & \\
{\scriptsize  2} & \cfl & \cfl & \cfl & & & & & & & & \cfl & & & & & & \cfl & & & & \cfl & \\
{\scriptsize  3} & \cfl & \cfl & \cfl & & & & & & & & & \cfl & & & \cfl & & & & & \cfl & & \\
{\scriptsize  4} & \cfl & \cfl & \cfl & & & & & & & & & & \cfl & & & \cfl & & & & & & \cfl \\
{\scriptsize  5} & \cfl & \cfl & \cfl & & & & & & & & & & & \cfl & & & & \cfl & \cfl & & & \\
\hhline{~----------------------}
{\scriptsize  6} & & & \cfl & \cfl & & & \cfl & \cfl & & & & & & & \cfl & \cfl & & & & & & \\
{\scriptsize  9} & & & \cfl & \cfl & & & \cfl & & & & \cfl & & & & & & & \cfl & & & & \cfl \\
{\scriptsize 13} & & & \cfl & \cfl & & & \cfl & & & \cfl & & \cfl & \cfl & & & & & & & & & \\
{\scriptsize 17} & & & \cfl & \cfl & & & \cfl & & \cfl & & & & & & & & & & \cfl & & \cfl & \\
{\scriptsize 35} & & & \cfl & \cfl & & & \cfl & & & & & & & \cfl & & & \cfl & & & \cfl & & \\
\hhline{~----------------------}
{\scriptsize  7} & \cfl & & & & & \cfl & \cfl & & \cfl & & & \cfl & & \cfl & & & & & & & & \\
{\scriptsize 14} & \cfl & & & & & \cfl & \cfl & & & & \cfl & & & & \cfl & & & & \cfl & & & \\
{\scriptsize 18} & \cfl & & & & & \cfl & \cfl & & & \cfl & & & & & & \cfl & \cfl & & & & & \\
{\scriptsize 19} & \cfl & & & & & \cfl & \cfl & & & & & & \cfl & & & & & \cfl & & \cfl & & \\
{\scriptsize 31} & \cfl & & & & & \cfl & \cfl & \cfl & & & & & & & & & & & & & \cfl & \cfl \\
\hhline{~----------------------}
{\scriptsize  8} & & \cfl & & & \cfl & & \cfl & & & \cfl & & & & & & & & & & \cfl & \cfl & \\
{\scriptsize 11} & & \cfl & & & \cfl & & \cfl & \cfl & & & & & & & & & \cfl & \cfl & & & & \\
{\scriptsize 20} & & \cfl & & & \cfl & & \cfl & & & & & & & \cfl & \cfl & & & & & & & \cfl \\
{\scriptsize 32} & & \cfl & & & \cfl & & \cfl & & \cfl & & \cfl & & \cfl & & & & & & & & & \\
{\scriptsize 34} & & \cfl & & & \cfl & & \cfl & & & & & \cfl & & & & \cfl & & & \cfl & & & \\
\hhline{~----------------------}
{\scriptsize 10} & \cfl & & & \cfl & \cfl & & & & & & & & \cfl & & & & \cfl & & \cfl & & & \\
{\scriptsize 12} & \cfl & & & \cfl & \cfl & & & & \cfl & & & & & & & & & & & \cfl & & \cfl \\
{\scriptsize 15} & \cfl & & & \cfl & \cfl & & & & & & & & & \cfl & & \cfl & & & & & \cfl & \\
{\scriptsize 16} & \cfl & & & \cfl & \cfl & & & \cfl & & & \cfl & \cfl & & & & & & & & & & \\
{\scriptsize 33} & \cfl & & & \cfl & \cfl & & & & & \cfl & & & & & \cfl & & & \cfl & & & & \\
\hhline{~----------------------}
{\scriptsize 21} & & & \cfl & & \cfl & \cfl & & \cfl & & & & & \cfl & \cfl & & & & & & & & \\
{\scriptsize 22} & & & \cfl & & \cfl & \cfl & & & \cfl & & & & & & \cfl & & \cfl & & & & & \\
{\scriptsize 23} & & & \cfl & & \cfl & \cfl & & & & \cfl & & & & & & & & & \cfl & & & \cfl \\
{\scriptsize 24} & & & \cfl & & \cfl & \cfl & & & & & \cfl & & & & & \cfl & & & & \cfl & & \\
{\scriptsize 25} & & & \cfl & & \cfl & \cfl & & & & & & \cfl & & & & & & \cfl & & & \cfl & \\
\hhline{~----------------------}
{\scriptsize 26} & & \cfl & & \cfl & & \cfl & & \cfl & & & & & & & & & & & \cfl & \cfl & & \\
{\scriptsize 27} & & \cfl & & \cfl & & \cfl & & & \cfl & & & & & & & \cfl & & \cfl & & & & \\
{\scriptsize 28} & & \cfl & & \cfl & & \cfl & & & & \cfl & \cfl & & & \cfl & & & & & & & & \\
{\scriptsize 29} & & \cfl & & \cfl & & \cfl & & & & & & \cfl & & & & & \cfl & & & & & \cfl \\
{\scriptsize 30} & & \cfl & & \cfl & & \cfl & & & & & & & \cfl & & \cfl & & & & & & \cfl & \\
\hhline{~----------------------}
\end{tabular}
\end{tabular}
\end{center}
\caption{The $(7 \times 15, 35)$-triple array from~\cite{mcsorleyDoubleArraysTriple2005a,yucas_structure_7x15} and a representation of its underlying unordered triple array. 
}
\label{fig:MPWY-71535}
\end{figure}

Yucas~\cite{yucas_structure_7x15} provided an analysis of the internal 
structure of a particular $(7 \times 15, 35)$-triple array, by examining its component designs.
The array in question is given in Figure~\ref{fig:MPWY-71535}, together with a representation of its underlying unordered triple array which highlights its internal structure: the symbols are partitioned into 7 groups of 5 symbols, with symbols in one group all appearing together in the same row-sets, and each column-set containing precisely one symbol from each group.
To make this into a general construction, the analysis has to be reversed, as in the following construction, whose name reflects the fact that we call the resulting design a resolvable unordered triple array (RUTA; see Definition~\ref{def:RUTA} below).

\begin{construction}[RUTA construction]\label{constr:res}
Let $(r \times c, v)$ be a parameter set admissible for triple arrays such that, in addition to $e, \lrc, \lrr, \lcc$, the following two parameters are integers:
\begin{equation}\label{eq:lrrc-k}
\lrrc := \frac{e(e - 1)}{r - 1} \in \Z,\quad k := \frac{c}{e} \in \Z.    
\end{equation}

Start with a symmetric 2-$(r, r, e, e, \lrrc)$ design $\mathcal{S}$ on the point set $\{1, \dots, r\}$ with blocks labeled $S_1, \dots, S_r$, and a resolution $\mathcal{B}$ of a 2-$(c, v, r, e, \lcc)$ design on the point set $\{1, \dots, c\}$.
Denote the parallel classes of the resolution by $\mathcal{B}_1, \dots, \mathcal{B}_r$, and denote the blocks comprising a parallel class by $\mathcal{B}_x = \{B_{x1}, \dots, B_{xk}\}$.
Treating the blocks $B_{xy}$ as symbols, define row-sets $R_i$ and column-sets $C_j$ as
\[
R_i := \bigcup_{x\,:\,i \in S_x} \mathcal{B}_x,\quad C_j := \{ B_{xy}\ :\ j \in B_{xy} \}.
\]
\end{construction}

As in \nameref{constr:Agrawal}, it is clear that the construction can be carried out, given a symmetric 2-design and a resolvable 2-design.

\begin{proposition}\label{prop:res-correct}
    The sets $R_i$ and $C_j$ from the \nameref{constr:res} form an $(r \times c, v)$-unordered triple array.
\end{proposition}
\begin{proof}
Indeed, by definition, $|R_i| = ek = c$ and $|C_j| = r$, each symbol $B_{xy}$ appears in $e$ row-sets and $e$ column-sets, and, for all $j \neq t$ and $i \neq s$, we have
\[
|C_j \cap C_t| = \lcc,\quad |R_i \cap R_s| = \sum_{x\,:\,\{i,s\} \subseteq S_x} |\mathcal{B}_x| = \lrrc k = \frac{(e - 1)c}{r - 1} = \lrr.
\]
Finally, for each parallel class $\mathcal{B}_x$ we have $|\mathcal{B}_x \cap C_j| = 1$, thus, for all $i, j$,
\[
|R_i \cap C_j| = \sum_{x\,:\,i \in S_x}|\mathcal{B}_x \cap C_j| = e = \lrc.\qedhere
\]
\end{proof}

The choice of notation for the parameter $\lrrc$ is explained in Section~\ref{subsec:quad}.
Note that the column design $CD_U$ is the same as $\mathcal{B}$, and the row design $RD_U$ is the same as $k\mathcal{S}$, the $k$-multiple of $\mathcal{S}$.
More precisely, for each $1 \leq x \leq r$, the blocks $CD_U(B_{x1}), \dots, CD_U(B_{xk}) \in \mathcal{B}_x$ form a parallel class in $CD_U$, and the blocks $RD_U(B_{x1}), \dots, RD_U(B_{xk})$ are identical, each corresponding to the block $\mathcal{S}_x$ in $\mathcal{S}$.
This motivates the following definition.

\begin{definition}\label{def:RUTA}
An $(r \times c, v)$-unordered triple array is \textit{resolvable} (RUTA) if $\lrrc$ and $k$ defined as in~\eqref{eq:lrrc-k} are integers, and the symbols can be partitioned into $r$ groups of size $k$, with each column-set containing precisely one symbol from each group, and symbols in one group all appearing in the same row-sets.
An $(r \times c, v)$-triple array is \textit{resolvable} (RTA) if its underlying unordered triple array is resolvable.

A parameter set $(r \times c, v)$ is \textit{admissible for resolvable triple arrays} if it is admissible for triple arrays and, additionally, $\lrrc$ and $k$ are integers.
\end{definition}

\begin{remark}
Any unordered triple array obtained from the \nameref{constr:res} is clearly resolvable. 
On the other hand, for any resolvable unordered triple array $T$, its row design $RD_T$ is the $k$-multiple of some symmetric 2-$(r, e, \lrrc)$ design $\mathcal{S}$, i.e. $RD_T = k\mathcal{S}$, and the groups of symbols induce a resolution $\mathcal{B}$ on its column design $CD_T$.
Then $T$ can be obtained by using the \nameref{constr:res} with these $\mathcal{S}$ and $\mathcal{B}$.
\end{remark}

Rows and columns are assigned different roles in Definition~\ref{def:RUTA}, which raises the question of the possibility of transposing a resolvable (unordered) triple array and again receiving such an array. In the case of (unordered) triple arrays, the transpose is always an (unordered) triple array, but it follows from the next lemma that an (unordered) triple array can only be resolvable with regards to one orientation.

\begin{lemma}\label{lm:res-transpose}
If a non-trivial parameter set $(r \times c, v)$ is admissible for resolvable triple arrays, then $(c \times r, v)$ is not admissible for resolvable triple arrays.
\end{lemma}
\begin{proof}
Assume, for a contradiction, that $(c \times r, v)$ is admissible. Then $e$ must divide $r$, so $r - 1$ is coprime with $e$. Thus $e$ divides $\lrrc = \frac{e(e - 1)}{r - 1}$, and consequently $e \leq \lrrc$.
On the other hand, since $e < r$, and equivalently $e-1 < r-1$, it follows that $\lrrc = \frac{e(e - 1)}{r - 1} < e$, a contradiction.
\end{proof}

\begin{figure}[!ht]
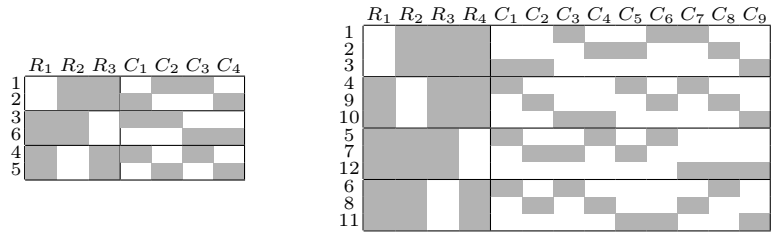

\begin{center}
\begin{tabular}{c c c}
\renewcommand{\arraystretch}{0.5}
\setlength\tabcolsep{1pt}
\begin{tabular}{c|*{3}{c}|*{4}{c}|}
\multicolumn{1}{c}{} & \multicolumn{1}{c}{\scriptsize $R_{1}$} & \multicolumn{1}{c}{\scriptsize $R_{2}$} & \multicolumn{1}{c}{\scriptsize $R_{3}$} & \multicolumn{1}{c}{\scriptsize $C_{1}$} & \multicolumn{1}{c}{\scriptsize $C_{2}$} & \multicolumn{1}{c}{\scriptsize $C_{3}$} & \multicolumn{1}{c}{\scriptsize $C_{4}$} \\
\hhline{~-------}
{\scriptsize  1} & & \cfl & \cfl & & \cfl & \cfl & \\
{\scriptsize  2} & & \cfl & \cfl & \cfl & & & \cfl \\
\hhline{~-------}
{\scriptsize  3} & \cfl & \cfl & & \cfl & \cfl & & \\
{\scriptsize  6} & \cfl & \cfl & & & & \cfl & \cfl \\
\hhline{~-------}
{\scriptsize  4} & \cfl & & \cfl & \cfl & & \cfl & \\
{\scriptsize  5} & \cfl & & \cfl & & \cfl & & \cfl \\
\hhline{~-------}
\end{tabular}
& &
\renewcommand{\arraystretch}{0.5}
\setlength\tabcolsep{1pt}
\begin{tabular}{c|*{4}{c}|*{9}{c}|}
\multicolumn{1}{c}{} & \multicolumn{1}{c}{\scriptsize $R_{1}$} & \multicolumn{1}{c}{\scriptsize $R_{2}$} & \multicolumn{1}{c}{\scriptsize $R_{3}$} & \multicolumn{1}{c}{\scriptsize $R_{4}$} & \multicolumn{1}{c}{\scriptsize $C_{1}$} & \multicolumn{1}{c}{\scriptsize $C_{2}$} & \multicolumn{1}{c}{\scriptsize $C_{3}$} & \multicolumn{1}{c}{\scriptsize $C_{4}$} & \multicolumn{1}{c}{\scriptsize $C_{5}$} & \multicolumn{1}{c}{\scriptsize $C_{6}$} & \multicolumn{1}{c}{\scriptsize $C_{7}$} & \multicolumn{1}{c}{\scriptsize $C_{8}$} & \multicolumn{1}{c}{\scriptsize $C_{9}$} \\
\hhline{~-------------}
{\scriptsize  1} & & \cfl & \cfl & \cfl & & & \cfl & & & \cfl & \cfl & & \\
{\scriptsize  2} & & \cfl & \cfl & \cfl & & & & \cfl & \cfl & & & \cfl & \\
{\scriptsize  3} & & \cfl & \cfl & \cfl & \cfl & \cfl & & & & & & & \cfl \\
\hhline{~-------------}
{\scriptsize  4} & \cfl & & \cfl & \cfl & \cfl & & & & \cfl & & \cfl & & \\
{\scriptsize  9} & \cfl & & \cfl & \cfl & & \cfl & & & & \cfl & & \cfl & \\
{\scriptsize 10} & \cfl & & \cfl & \cfl & & & \cfl & \cfl & & & & & \cfl \\
\hhline{~-------------}
{\scriptsize  5} & \cfl & \cfl & \cfl & & \cfl & & & \cfl & & \cfl & & & \\
{\scriptsize  7} & \cfl & \cfl & \cfl & & & \cfl & \cfl & & \cfl & & & & \\
{\scriptsize 12} & \cfl & \cfl & \cfl & & & & & & & & \cfl & \cfl & \cfl \\
\hhline{~-------------}
{\scriptsize  6} & \cfl & \cfl & & \cfl & \cfl & & \cfl & & & & & \cfl & \\
{\scriptsize  8} & \cfl & \cfl & & \cfl & & \cfl & & \cfl & & & \cfl & & \\
{\scriptsize 11} & \cfl & \cfl & & \cfl & & & & & \cfl & \cfl & & & \cfl \\
\hhline{~-------------}
\end{tabular}
\end{tabular}
\end{center}
\caption{A $(3 \times 4, 6)$ and a $(4 \times 9, 12)$ resolvable unordered triple arrays.}
\label{fig:UTA-Res-346-4912}
\end{figure}

\begin{example}\label{ex:346-4912-res}
The unordered triple arrays from Figures~\ref{fig:4-9-12-TA} and~\ref{fig:UTA-no-sol} are both resolvable.
Alternative representations with symbols rearranged to highlight their structure are given in Figure~\ref{fig:UTA-Res-346-4912}. 
\end{example}

\begin{figure}[!ht]
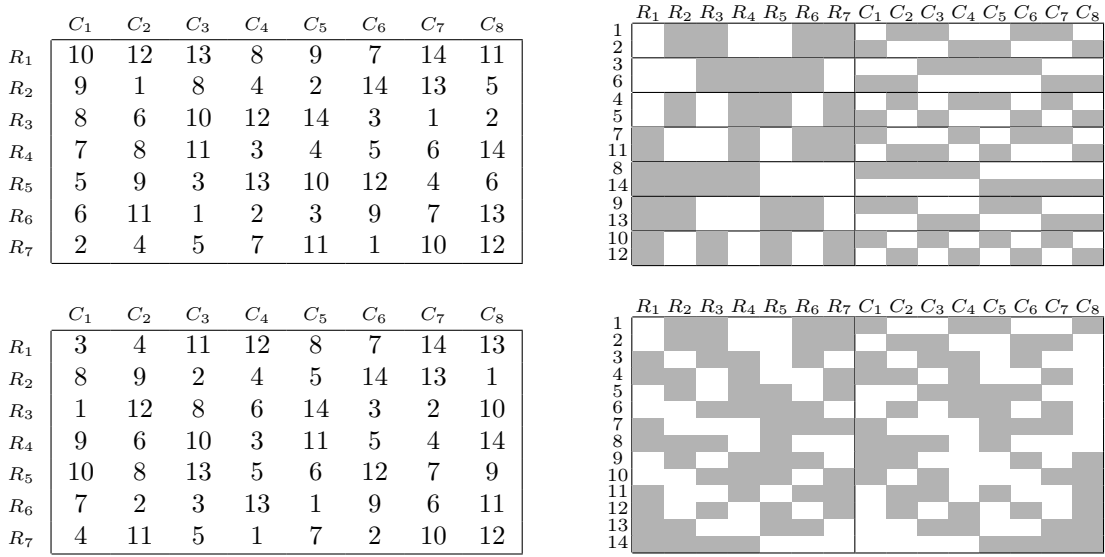

\begin{center}
\begin{tabular}{c c c}
\begin{tabular}{c|*{8}{c}|}
\multicolumn{1}{c}{} & \multicolumn{1}{c}{\scriptsize $C_{1}$} & \multicolumn{1}{c}{\scriptsize $C_{2}$} & \multicolumn{1}{c}{\scriptsize $C_{3}$} & \multicolumn{1}{c}{\scriptsize $C_{4}$} & \multicolumn{1}{c}{\scriptsize $C_{5}$} & \multicolumn{1}{c}{\scriptsize $C_{6}$} & \multicolumn{1}{c}{\scriptsize $C_{7}$} & \multicolumn{1}{c}{\scriptsize $C_{8}$} \\
\hhline{~--------}
{\scriptsize $R_{ 1}$} & $10$ & $12$ & $13$ & $ 8$ & $ 9$ & $ 7$ & $14$ & $11$ \\
{\scriptsize $R_{ 2}$} & $ 9$ & $ 1$ & $ 8$ & $ 4$ & $ 2$ & $14$ & $13$ & $ 5$ \\
{\scriptsize $R_{ 3}$} & $ 8$ & $ 6$ & $10$ & $12$ & $14$ & $ 3$ & $ 1$ & $ 2$ \\
{\scriptsize $R_{ 4}$} & $ 7$ & $ 8$ & $11$ & $ 3$ & $ 4$ & $ 5$ & $ 6$ & $14$ \\
{\scriptsize $R_{ 5}$} & $ 5$ & $ 9$ & $ 3$ & $13$ & $10$ & $12$ & $ 4$ & $ 6$ \\
{\scriptsize $R_{ 6}$} & $ 6$ & $11$ & $ 1$ & $ 2$ & $ 3$ & $ 9$ & $ 7$ & $13$ \\
{\scriptsize $R_{ 7}$} & $ 2$ & $ 4$ & $ 5$ & $ 7$ & $11$ & $ 1$ & $10$ & $12$ \\
\hhline{~--------}
\end{tabular}
& &
\renewcommand{\arraystretch}{0.5}
\setlength\tabcolsep{1pt}
\begin{tabular}{c|*{7}{c}|*{8}{c}|}
\multicolumn{1}{c}{} & \multicolumn{1}{c}{\scriptsize $R_{1}$} & \multicolumn{1}{c}{\scriptsize $R_{2}$} & \multicolumn{1}{c}{\scriptsize $R_{3}$} & \multicolumn{1}{c}{\scriptsize $R_{4}$} & \multicolumn{1}{c}{\scriptsize $R_{5}$} & \multicolumn{1}{c}{\scriptsize $R_{6}$} & \multicolumn{1}{c}{\scriptsize $R_{7}$} & \multicolumn{1}{c}{\scriptsize $C_{1}$} & \multicolumn{1}{c}{\scriptsize $C_{2}$} & \multicolumn{1}{c}{\scriptsize $C_{3}$} & \multicolumn{1}{c}{\scriptsize $C_{4}$} & \multicolumn{1}{c}{\scriptsize $C_{5}$} & \multicolumn{1}{c}{\scriptsize $C_{6}$} & \multicolumn{1}{c}{\scriptsize $C_{7}$} & \multicolumn{1}{c}{\scriptsize $C_{8}$} \\
\hhline{~---------------}
{\scriptsize  1} & & \cfl & \cfl & & & \cfl & \cfl & & \cfl & \cfl & & & \cfl & \cfl & \\
{\scriptsize  2} & & \cfl & \cfl & & & \cfl & \cfl & \cfl & & & \cfl & \cfl & & & \cfl \\
\hhline{~---------------}
{\scriptsize  3} & & & \cfl & \cfl & \cfl & \cfl & & & & \cfl & \cfl & \cfl & \cfl & & \\
{\scriptsize  6} & & & \cfl & \cfl & \cfl & \cfl & & \cfl & \cfl & & & & & \cfl & \cfl \\
\hhline{~---------------}
{\scriptsize  4} & & \cfl & & \cfl & \cfl & & \cfl & & \cfl & & \cfl & \cfl & & \cfl & \\
{\scriptsize  5} & & \cfl & & \cfl & \cfl & & \cfl & \cfl & & \cfl & & & \cfl & & \cfl \\
\hhline{~---------------}
{\scriptsize  7} & \cfl & & & \cfl & & \cfl & \cfl & \cfl & & & \cfl & & \cfl & \cfl & \\
{\scriptsize 11} & \cfl & & & \cfl & & \cfl & \cfl & & \cfl & \cfl & & \cfl & & & \cfl \\
\hhline{~---------------}
{\scriptsize  8} & \cfl & \cfl & \cfl & \cfl & & & & \cfl & \cfl & \cfl & \cfl & & & & \\
{\scriptsize 14} & \cfl & \cfl & \cfl & \cfl & & & & & & & & \cfl & \cfl & \cfl & \cfl \\
\hhline{~---------------}
{\scriptsize  9} & \cfl & \cfl & & & \cfl & \cfl & & \cfl & \cfl & & & \cfl & \cfl & & \\
{\scriptsize 13} & \cfl & \cfl & & & \cfl & \cfl & & & & \cfl & \cfl & & & \cfl & \cfl \\
\hhline{~---------------}
{\scriptsize 10} & \cfl & & \cfl & & \cfl & & \cfl & \cfl & & \cfl & & \cfl & & \cfl & \\
{\scriptsize 12} & \cfl & & \cfl & & \cfl & & \cfl & & \cfl & & \cfl & & \cfl & & \cfl \\
\hhline{~---------------}
\end{tabular}\\
\\
\begin{tabular}{c|*{8}{c}|}
\multicolumn{1}{c}{} & \multicolumn{1}{c}{\scriptsize $C_{1}$} & \multicolumn{1}{c}{\scriptsize $C_{2}$} & \multicolumn{1}{c}{\scriptsize $C_{3}$} & \multicolumn{1}{c}{\scriptsize $C_{4}$} & \multicolumn{1}{c}{\scriptsize $C_{5}$} & \multicolumn{1}{c}{\scriptsize $C_{6}$} & \multicolumn{1}{c}{\scriptsize $C_{7}$} & \multicolumn{1}{c}{\scriptsize $C_{8}$} \\
\hhline{~--------}
{\scriptsize $R_{ 1}$} & $ 3$ & $ 4$ & $11$ & $12$ & $ 8$ & $ 7$ & $14$ & $13$ \\
{\scriptsize $R_{ 2}$} & $ 8$ & $ 9$ & $ 2$ & $ 4$ & $ 5$ & $14$ & $13$ & $ 1$ \\
{\scriptsize $R_{ 3}$} & $ 1$ & $12$ & $ 8$ & $ 6$ & $14$ & $ 3$ & $ 2$ & $10$ \\
{\scriptsize $R_{ 4}$} & $ 9$ & $ 6$ & $10$ & $ 3$ & $11$ & $ 5$ & $ 4$ & $14$ \\
{\scriptsize $R_{ 5}$} & $10$ & $ 8$ & $13$ & $ 5$ & $ 6$ & $12$ & $ 7$ & $ 9$ \\
{\scriptsize $R_{ 6}$} & $ 7$ & $ 2$ & $ 3$ & $13$ & $ 1$ & $ 9$ & $ 6$ & $11$ \\
{\scriptsize $R_{ 7}$} & $ 4$ & $11$ & $ 5$ & $ 1$ & $ 7$ & $ 2$ & $10$ & $12$ \\
\hhline{~--------}
\end{tabular}
& &
\renewcommand{\arraystretch}{0.5}
\setlength\tabcolsep{1pt}
\begin{tabular}{c|*{7}{c}|*{8}{c}|}
\multicolumn{1}{c}{} & \multicolumn{1}{c}{\scriptsize $R_{1}$} & \multicolumn{1}{c}{\scriptsize $R_{2}$} & \multicolumn{1}{c}{\scriptsize $R_{3}$} & \multicolumn{1}{c}{\scriptsize $R_{4}$} & \multicolumn{1}{c}{\scriptsize $R_{5}$} & \multicolumn{1}{c}{\scriptsize $R_{6}$} & \multicolumn{1}{c}{\scriptsize $R_{7}$} & \multicolumn{1}{c}{\scriptsize $C_{1}$} & \multicolumn{1}{c}{\scriptsize $C_{2}$} & \multicolumn{1}{c}{\scriptsize $C_{3}$} & \multicolumn{1}{c}{\scriptsize $C_{4}$} & \multicolumn{1}{c}{\scriptsize $C_{5}$} & \multicolumn{1}{c}{\scriptsize $C_{6}$} & \multicolumn{1}{c}{\scriptsize $C_{7}$} & \multicolumn{1}{c}{\scriptsize $C_{8}$} \\
\hhline{~---------------}
{\scriptsize  1} & & \cfl & \cfl & & & \cfl & \cfl & \cfl & & & \cfl & \cfl & & & \cfl \\
{\scriptsize  2} & & \cfl & \cfl & & & \cfl & \cfl & & \cfl & \cfl & & & \cfl & \cfl & \\
{\scriptsize  3} & \cfl & & \cfl & \cfl & & \cfl & & \cfl & & \cfl & \cfl & & \cfl & & \\
{\scriptsize  4} & \cfl & \cfl & & \cfl & & & \cfl & \cfl & \cfl & & \cfl & & & \cfl & \\
{\scriptsize  5} & & \cfl & & \cfl & \cfl & & \cfl & & & \cfl & \cfl & \cfl & \cfl & & \\
{\scriptsize  6} & & & \cfl & \cfl & \cfl & \cfl & & & \cfl & & \cfl & \cfl & & \cfl & \\
{\scriptsize  7} & \cfl & & & & \cfl & \cfl & \cfl & \cfl & & & & \cfl & \cfl & \cfl & \\
{\scriptsize  8} & \cfl & \cfl & \cfl & & \cfl & & & \cfl & \cfl & \cfl & & \cfl & & & \\
{\scriptsize  9} & & \cfl & & \cfl & \cfl & \cfl & & \cfl & \cfl & & & & \cfl & & \cfl \\
{\scriptsize 10} & & & \cfl & \cfl & \cfl & & \cfl & \cfl & & \cfl & & & & \cfl & \cfl \\
{\scriptsize 11} & \cfl & & & \cfl & & \cfl & \cfl & & \cfl & \cfl & & \cfl & & & \cfl \\
{\scriptsize 12} & \cfl & & \cfl & & \cfl & & \cfl & & \cfl & & \cfl & & \cfl & & \cfl \\
{\scriptsize 13} & \cfl & \cfl & & & \cfl & \cfl & & & & \cfl & \cfl & & & \cfl & \cfl \\
{\scriptsize 14} & \cfl & \cfl & \cfl & \cfl & & & & & & & & \cfl & \cfl & \cfl & \cfl \\
\hhline{~---------------}
\end{tabular}
\end{tabular}

\end{center}
\caption{Two examples of $(7 \times 8, 14)$-triple arrays and representations of their underlying unordered triple arrays.
The top triple array is resolvable, while the bottom one is not.}
\label{fig:ex7-8-14}
\end{figure}

\begin{example}\label{ex:7814-res-or-not}
Not all (unordered) triple arrays are resolvable, even when the corresponding $\lrrc$ and $k$ are integers.
Indeed, for the parameter set $(7 \times 8, 14)$, we have $\lrrc = 2$ and $k = 2$, and two examples of $(7 \times 8, 14)$-triple arrays, one resolvable and one not, are given in Figure~\ref{fig:ex7-8-14}.
\end{example}

\subsection{Quad arrays}
\label{subsec:quad}

Resolvable (unordered) triple arrays turn out to have an interesting additional intersection property.
Returning to the notation from the \nameref{constr:res}, for all $i \neq s$ and $j$, we have
\[
|R_i \cap R_s \cap C_j| = \sum_{x\,:\,\{i,s\} \subseteq S_x} |\mathcal{B}_x \cap C_j| = \lrrc.
\]
This motivates the following definitions, where we note that the parameter $\lrrc$ can be directly calculated from the parameters $(r \times c,v)$ as in~\eqref{eq:lrrc-k}.

\begin{definition}\label{def:UQA}
An \textit{$(r \times c, v)$-unordered quad array} (UQA) is an $(r \times c, v)$-unordered triple array with row-sets $R_1,\dots, R_r$ and column-sets $C_1, \dots, C_c$, in which, for some integer $\lrrc$,
\[
|R_i \cap R_s \cap C_j| = \lrrc \text{ for all } i \neq s \text{ and } j.
\]

An \textit{$(r \times c, v)$-quad array} (QA) is an $(r \times c, v)$-triple array satisfying, for some integer $\lrrc$, an additional property:
\begin{enumerate}
    \customitem{(RRC)}\label{QA:rrc} any triple of two distinct rows and a column shares $\lrrc$ common symbols.
\end{enumerate}

A parameter set $(r \times c, v)$ is \textit{admissible for quad arrays} if it is admissible for triple arrays and $\lrrc$ is an integer.
\end{definition}

The name `quad array' continues the established naming convention: recalling Definition~\ref{def:TA}, a triple array is a binary equireplicate row-column design that satisfies three intersection properties~\ref{TA:rc},~\ref{TA:rr} and~\ref{TA:cc}.
Similarly, a binary equireplicate row-column design is called a \textit{mono array}~\cite[Definition~2.4]{jagerEnumerationConstructionRowColumn2025} if it satisfies~\ref{TA:cc}, and a \textit{double array}~\cite[Section~2]{mcsorleyDoubleArraysTriple2005a} if it satisfies~\ref{TA:rr} and~\ref{TA:cc}.

\begin{example}\label{ex:7814-quad-or-not}
Not all triple arrays are quad arrays, even when $\lrrc$ is an integer, as evidenced by the examples in Figure~\ref{fig:ex7-8-14}.
The top triple array is resolvable and is therefore a quad array.
On the other hand, in the bottom triple array, for instance, $|R_1 \cap R_2 \cap C_6| = |\{14\}| = 1$ while $|R_1 \cap R_2 \cap C_7| = |\{4, 13, 14\}| = 3$.
\end{example}

\begin{remark}
As noted above, any resolvable (unordered) triple array is an (unordered) quad array.
\end{remark}

On the other hand, we were not able to find an (unordered) quad array which is not resolvable.
This is somewhat surprising as the definition of a quad array does not require that $k = \frac{c}{e}$ is an integer, and even when this is the case, there is no obvious reason why every quad array should be resolvable.
For example, the parameter set $(16 \times 9, 24)$ is admissible for quad arrays but not for resolvable triple arrays.
However, in Observation~\ref{obs:quad-non-res} of Section~\ref{subsec:enum-extr} we will see that there are $(16 \times 9, 24)$-unordered triple arrays but no $(16 \times 9, 24)$-unordered quad arrays.

\subsection{Infinite families of resolvable unordered triple arrays}
\label{subsec:res-families}

In this section, we use the \nameref{constr:res} to produce three infinite families of resolvable unordered triple arrays, two extremal and one non-extremal.

\begin{theorem}\label{thm:AGnq->TA}
For any prime power $q$ and $n \geq 2$, there is a resolvable $(\frac{q^n - 1}{q - 1} \times q^n, \frac{q(q^n - 1)}{q - 1})$-unordered triple array.
\end{theorem}
\begin{proof}
The corresponding parameters are $e = q^{n - 1}$, $\lcc = \frac{q^{n - 1} - 1}{q - 1}$, $\lrr = q^{n - 1}(q - 1)$, $\lrrc = q^{n - 2}(q - 1)$, $k = q$.
When $n = 2$, let $\mathcal{S}$ be the trivial symmetric 2-$(q + 1, q, q - 1)$ design.
When $n > 2$, $\PG_{n - 2}(n - 1, q)$ is a symmetric 2-$(\frac{q^n - 1}{q - 1}, \frac{q^{n - 1} - 1}{q - 1}, \frac{q^{n - 2} - 1}{q - 1})$ design, so let $\mathcal{S}$ be its complement, that is, a symmetric 2-$(\frac{q^n - 1}{q - 1}, q^{n - 1}, q^{n - 2}(q - 1))$ design.
In any case, $\mathcal{S}$ has parameters $(r, e, \lrrc)$.
Let $\mathcal{B}$ be the unique resolution of the 2-$(c, e, \lcc)$ design $\AG_{n - 1}(n, q)$.
The theorem follows by applying the \nameref{constr:res} with these $\mathcal{S}$ and $\mathcal{B}$.
\end{proof}

We have already seen in Figure~\ref{fig:UTA-Res-346-4912} the first two unordered triple arrays from this family, corresponding to $q = 2$, $n = 2$ and $q = 3$, $n = 2$.

\begin{remark}\label{rm:aff->TA}
Similarly, \emph{any} affine plane of order $q$ gives rise to a resolvable $((q + 1) \times q^2, q(q + 1))$-unordered triple array: apply the \nameref{constr:res} with the trivial symmetric 2-$(q + 1, q, q - 1)$ design as $\mathcal{S}$ and the unique resolution of this affine plane as $\mathcal{B}$.
The current state of knowledge is that there are non-Galois affine planes for some prime power orders, but no known such planes for prime orders.
Also, there are no known affine planes for non-prime-power orders, but the construction presented here would in principle work for any finite affine plane.
\end{remark}

In Section~\ref{sec:AG->TA}, we further investigate parameters sets corresponding to Remark~\ref{rm:aff->TA} and the case $n = 2$ in Theorem~\ref{thm:AGnq->TA}.
In particular, in Theorem~\ref{thm:TA->AG} we show that all triple arrays on such parameter sets are resolvable.

\begin{theorem}\label{thm:HD->TA}
If there is a symmetric 2-$(4m - 1, 2m - 1, m - 1)$ design, then there is a resolvable $((4m - 1) \times 4m, 8m - 2)$-unordered triple array.
\end{theorem}
\begin{proof}
The corresponding parameters are $e = 2m$, $\lcc = 2m - 1$, $\lrr = 2m$, $\lrrc = m$, $k = 2$.
Let $\mathcal{D}$ be a symmetric 2-$(4m - 1, 2m - 1, m - 1)$ design on the symbol set $V$ with blocks $A_1, \dots, A_{4m - 1}$.
Introduce a new symbol $\sigma$ and let $\mathcal{D}'$ be the 2-$(c, e, \lrr)$ design on the symbol set $V \cup \{\sigma\}$ with blocks
\[
A_1 \cup \{\sigma\},\ V \setminus A_1,\ A_2 \cup \{\sigma\},\ V \setminus A_2,\ \dots,\ A_{4m - 1} \cup \{\sigma\},\ V \setminus A_{4m - 1}.
\]
Note that $\mathcal{D'}$ has a resolution $\mathcal{B}$ with parallel classes $\{A_i \cup \{\sigma\}, V \setminus A_i\}$.
Further, let $\mathcal{S}$ be the complement of $\mathcal{D}$, a symmetric 2-$(r, e, \lrrc)$ design.
The theorem follows by applying Construction~\ref{constr:res} with these $\mathcal{S}$ and $\mathcal{B}$.
\end{proof}

The symmetric designs used in Theorem~\ref{thm:HD->TA} belong to a well-known class of designs based on Hadamard matrices.
An \textit{Hadamard matrix of order $m$} is an $m \times m$ $(\pm 1)$-matrix with pairwise orthogonal rows.
For a general overview on Hadamard matrices as well as for further details on the facts listed below, see the recent book by Seberry and Yamada~\cite{seberryHadamardMatricesConstructions2020}.
It is known that a symmetric 2-$(4m - 1, 2m - 1, m - 1)$ design exists if and only if there exists an Hadamard matrix of order $4m$.
The long-standing \textit{Hadamard conjecture} states that Hadamard matrices of order $4m$ exist for all $m \geq 1$.
The smallest order for which no Hadamard matrix is known is $668$, corresponding to $m = 167$. This gives the following corollary to Theorem~\ref{thm:HD->TA}.

\begin{corollary}\label{cor:HD->TA-small}
For any $1 \leq m < 167$, there is a resolvable $((4m - 1) \times 4m, 8m - 2)$-unordered triple array.
\end{corollary}

Paley~\cite{paleyOrthogonalMatrices1933} constructed two infinite series of Hadamard matrices, one with orders $q + 1$ for prime powers $q \equiv 3 \Mod{4}$, and one with orders $2(q + 1)$ for prime powers $q \equiv 1 \Mod{4}$. This implies the following corollary.

\begin{corollary}\label{cor:HD->TA-Paley}
For any prime power $q \equiv 3 \Mod{4}$, there is a resolvable $(q \times (q + 1), 2q)$-unordered triple array.
For any prime power $q \equiv 1 \Mod{4}$, there is a resolvable $((2q + 1) \times (2q + 2), 4q + 2)$-unordered triple array.
\end{corollary}

In Theorem~\ref{thm:Paley} of Section~\ref{sec:paley}, we show that for any prime power $q \equiv 3 \Mod{4}$ there in fact exists a resolvable $(q \times (q + 1), 2q)$-triple array.

Curiously, Theorem~\ref{thm:AGnq->TA}, Remark~\ref{rm:aff->TA} and Theorem~\ref{thm:HD->TA} seem to cover most extremal parameter sets admissible for resolvable triple arrays, see Table~\ref{tbl:extr-res-comp}.
The smallest extremal parameter set not covered by them and for which the existence of unordered triple arrays is not ruled out by Theorem~\ref{thm:Agrawal-univ} and the non-existence of symmetric 2-$(r + c, r, \lcc)$ designs, is $(37 \times 112, 148)$.
It should be noted that there are no known examples of symmetric 2-$(149, 37, 9)$ designs, corresponding to \nameref{constr:Agrawal}, nor of resolvable 2-$(112, 28, 9)$ designs, corresponding to the \nameref{constr:res}.

\begin{table}[!ht]
\begin{center}
\begin{tabular}{|rr|rrrrrr|r|r|r|}
\hline
\multicolumn{8}{|c|}{Parameter set} & T\ref{thm:AGnq->TA} & R\ref{rm:aff->TA} & T\ref{thm:HD->TA}\\
\hline
$r \times c$ & $v$ & $e$ & $\lrc$ & $\lrr$ & $\lcc$ & $\lrrc$ & $k$ & $q, n$ & $q$ & $m$ \\
\hline
$3 \times 4$ & $6$ & $2$ & $2$ & $2$ & $1$ & $1$ & $2$ & $2, 2$ & $2$ & $1$\\
$4 \times 9$ & $12$ & $3$ & $3$ & $6$ & $1$ & $2$ & $3$ & $3, 2$ & $3$ & $-$\\
$5 \times 16$ & $20$ & $4$ & $4$ & $12$ & $1$ & $3$ & $4$ & $4, 2$ & $4$ & $-$\\
$6 \times 25$ & $30$ & $5$ & $5$ & $20$ & $1$ & $4$ & $5$ & $5, 2$ & $5$ & $-$\\
$7 \times 8$ & $14$ & $4$ & $4$ & $4$ & $3$ & $2$ & $2$ & $2, 3$ & $-$ & $2$\\
$8 \times 49$ & $56$ & $7$ & $7$ & $42$ & $1$ & $6$ & $7$ & $7, 2$ & $7$ & $-$\\
$9 \times 64$ & $72$ & $8$ & $8$ & $56$ & $1$ & $7$ & $8$ & $8, 2$ & $8$ & $-$\\
$10 \times 81$ & $90$ & $9$ & $9$ & $72$ & $1$ & $8$ & $9$ & $9, 2$ & $9$ & $-$\\
$11 \times 12$ & $22$ & $6$ & $6$ & $6$ & $5$ & $3$ & $2$ & $-$ & $-$ & $3$\\
$12 \times 121$ & $132$ & $11$ & $11$ & $110$ & $1$ & $10$ & $11$ & $11, 2$ & $11$ & $-$\\
$13 \times 27$ & $39$ & $9$ & $9$ & $18$ & $4$ & $6$ & $3$ & $3, 3$ & $-$ & $-$\\
$13 \times 144$ & $156$ & $12$ & $12$ & $132$ & $1$ & $11$ & $12$ & $-$ & $12$ & $-$\\
$14 \times 169$ & $182$ & $13$ & $13$ & $156$ & $1$ & $12$ & $13$ & $13, 2$ & $13$ & $-$\\
$15 \times 16$ & $30$ & $8$ & $8$ & $8$ & $7$ & $4$ & $2$ & $2, 4$ & $-$ & $4$\\
$16 \times 225$ & $240$ & $15$ & $15$ & $210$ & $1$ & $14$ & $15$ & $-$ & $15$ & $-$\\
$17 \times 256$ & $272$ & $16$ & $16$ & $240$ & $1$ & $15$ & $16$ & $16, 2$ & $16$ & $-$\\
$18 \times 289$ & $306$ & $17$ & $17$ & $272$ & $1$ & $16$ & $17$ & $17, 2$ & $17$ & $-$\\
$19 \times 20$ & $38$ & $10$ & $10$ & $10$ & $9$ & $5$ & $2$ & $-$ & $-$ & $5$\\
$19 \times 324$ & $342$ & $18$ & $18$ & $306$ & $1$ & $17$ & $18$ & $-$ & $18$ & $-$\\
$20 \times 361$ & $380$ & $19$ & $19$ & $342$ & $1$ & $18$ & $19$ & $19, 2$ & $19$ & $-$\\
\hline
\end{tabular}
\end{center}
\caption{The list of all extremal parameter sets with $r \leq 20$ and $r \leq c$ admissible for resolvable triple arrays, excluding those where there are no symmetric 2-$(r + c, r, \lcc)$ designs and thus, by Theorem~\ref{thm:Agrawal-univ}, no unordered triple arrays.
When the parameter set is covered by Theorem~\ref{thm:AGnq->TA}, Remark~\ref{rm:aff->TA} or Theorem~\ref{thm:HD->TA}, the corresponding column lists the relevant parameter values, otherwise, it contains the symbol $-$.}
\label{tbl:extr-res-comp}
\end{table}

So far, all the families of unordered triple arrays presented have been extremal, that is, with $v=r+c-1$.
Importantly, the \nameref{constr:res} also works for some non-extremal parameter sets, as proven in the following theorem.

\begin{theorem}\label{thm:PG3q}
For any prime power $q$, there is a $(\frac{q^3 - 1}{q - 1}, \frac{q^4 -1}{q - 1}, \frac{(q^4 - 1)(q^3 - 1)}{(q^2 - 1)(q - 1)})$-unordered triple array.
\end{theorem}
\begin{proof}
The corresponding parameters are $e = q + 1$,
$\lcc = 1$, $\lrr = q^2 + 1$, $\lrrc = 1$, $k = q^2 + 1$.
Denniston~\cite{DennistonPackings}, and later  Beutelspacher~\cite{BeutelspacherParallelisms} proved that $\PG(3, q)$  is resolvable for any prime power $q$.
The theorem follows by applying the \nameref{constr:res} with the symmetric 2-$(r, e, \lrrc)$ design $\mathcal{S} := \PG(2, q)$ and a resolution $\mathcal{B}$ of the 2-$(c, e, \lcc)$ design $\PG(3, q)$.
\end{proof}

For $q=2$, Theorem~\ref{thm:PG3q} gives $(7 \times 15, 35)$-unordered triple arrays, an example of which is given in Figure~\ref{fig:MPWY-71535}.
A more detailed investigation of this case is presented in Section~\ref{subsec:7x15}, where we enumerate and analyse the structure of all resolvable $(7 \times 15, 35)$-(unordered) triple arrays.
For $q=3$, the resulting unordered triple array has parameters $(13 \times 40, 130)$, for which there is no previously known example, as is also the case for larger $q$.

Table~\ref{tbl:non-extr-res} lists all non-extremal parameter sets with $r \leq 30$ admissible for resolvable triple arrays.
Besides the parameter sets covered by Theorem~\ref{thm:PG3q}, $(21 \times 15, 63)$ is the only non-extremal parameter set $(r \times c, v)$ with known examples of resolvable 2-$(c, e, \lcc)$ designs.
A number of resolutions of 2-$(15, 5, 6)$ designs have been constructed by Mathon and Rosa~\cite[Table 6]{mathon15FamilyBIBDs1989}, and in Section~\ref{subsec:21x15} we use these resolutions to find $(21 \times 15, 63)$-triple arrays, no examples of which were known before.

\begin{table}[!ht]
\begin{center}
\begin{tabular}{|rr|rrrrrr|r|}
\hline
\multicolumn{8}{|c|}{Parameter set} & T\ref{thm:PG3q}\\
\hline
$r \times c$ & $v$ & $e$ & $\lrc$ & $\lrr$ & $\lcc$ & $\lrrc$ & $k$ & $q$ \\
\hline
$7 \times 15$ & $35$ & $3$ & $3$ & $5$ & $1$ & $1$ & $5$ & $2$\\
$11 \times 45$ & $99$ & $5$ & $5$ & $18$ & $1$ & $2$ & $9$ & $-$\\
$13 \times 40$ & $130$ & $4$ & $4$ & $10$ & $1$ & $1$ & $10$ & $3$\\
$15 \times 91$ & $195$ & $7$ & $7$ & $39$ & $1$ & $3$ & $13$ & $-$\\
$19 \times 153$ & $323$ & $9$ & $9$ & $68$ & $1$ & $4$ & $17$ & $-$\\
$21 \times 15$ & $63$ & $5$ & $5$ & $3$ & $6$ & $1$ & $3$ & $-$\\
$21 \times 85$ & $357$ & $5$ & $5$ & $17$ & $1$ & $1$ & $17$ & $4$\\
$22 \times 133$ & $418$ & $7$ & $7$ & $38$ & $1$ & $2$ & $19$ & $-$\\
$23 \times 231$ & $483$ & $11$ & $11$ & $105$ & $1$ & $5$ & $21$ & $-$\\
$27 \times 325$ & $675$ & $13$ & $13$ & $150$ & $1$ & $6$ & $25$ & $-$\\
\hline
\end{tabular}
\end{center}
\caption{The list of all non-extremal parameter sets with $r \leq 30$ admissible for resolvable triple arrays.
When the parameter set is covered by Theorem~\ref{thm:PG3q}, the column T\ref{thm:PG3q} lists the relevant value of $q$, otherwise, it contains the symbol $-$.}
\label{tbl:non-extr-res}
\end{table}

\section{Resolvable Paley triple arrays}
\label{sec:paley}

As we have seen, the \nameref{constr:res} provides a novel way of solving the \nameref{pbm:exUTA} for a number of parameter sets.
In particular, it is the first general method of doing this for non-extremal parameter sets.
Nevertheless, similar to \nameref{constr:Agrawal}, it gives no indication of how to then solve the \nameref{pbm:UTAtoTA}.
There are two known infinite families of triple arrays~\cite{preecePaleyTripleArrays2005,nilsonTripleArraysDifference2017}, both of which are extremal.
They can each be seen as an implicit way of solving the \nameref{pbm:UTAtoTA} for certain parameter sets, though it should be mentioned that in both cases the triple array is constructed directly, without an intermediate step of building an unordered triple array.
In this section we investigate these two families and show that one of them contains an infinite number of resolvable triple arrays as a subfamily, while the other family contains no resolvable triple arrays.

Preece, Wallis and Yucas~\cite{preecePaleyTripleArrays2005} used the Hadamard matrices constructed by Paley in a manner somewhat different from Theorem~\ref{thm:HD->TA} to construct an infinite family of $(q \times (q + 1), 2q)$-triple arrays for all odd prime powers $q \geq 5$, which they called Paley triple arrays.
When $q \equiv 3 \Mod{4}$, their construction gives $6$ different types of arrays, and examining them closer, we can show that the first four of these types are resolvable, giving an infinite family of resolvable triple arrays.

To prepare for this result, we follow notation used in~\cite[p.241]{preecePaleyTripleArrays2005}, and let $\F_q = \{0 = w_1, w_2, \dots, w_q\}$, $\F_q' = \{0' = w_1', w_2', \dots, w_q'\}$ be two duplicate copies of the finite field with $q$ elements.
Let $Q$ denote the set of non-zero squares of elements of $\F_q$, let $N$ denote the set of non-squares, and let $Q_0 := Q \cup \{0\}$, $N_0 := N \cup \{0\}$.
For two non-zero elements $a, b \in F_q$, define the $q \times (q + 1)$ array $T$ as follows:
\[
T(i, j) = \begin{cases}
w_i - \frac{w_i - w_j}{a} & \text{if } j \leq q \text{ and } w_i - w_j \in Q,\\
(w_i + \frac{w_i - w_j}{b})' & \text{if } j \leq q \text{ and } w_i - w_j \in N_0,\\
w_i & \text{if } j = q + 1.
\end{cases}
\]
Assuming that $q \equiv 3 \Mod{4}$, $(a - 1)(b + 1) \in Q$, and that if $a - 1 \in N$ then $ab \in Q$, it is shown in~\cite[Theorem~10]{preecePaleyTripleArrays2005} that $T$ is a triple array.
We will assume that both $(a - 1)(b + 1) \in Q$ and $ab \in Q$. This corresponds to Paley triple arrays of types 1--4 out of the 6 types given in~\cite[proof of Theorem~10]{preecePaleyTripleArrays2005}, which, by~\cite[Corollary~12]{preecePaleyTripleArrays2005} all exist for every $q \geq 7$, $q \equiv 3 \Mod{4}$.

\begin{theorem}\label{thm:Paley}
For any prime power $q \geq 7$, $q \equiv 3 \Mod{4}$, Paley $(q \times (q + 1), 2q)$-triple arrays of types 1--4 are resolvable.
\end{theorem}
\begin{proof}
Partition $\F_q \cup \F_q'$ into $q$ pairs $\{w_s, w_s'\}$, $1 \leq s \leq q$.
To prove that a Paley triple array $T$ is resolvable, we will show that, for each $s$, the symbols $w_s$ and $w_s'$ appear in the same rows, and each column contains precisely one of $w_s, w_s'$.

Denote by $R_i$ and $C_j$ the sets of symbols appearing in row $i$ and column $j$ of $T$, respectively.
By definition, $R_i = ( w_i - \frac{Q_0}{a} ) \cup ( w_i + \frac{N_0}{b} )'$, and since $-1$ is a non-square in $\F_q$ when $q \equiv 3 \Mod{4}$, it holds that $( w_i - \frac{Q_0}{a} ) = ( w_i + \frac{N_0}{a} )$.
Since $ab \in Q$, either both $a$ and $b$ are square, or they are both non-square. 
We therefore have either $R_i = (w_i + N_0) \cup (w_i + N_0)'$ or $R_i = (w_i + Q_0) \cup (w_i + Q_0)'$.
In both cases, for each $s$, the symbols $w_s$ and $w_s'$ either both occur, or both do not occur in $R_i$.

Similarly, using $w_i - \frac{w_i - w_j}{a} = w_j + \frac{a - 1}{a}(w_i - w_j)$ and $(w_i + \frac{w_i - w_j}{b})' = (w_j + \frac{b + 1}{b}(w_i - w_j))'$, we have $C_j = ( w_j + \frac{a - 1}{a}Q ) \cup ( w_j + \frac{b + 1}{b}N_0 )'$.
Since $ab \in Q$ and $(a - 1)(b + 1) \in Q$, we have either $C_j = (w_j + Q) \cup (w_j + N_0)'$ or $C_j = (w_j + N) \cup (w_j + Q_0)'$.
In any case, $C_j$ contains exactly one of $w_s, w_s'$ for each $s$.
\end{proof}

Nilson and Cameron~\cite[Theorem 4.8]{nilsonTripleArraysDifference2017} constructed $((2u^2 - u) \times (2u^2 + u), 4u^2 - 1)$-triple arrays for all positive integers $u$ with the square-free part dividing 6.
When $u = 1$, the construction produces trivial $(1 \times 3, 3)$-triple arrays, and when $u > 1$, the corresponding parameter $e = u^2$ divides neither $r = 2u^2 - u$ nor $c = 2u^2 + u$, so these arrays are never resolvable.

\section{Enumeration and examples}
\label{sec:enum}

In this section, we first briefly describe the notions of equivalence we  used when enumerating various objects.
We then give an overview over the enumerative problems treated, how we have handled them algorithmically, and our results.

\subsection{Notions of equivalence}
\label{subseq:isomorph}

An \textit{isomorphism} between two block designs is a pair of bijective maps, one between their point sets and one between their block collections, which preserve the point-block inclusion relationship.
An \textit{isomorphism} between two resolutions of 2-designs is an isomorphism between these 2-designs which maps the parallel classes of one resolution to the parallel classes of the other.
An \textit{isomorphism} between two unordered triple arrays consists of three bijective maps, between their symbol sets, collections of row-sets and collections of column-sets, which preserve the inclusion relationships.
Two objects are \textit{isomorphic} if there exists an isomorphism between them.
An \textit{automorphism} is an isomorphism of the object with itself.
All automorphisms of $X$ form its \textit{(full) automorphism group} $\Aut X$.
Note that a 2-design may have multiple non-isomorphic resolutions.
The automorphism group of a resolution $\mathcal{R}$ of a 2-design $\mathcal{D}$ is a subgroup of the automorphism group of the 2-design: $\Aut \mathcal{R} \leq \Aut \mathcal{D}$.

Two triple arrays $T$ and $T'$ are \textit{isotopic} if there exist bijective maps $\varphi_R, \varphi_C, \varphi_V$ between their sets of rows, sets of columns, and symbol sets, respectively, such that if the cell $(i, j)$ of $T$ contains symbol $a$, then the cell $(\varphi_R(i), \varphi_C(j))$ of $T'$ contains symbol $\varphi_V(a)$.
The \textit{autotopism group} $\Aut T$ of a triple array $T$ consists of all its \textit{autotopisms}, that is, isotopisms of $T$ with itself.
If two triple arrays $T$ and $T'$ are isotopic, then their underlying unordered triple arrays $U_{T}$ and $U_{T'}$ are isomorphic.
On the other hand, an unordered triple array can, in general, be ordered in multiple non-isotopic ways, see Figure~\ref{fig:2ex-5-6-10} and the $(5 \times 6, 10)$ entry in Table~\ref{tbl:extr-UTA-counts} below.
For a triple array $T$, we have $\Aut T \leq \Aut U_T$.

\begin{figure}[!ht]
\begin{center}
\begin{tabular}{c c c}
\begin{tabular}{c|*{6}{c}|}
\multicolumn{1}{c}{} & \multicolumn{1}{c}{\scriptsize $C_{1}$} & \multicolumn{1}{c}{\scriptsize $C_{2}$} & \multicolumn{1}{c}{\scriptsize $C_{3}$} & \multicolumn{1}{c}{\scriptsize $C_{4}$} & \multicolumn{1}{c}{\scriptsize $C_{5}$} & \multicolumn{1}{c}{\scriptsize $C_{6}$} \\
\hhline{~------}
{\scriptsize $R_{ 1}$} & $ 1$ & $ 5$ & $ 9$ & $10$ & $ 4$ & $ 7$ \\
{\scriptsize $R_{ 2}$} & $ 4$ & $ 8$ & $ 3$ & $ 1$ & $10$ & $ 2$ \\
{\scriptsize $R_{ 3}$} & $ 6$ & $ 7$ & $ 2$ & $ 5$ & $ 3$ & $10$ \\
{\scriptsize $R_{ 4}$} & $ 7$ & $ 3$ & $ 1$ & $ 8$ & $ 6$ & $ 9$ \\
{\scriptsize $R_{ 5}$} & $ 2$ & $ 4$ & $ 5$ & $ 6$ & $ 9$ & $ 8$ \\
\hhline{~------}
\end{tabular}
& &
\begin{tabular}{c|*{6}{c}|}
\multicolumn{1}{c}{} & \multicolumn{1}{c}{\scriptsize $C_{1}$} & \multicolumn{1}{c}{\scriptsize $C_{2}$} & \multicolumn{1}{c}{\scriptsize $C_{3}$} & \multicolumn{1}{c}{\scriptsize $C_{4}$} & \multicolumn{1}{c}{\scriptsize $C_{5}$} & \multicolumn{1}{c}{\scriptsize $C_{6}$} \\
\hhline{~------}
{\scriptsize $R_{ 1}$} & $ 7$ & $ 4$ & $ 5$ & $ 1$ & $10$ & $ 9$ \\
{\scriptsize $R_{ 2}$} & $ 2$ & $ 3$ & $ 1$ & $10$ & $ 4$ & $ 8$ \\
{\scriptsize $R_{ 3}$} & $ 6$ & $ 7$ & $ 2$ & $ 5$ & $ 3$ & $10$ \\
{\scriptsize $R_{ 4}$} & $ 1$ & $ 8$ & $ 3$ & $ 6$ & $ 9$ & $ 7$ \\
{\scriptsize $R_{ 5}$} & $ 4$ & $ 5$ & $ 9$ & $ 8$ & $ 6$ & $ 2$ \\
\hhline{~------}
\end{tabular}\\
\\
\multicolumn{3}{c}{\renewcommand{\arraystretch}{0.5}
\setlength\tabcolsep{1pt}
\begin{tabular}{c|*{5}{c}|*{6}{c}|}
\multicolumn{1}{c}{} & \multicolumn{1}{c}{\scriptsize $R_{1}$} & \multicolumn{1}{c}{\scriptsize $R_{2}$} & \multicolumn{1}{c}{\scriptsize $R_{3}$} & \multicolumn{1}{c}{\scriptsize $R_{4}$} & \multicolumn{1}{c}{\scriptsize $R_{5}$} & \multicolumn{1}{c}{\scriptsize $C_{1}$} & \multicolumn{1}{c}{\scriptsize $C_{2}$} & \multicolumn{1}{c}{\scriptsize $C_{3}$} & \multicolumn{1}{c}{\scriptsize $C_{4}$} & \multicolumn{1}{c}{\scriptsize $C_{5}$} & \multicolumn{1}{c}{\scriptsize $C_{6}$} \\
\hhline{~-----------}
{\scriptsize  1} & \cfl & \cfl & & \cfl & & \cfl & & \cfl & \cfl & & \\
{\scriptsize  2} & & \cfl & \cfl & & \cfl & \cfl & & \cfl & & & \cfl \\
{\scriptsize  3} & & \cfl & \cfl & \cfl & & & \cfl & \cfl & & \cfl & \\
{\scriptsize  4} & \cfl & \cfl & & & \cfl & \cfl & \cfl & & & \cfl & \\
{\scriptsize  5} & \cfl & & \cfl & & \cfl & & \cfl & \cfl & \cfl & & \\
{\scriptsize  6} & & & \cfl & \cfl & \cfl & \cfl & & & \cfl & \cfl & \\
{\scriptsize  7} & \cfl & & \cfl & \cfl & & \cfl & \cfl & & & & \cfl \\
{\scriptsize  8} & & \cfl & & \cfl & \cfl & & \cfl & & \cfl & & \cfl \\
{\scriptsize  9} & \cfl & & & \cfl & \cfl & & & \cfl & & \cfl & \cfl \\
{\scriptsize 10} & \cfl & \cfl & \cfl & & & & & & \cfl & \cfl & \cfl \\
\hhline{~-----------}
\end{tabular}}
\end{tabular}

\end{center}
\caption{Two examples of non-isotopic $(5 \times 6, 10)$-triple arrays sharing the same underlying unordered triple array.}
\label{fig:2ex-5-6-10}
\end{figure}

We can therefore define an equivalence relation on triple arrays based on whether they have the same underlying unordered triple array. 
For some further context on this relation, an \emph{intercalate} in a triple array is a $2\times 2$ subsquare on two symbols, and flipping an intercalate means exchanging positions of the two symbols involved. This operation does not change the contents of any row-set or column-set, and so does not change the underlying unordered triple array. More generally, a \emph{trade} in a triple array is a set of cells such that there is a rearrangement of the symbols in these cells such that the row-sets and the column-sets do not change. If two triple arrays share the same underlying unordered triple array, they may be seen as both having trivial trades, carrying the one to the other by rearranging symbols while still leaving row-sets and column-sets intact.

\subsection{Computational methods}
\label{subseq:comput}

In this section we briefly describe how we approach the \nameref{pbm:exUTA} and the \nameref{pbm:UTAtoTA} computationally.
In all cases, we use the package \textrm{nauty} by McKay, see~\cite{mckayPracticalGraphIsomorphism2014}, to detect isomorphisms/isotopisms and calculate automorphism/autotopism group sizes.
All algorithms described below were implemented in \textrm{C++} and run on a personal laptop.
The total run time for all cases was around a hundred core hours.
The source code and all generated data are available at~\cite{gordeevZenodoLibrary25}.

\subsubsection*{Ordering unordered triple arrays}

In order to search for solutions to the \nameref{pbm:UTAtoTA}, we restate it as a special case of the following well-known NP-complete problem.

\begin{problem}[Exact Cover Problem]\label{pbm:exact-cover}
Given a set $X$ and a collection of its subsets $\mathcal{A}$, find an \textit{exact cover} of $X$, that is, a subcollection $\mathcal{A}^* \subseteq \mathcal{A}$ of pairwise disjoint subsets such that $\bigcup_{A \in \mathcal{A}^*} A = X$.
\end{problem}

For an $(r \times c, v)$-unordered triple array $U$ with row-sets $R_1, \dots, R_r$ and columns-sets $C_1, \dots, C_c$, the \nameref{pbm:UTAtoTA} is a special instance of the \nameref{pbm:exact-cover} with
\[
X := \{x_{ij}\,:\,1 \leq i \leq r, 1 \leq j \leq c\} \cup \{y_{ia}\,:\,1 \leq i \leq r, a \in R_i\} \cup \{z_{ja}\,:\, 1 \leq j \leq c, a \in C_j\}
\]
and
\[
\mathcal{A} := \{\{x_{ij}, y_{ia}, z_{ja}\}\,:\, 1 \leq i \leq r, 1 \leq j \leq c, a \in R_i \cap C_j\}.
\]

The \textit{Dancing links}~\cite{knuthDancingLinks2000} and \textit{Dancing cells}~\cite[Section ``Dancing cells'']{knuth2025fascicle7} backtracking algorithms provide a standard and quite efficient approach to the \nameref{pbm:exact-cover}.
The special case we are interested in has a number of properties which can be used to streamline and speed up the algorithms. For example, each element of $X$ belongs to at most $e$ subsets from the collection $\mathcal{A}$, where $e$, the replication number of $U$, is quite small compared to the total size of $\mathcal{A}$. Similarly, each set $A \in \mathcal{A}$ has constant size 3, much smaller than the total size of $X$.
To fully utilize these properties and maximize the search efficiency, we have implemented our own specialized version of the Dancing cells algorithm.

For a given unordered triple array $U$, we can find all non-isotopic triple arrays $T$ with $U_T = U$ as follows: first, find all solutions to the corresponding instance of the \nameref{pbm:exact-cover} using the Dancing cells algorithm, and then remove isotopic copies from the resulting list of triple arrays.

To verify the correctness of our computations, we perform a consistency check based on the orbit-stabilizer theorem: for each found triple array $T$, one must have $|\Aut T| \cdot N(T) = |\Aut U|$, where $N(T)$ is the total number of found triple arrays isotopic to $T$.
We verify that this identity holds for each triple array $T$ obtained during the search.

In some cases where enumerating all solutions to the \nameref{pbm:exact-cover} is not feasible, usually due to space constraints, we use the Dancing cells algorithm to only look for examples of solutions. 

\subsubsection*{Enumerating extremal unordered triple arrays}

If two $(r \times c, r + c - 1)$-unordered triple arrays are isomorphic, they correspond to isomorphic symmetric 2-$(r + c, r, \lcc)$ designs in \nameref{constr:Agrawal}.
Therefore, if the complete list of non-isomorphic symmetric 2-$(r + c, r, \lcc)$ designs is known, then, by Theorem~\ref{thm:Agrawal-univ}, we can generate all non-isomorphic $(r \times c, r + c - 1)$-unordered triple arrays as follows: for each symmetric 2-$(r + c, r, \lcc)$ design $\mathcal{S}$ on the point set $V$ and each point $\sigma \in V$, apply \nameref{constr:Agrawal} to obtain the unordered triple array $U_{\mathcal{S},\sigma}$, and then remove isomorphic copies from the resulting list of unordered triple arrays.

To verify the correctness of our computations, we again perform a consistency check based on the orbit-stabilizer theorem: for each $U_{\mathcal{S},\sigma}$ we verify that $|\Aut U_{\mathcal{S},\sigma}| \cdot N(U_{\mathcal{S},\sigma}) = |\Aut \mathcal{S}|$, where $N(U_{\mathcal{S},\sigma})$ is the total number of constructed unordered triple arrays isomorphic to $U_{\mathcal{S},\sigma}$.

\subsubsection*{Enumerating resolvable unordered triple arrays}

To generate all non-isomorphic resolvable $(r \times c, v)$-unordered triple arrays corresponding to a given symmetric 2-$(r, e, \lrrc)$ design $\mathcal{S}$ with blocks $S_1, \dots, S_r$ and a given resolution $\mathcal{Z}$ of a 2-$(c, e, \lcc)$ design with parallel classes $\mathcal{Z}_1, \dots, \mathcal{Z}_r$, for each permutation $\pi \in \mathfrak{S}_r$, where $\mathfrak{S}_r$ is the symmetric group on $r$ elements, we apply the \nameref{constr:res} with parallel classes labeled $\mathcal{B}_i := \mathcal{Z}_{\pi(i)}$ to obtain the unordered triple array $U_{\pi}$, and then remove isomorphic copies from the resulting list of unordered triple arrays.
Note that the choice of $\pi$ was implicit in the description of the \nameref{constr:res} as the parallel classes were already assumed to be labeled, but different permutations $\pi$ may lead to non-isomorphic unordered triple arrays.

Similar to previous cases, we verify the correctness of our computations by performing a consistency check based on the orbit-stabilizer theorem: for each $U_{\pi}$ we verify that $|\Aut U_{\pi}| \cdot N(U_{\pi}) = |\Aut \mathcal{S}| \cdot |\Aut \mathcal{R} |$, where $N(U_{\pi})$ is the total number of constructed unordered triple arrays isomorphic to $U_{\pi}$.

\subsection{The extremal case}
\label{subsec:enum-extr}

Heinlein, Ivanov, McKay and {\"O}sterg\r{a}rd~\cite{heinlein2023} provide complete lists of non-isomorphic symmetric 2-designs for a number of parameters.
Using this data, we employ \nameref{constr:Agrawal} to generate all corresponding unordered triple arrays.
We excluded symmetric 2-$(31, 15, 7)$ designs from consideration, as their total number exceeds 10 billion, and only some are given in~\cite{heinlein2023}.
On the other hand, to cover all known symmetric 2-designs with $\lambda = 2$, also known as \textit{biplanes}, we applied \nameref{constr:Agrawal} to the two symmetric 2-$(79, 13, 2)$ designs produced by SageMath~\cite{sagemath25}, which were first constructed by Aschbacher~\cite{aschbacherCollineationGroupsSymmetric1971}.
To the best of our knowledge, full enumerations of symmetric 2-designs are not known for any parameters not present in~\cite{heinlein2023}.

For some parameter sets, we use the obtained complete lists of non-isomorphic unordered triple arrays to generate all non-isotopic triple arrays.
This gives a further correctness check, as we verify that the number of found triple arrays matches our previous work \cite[Table C.1]{gordeevTripleArrays2026}, where a completely different algorithm was used.
For larger parameter sets, the enumeration of all non-isotopic triple arrays turned out to be unattainable with current methods, mainly due to space limitations.
In each such case, we empirically investigate Conjecture~\ref{conj:Agrawal-gen} by finding at least one triple array for each unordered triple array.
The summary of results is given in Table~\ref{tbl:extr-UTA-counts}.
For more details, see Appendix~\ref{ap:extr}.

\begin{table}[!ht]
\begin{center}
\begin{tabular}{|rr|rrrr|rr|r|r|}
\hline
\multicolumn{6}{|c|}{Parameter set} & \multicolumn{2}{c|}{Symmetric 2-design} & \multicolumn{1}{c|}{UTA} & \multicolumn{1}{c|}{TA} \\
\hline
$r \times c$ & $v$ & $e$ & $\lrc$ & $\lrr$ & $\lcc$ & $(r + c, r, \lcc)$ & Total \# & Total \# & Total \# \\
\hline
$3 \times 4$ & $6$ & $2$ & $2$ & $2$ & $1$ & $(7, 3, 1)$ & $1$ & $1$ & $0$\\
$4 \times 9$ & $12$ & $3$ & $3$ & $6$ & $1$ & $(13, 4, 1)$ & $1$ & $1$ & $1$\\
$5 \times 6$ & $10$ & $3$ & $3$ & $3$ & $2$ & $(11, 5, 2)$ & $1$ & $1$ & $7$\\
$5 \times 16$ & $20$ & $4$ & $4$ & $12$ & $1$ & $(21, 5, 1)$ & $1$ & $1$ & $26804$\\
$6 \times 10$ & $15$ & $4$ & $4$ & $6$ & $2$ & $(16, 6, 2)$ & $3$ & $3$ & $270119$\\
$6 \times 25$ & $30$ & $5$ & $5$ & $20$ & $1$ & $(31, 6, 1)$ & $1$ & $1$ & $+$\\
$7 \times 8$ & $14$ & $4$ & $4$ & $4$ & $3$ & $(15, 7, 3)$ & $5$ & $10$ & $684782$\\
$8 \times 49$ & $56$ & $7$ & $7$ & $42$ & $1$ & $(57, 8, 1)$ & $1$ & $1$ & $+$\\
$9 \times 10$ & $18$ & $5$ & $5$ & $5$ & $4$ & $(19, 9, 4)$ & $6$ & $22$ & $+$\\
$9 \times 16$ & $24$ & $6$ & $6$ & $10$ & $3$ & $(25, 9, 3)$ & $78$ & $1382$ & $+$\\
$9 \times 28$ & $36$ & $7$ & $7$ & $21$ & $2$ & $(37, 9, 2)$ & $4$ & $8$ & $+$\\
$9 \times 64$ & $72$ & $8$ & $8$ & $56$ & $1$ & $(73, 9, 1)$ & $1$ & $1$ & $+$\\
$10 \times 21$ & $30$ & $7$ & $7$ & $14$ & $3$ & $(31, 10, 3)$ & $151$ & $3809$ & $+$\\
$10 \times 81$ & $90$ & $9$ & $9$ & $72$ & $1$ & $(91, 10, 1)$ & $4$ & $7$ & $+$\\
$11 \times 12$ & $22$ & $6$ & $6$ & $6$ & $5$ & $(23, 11, 5)$ & $1106$ & $23360$ & $+$\\
$11 \times 45$ & $55$ & $9$ & $9$ & $36$ & $2$ & $(56, 11, 2)$ & $5$ & $16$ & $+$\\
$13 \times 14$ & $26$ & $7$ & $7$ & $7$ & $6$ & $(27, 13, 6)$ & $208310$ & $5606594$ & $+$\\
$13 \times 66$ & $78$ & $11$ & $11$ & $55$ & $2$ & $(79, 13, 2)$ & $\geq 2$ & $\geq 8$ & $+$\\
\hline
\end{tabular}
\end{center}
\caption{The number of non-isomorphic unordered triple arrays (UTA) and non-isotopic triple arrays (TA) for all extremal parameter sets where the complete list of corresponding symmetric 2-designs is available at~\cite{heinlein2023}, and for the parameter set $(13 \times 66, 78)$. The symbol $+$ indicates that for each unordered triple array $U$ we found at least one triple array $T$ with $U_T = U$.}
\label{tbl:extr-UTA-counts}
\end{table}

\begin{observation}
    Every unordered triple array with $e > 2$ in Table~\ref{tbl:extr-UTA-counts} corresponds to at least one triple array, which corroborates Conjecture~\ref{conj:Agrawal-gen}.
\end{observation}\label{obs:UTAgivesTA}

Having generated all unordered triple arrays, and in some cases all triple arrays, we can also easily check which of them are resolvable and which of them are (unordered) quad arrays.
The summary is given in Table~\ref{tbl:extr-QTA-counts}. 

\begin{table}[!ht]
\begin{center}
\begin{tabular}{|rr|rr|rrr|rrr|}
\hline
$r \times c$ & $v$ & $\lrrc$ & $k$ & \multicolumn{1}{c}{UTA} & \multicolumn{1}{c}{UQA} & \multicolumn{1}{c|}{RUTA} & \multicolumn{1}{c}{TA} & \multicolumn{1}{c}{QA} & \multicolumn{1}{c|}{RTA} \\
\hline
$3 \times 4$ & $6$ & $1$ & $2$ & $1$ & $1$ & $1$ & $0$ & $0$ & $0$\\
$4 \times 9$ & $12$ & $2$ & $3$ & $1$ & $1$ & $1$ & $1$ & $1$ & $1$\\
$5 \times 16$ & $20$ & $3$ & $4$ & $1$ & $1$ & $1$ & $26804$ & $26804$ & $26804$\\
$6 \times 25$ & $30$ & $4$ & $5$ & $1$ & $1$ & $1$ & $+$ & $+$ & $+$\\
$7 \times 8$ & $14$ & $2$ & $2$ & $10$ & $4$ & $4$ & $684782$ & $43388$ & $43388$\\
$8 \times 49$ & $56$ & $6$ & $7$ & $1$ & $1$ & $1$ & $+$ & $+$ & $+$\\
$9 \times 64$ & $72$ & $7$ & $8$ & $1$ & $1$ & $1$ & $+$ & $+$ & $+$\\
$10 \times 81$ & $90$ & $8$ & $9$ & $7$ & $7$ & $7$ & $+$ & $+$ & $+$\\
$11 \times 12$ & $22$ & $3$ & $2$ & $23360$ & $20$ & $20$ & $+$ & $+$ & $+$\\
$16 \times 9$ & $24$ & $2$ & $3/2$ & $1382$ & $0$ & $-$ & $+$ & $0$ & $-$\\
\hline
\end{tabular}
\end{center}
\caption{The number of non-isomorphic unordered triple arrays (UTA), unordered quad arrays (UQA) and resolvable unordered triple arrays (RUTA), non-isotopic triple arrays (TA), quad arrays (QA) and resolvable triple arrays (RTA).
The parameter sets covered here are those from Table~\ref{tbl:extr-UTA-counts} with integer $\lrrc$. Note that $(16 \times 9, 24)$ is transposed.
The symbol $+$ indicates that we found examples of corresponding objects.
The symbol $-$ indicates that the parameter set is not admissible for resolvable triple arrays.}
\label{tbl:extr-QTA-counts}
\end{table}

\begin{observation}\label{obs:QA-vs-res}
    Every (unordered) quad array found in Table~\ref{tbl:extr-QTA-counts} is resolvable.
\end{observation}

\begin{observation}\label{obs:quad-non-res}
There are $(7 \times 8, 14)$-triple arrays and $(11 \times 12, 22)$-triple arrays that are neither resolvable nor quad arrays, despite these parameters being admissible for resolvable triple arrays.

None of the $(16 \times 9, 24)$-triple arrays are quad arrays, despite the parameter set being admissible for quad arrays.
\end{observation}

\begin{observation}\label{obs:7814-res}
All 4 resolvable $(7 \times 8, 14)$-unordered triple arrays have the same component designs: the column design is the unique resolvable 2-$(8, 4, 3)$ design with unique resolution, and the row design is the 2-multiple of the symmetric 2-$(7, 4, 2)$ design, the complement of the Fano plane.
Peculiarly, they come from pairwise non-isomorphic symmetric 2-designs in \nameref{constr:Agrawal}, see Table~\ref{tbl:7814}.

Similarly, all 20 resolvable $(11 \times 12, 22)$-unordered triple arrays have the same component designs: the column design is the unique resolvable 2-$(12, 6, 5)$ design with unique resolution, and the row design is the 2-multiple of the unique symmetric 2-$(11, 6, 3)$ design.
Again, they come from pairwise non-isomorphic symmetric 2-designs in \nameref{constr:Agrawal}.
\end{observation}

\subsection{Kirkman parades and resolvable \texorpdfstring{$(7 \times 15, 35)$}{(7 x 15, 35)}-triple arrays}
\label{subsec:7x15}

Among the admissible non-extremal parameter sets, $(7 \times 15, 35)$ is the smallest one by any measure (number of rows, number of columns, number of symbols), and, as already mentioned, the only parameter set for which a non-extremal triple array has been found so far, see Figure~\ref{fig:MPWY-71535}.
It is admissible for resolvable triple arrays, with $\lrrc = 1$ and $k = 5$.

To generate all resolvable $(7 \times 15, 35)$-unordered triple arrays, we apply the \nameref{constr:res} for all possible combinations of a symmetric 2-$(7, 3, 1)$ design $\mathcal{S}$ and a resolution $\mathcal{B}$ of a 2-$(15, 3, 1)$ design.
Up to isomorphism, there is only one choice for $\mathcal{S}$, the Fano plane $\PG(2, 2)$.
The resolutions of 2-$(15, 3, 1)$ designs are known as \emph{Kirkman parades}, named so by Cole~\cite{cole_kirkman_parades}.
For a more contemporary treatment on Kirkman parades, see the recent work of Pavone~\cite{pavoneSevenNonisomorphicSolutions2023}.
There are 7 non-isomorphic Kirkman parades.
Following~\cite[Table 1]{pavoneSevenNonisomorphicSolutions2023}, we will refer to them as `1a', `1b', `7a', `7b', `19a', `19b' and `61'.
As the names suggest, they stem from four non-isomorphic resolvable 2-$(15,3,1)$ designs, three of which correspond to two Kirkman parades each, and one to a single Kirkman parade `61'.
The 7 Kirkman parades give rise to a total of 42 $(7 \times 15, 35)$-unordered triple arrays and 85 $(7 \times 15, 35)$-triple arrays, see Table~\ref{tbl:71535} for a summary.
An extended Table~\ref{tbl:extended-7-15-35} is available in Appendix~\ref{ap:7-15-35}.

\begin{table}[ht]
\begin{center}
\begin{tabular}{|rr|*{8}{r|}}
\hline
\multicolumn{2}{|c|}{Parade} & 1a & 1b & 7a & 7b & 19a & 19b & 61 & All \\
\hline
\multicolumn{10}{c}{\multirow{2}{*}{Unordered triple arrays}} \\
\multicolumn{10}{c}{} \\
\hline
\multicolumn{2}{|c|}{Total \#} & $4$ & $4$ & $6$ & $6$ & $8$ & $8$ & $6$ & $42$ \\
\hline
$|\Aut|$ & $1$ &  &  &  &  & $1$ & $1$ &  & $2$ \\
 & $3$ &  &  & $2$ & $2$ & $3$ & $3$ & $4$ & $14$ \\
 & $4$ &  &  & $2$ & $2$ & $1$ & $1$ &  & $6$ \\
 & $12$ & $1$ & $1$ &  &  & $3$ & $3$ &  & $8$ \\
 & $21$ & $1$ & $1$ &  &  &  &  & $2$ & $4$ \\
 & $24$ & $1$ & $1$ & $2$ & $2$ &  &  &  & $6$ \\
 & $168$ & $1$ & $1$ &  &  &  &  &  & $2$ \\
\hline
\multicolumn{10}{c}{\multirow{2}{*}{Triple arrays}} \\
\multicolumn{10}{c}{} \\
\hline
\multicolumn{2}{|c|}{Total \#} & $0$ & $3$ & $24$ & $4$ & $21$ & $21$ & $12$ & $85$ \\
\hline
$|\Aut|$ & $1$ &  & $3$ & $12$ & $4$ & $21$ & $21$ & $12$ & $73$ \\
 & $3$ &  &  & $12$ &  &  &  &  & $12$ \\
\hline
\end{tabular}
\end{center}
\caption{The total number of non-isomorphic resolvable $(7 \times 15, 35)$-unordered triple arrays and non-isotopic resolvable $(7 \times 15, 35)$-triple arrays sorted by automorphism/autotopism group order, as well as the number of arrays corresponding to each of the Kirkman parades.
The naming convention for the Kirkman parades is borrowed from~\cite[Table 1]{pavoneSevenNonisomorphicSolutions2023}.}
\label{tbl:71535}
\end{table}

\begin{observation}\label{obs:71535-PG}
Kirkman parades `1a' and `1b' are the two non-isomorphic resolutions of $\PG(3, 2)$, corresponding to the case $q = 2$ in Theorem~\ref{thm:PG3q}.
Kirkman parade `1b' corresponds to the only previously known $(7 \times 15, 35)$-triple array from Figure~\ref{fig:MPWY-71535}.
This parade gives rise to 4 non-isomorphic $(7 \times 15, 35)$-unordered triple arrays.
Two of these cannot be ordered, the third has two non-isotopic orderings, one of which is the triple array from Figure~\ref{fig:MPWY-71535}, and the fourth can be ordered uniquely up to isotopism, see Table~\ref{tbl:extended-7-15-35} and Figure~\ref{fig:PG-71535}.    
\end{observation}

\begin{observation}\label{obs:71535-no-sol}
Kirkman parade `1a', one of two resolutions of $\PG(3, 2)$, corresponds to 4 non-isomorphic resolvable $(7 \times 15, 35)$-unordered triple arrays, none of which can be ordered.
Moreover, in total, 12 out of 42 resolvable $(7 \times 15, 35)$-unordered triple arrays cannot be ordered, with each Kirkman parade except `7a' giving rise to at least one these 12, see Table~\ref{tbl:extended-7-15-35}.
This is in some sense in contrast with Conjectures~\ref{conj:Agrawal} and~\ref{conj:Agrawal-gen} and Observation~\ref{obs:UTAgivesTA}.
\end{observation}

It should be noted that in the presentation given by Yucas~\cite{yucas_structure_7x15} on the structure of the $(7 \times 15, 35)$-triple array from Figure~\ref{fig:MPWY-71535}, the points in the blocks of the resolution had been ordered in such a way that the \nameref{pbm:UTAtoTA} was in essence already solved, but no indication was given as to how to find this ordering.
Yucas remarked that there was ``nothing special'' about the choice of $\PG(3, 2)$ as the column design used in the construction of this triple array.
In view of Observations~\ref{obs:71535-PG} and~\ref{obs:71535-no-sol}, this remark is not quite accurate.

\subsection{Resolvable \texorpdfstring{$(21 \times 15, 63)$}{(21 x 15, 63)}-triple arrays}
\label{subsec:21x15}

The parameter set $(21 \times 15, 63)$ is also admissible for resolvable triple arrays, with $\lrrc = 1$ and $k = 3$.
To produce a resolvable $(21 \times 15, 63)$-unordered triple array, the \nameref{constr:res} requires a symmetric 2-$(21, 5, 1)$ design $\mathcal{S}$ and a resolution $\mathcal{B}$ of a 2-$(15, 5, 6)$ design.
Up to isomorphism, the only choice for $\mathcal{S}$ is the Galois projective plane $\PG(2, 4)$.
The resolutions of 2-$(15, 5, 6)$ designs have not been fully classified, but 149 non-isomorphic resolutions on these parameters have been constructed by Mathon and Rosa~\cite[Table~6]{mathon15FamilyBIBDs1989}.

For each one of these 149 resolutions, we were able to produce a $(21 \times 15, 63)$-triple array by first applying the \nameref{constr:res} with parallel classes of the resolution labeled in random order, and then finding some ordering of the produced unordered triple array using the Dancing cells algorithm.
An example is given in Figure~\ref{fig:ex21-15-63}.
We emphasize that these are the first known triple arrays on this parameter set, and that this is just the second non-extremal parameter set on which triple arrays have been constructed.

We did not attempt to fully enumerate resolvable triple arrays corresponding to resolutions from~\cite{mathon15FamilyBIBDs1989}.
Empirically, the number of non-isomorphic unordered triple arrays corresponding to each given resolution as well as the number of non-isotopic triple arrays corresponding to each produced unordered triple array seem to be very large, so full enumeration via the current approach does not seem feasible due to both space and time limitations.
Moreover, the currently known list of resolutions is likely not exhaustive, so neither would be the lists of corresponding (unordered) triple arrays.

\subsection{Other non-extremal parameter sets}
\label{subsec:enum-non-extr}

As noted in Section~\ref{subsec:res-families}, all other non-extremal parameter sets $(r \times c, v)$ with known examples of resolvable 2-$(c, e, \lcc)$ designs are covered by Theorem~\ref{thm:PG3q}.
The smallest such parameter set besides $(7 \times 15, 35)$ is $(13 \times 40, 130)$, corresponding to $q = 3$.
We attempted to construct $(13 \times 40, 130)$-triple arrays by the same approach as in the previous section.
We obtained examples of $(13 \times 40, 130)$-unordered triple arrays using the \nameref{constr:res} with the symmetric 2-$(13, 4, 1)$ design $\mathcal{S} := \PG(2, 3)$ and two different resolutions $\mathcal{B}$ of 2-$(40, 4, 1)$ designs with randomly labeled blocks.
Both resolutions were produced by Sagemath~\cite{sagemath25}, the first is a resolution of $\PG(3, 3)$, while the second comes from~\cite[Theorem VII.7.4]{bethDesignTheory1999}.

Unfortunately, we could not find an ordering of any such unordered triple array, despite the search running, in some cases, for multiple core-days before being aborted.
This does not necessarily imply that such an ordering does not exist, as we deemed that an exhaustive search for orderings of any given $(13 \times 40, 130)$-unordered triple array would take an unreasonable amount of time.
It should be noted, however, that in all previously considered cases where the search was successful, at least some of the orderings were found almost immediately.

The complete list of 73343 non-isomorphic resolutions of $\PG(3, 3)$, enumerated by Betten~\cite{BettenPackings}, may prove useful in the search for $(13 \times 40, 130)$-triple arrays.
In view of the above though, it seems likely that some further theoretical guidance would be required for a computer search to successfully construct such a triple array.

\section{Triple arrays and finite affine planes}
\label{sec:AG->TA}

In this section, we further investigate extremal (unordered) triple arrays on the parameter set $((q + 1) \times q^2, q(q + 1))$, which we have already seen in the case $n = 2$ of Theorem~\ref{thm:AGnq->TA} and in Remark~\ref{rm:aff->TA}.
We show that all such (unordered) triple arrays are resolvable, with the corresponding resolvable 2-design being an affine plane of order $q$, and propose two problems equivalent to the \nameref{pbm:UTAtoTA} for these parameter sets.

\subsection{\texorpdfstring{$((q + 1) \times q^2, q(q + 1))$}{((q+1)xq2,q(q+1))}-unordered triple arrays}

Recall that Youden rectangles are triple arrays with $v = c$, or, equivalently, $r = e$.
Triple arrays with $r = e + 1$ are, in a sense, closest to Youden rectangles among all non-trivial triple arrays.
One can thus expect them to have simpler structure than triple arrays in general.
This is supported by the following theorem of J{\"a}ger, Markstr{\"o}m, Shcherbak and the second author, see~\cite[Theorem~3.2 and Remark~3.3]{jagerEnumerationConstructionRowColumn2025}.

\begin{theorem}
Any binary equireplicate $r \times c$ row-column design with replication number $e$ and $r = e + 1$ satisfies property~\ref{TA:rr} of triple arrays.
\end{theorem}

The next lemma shows that (unordered) triple arrays with $r = e + 1$ are precisely $((q + 1) \times q^2, q(q + 1))$-(unordered) triple arrays with $q = e$.

\begin{lemma}\label{lm:r=e+1}
Let $(r \times c, v)$ be a non-trivial admissible parameter set with $r = e + 1$.
Then $c = e^2$ and $v = e(e + 1)$.
\end{lemma}
\begin{proof}
Recall that $v = \frac{rc}{e} = \frac{(e + 1)c}{e}$, so, due to $e$ and $e + 1$ being coprime, $c$ must be divisible by $e$.
Further, $\lcc = \frac{r(e - 1)}{c - 1} \equiv \frac{1(-1)}{-1} \equiv 1 \pmod{e}$.
A non-trivial parameter set has $e < c$, thus $\lcc = \frac{r(e - 1)}{c - 1} < r = e + 1$, i.e. $\lcc = 1$.
This, in turn, implies $c = \frac{r(e - 1)}{\lcc} + 1 = e^2$ and $v = \frac{rc}{e} = e(e + 1)$.
\end{proof}

It turns out that all $((q + 1) \times q^2, q(q + 1)$-(unordered) triple arrays are resolvable.

\begin{theorem}\label{thm:TA->AG}
Any $((q + 1) \times q^2, q(q + 1))$-unordered triple array is resolvable.
\end{theorem}
\begin{proof}
Let $U$ be a $((q + 1) \times q^2, q(q + 1))$-unordered triple array.
Its parameters are $e = \lrc = q$, $\lcc = 1$, $\lrr = q(q - 1)$.
Let $R_1, \dots, R_r$ denote the row-sets and $C_1, \dots, C_c$ the column-sets of $U$.
Further, let $V$ be the symbol set of $U$, and define $V_i := V \setminus R_i$, that is, $V_i$ is the set of symbols missing from row $i$.
Each symbol appears in all rows but one, so the sets $V_1, \dots, V_r$ form a partition of $V$ with $|V_i|=q$.

The row design $RD_U$ is a 2-$(q + 1, q, q(q - 1))$ design, that is, the $q$-multiple of the trivial symmetric 2-$(q + 1, q, q - 1)$ design.
The column design $CD_U$ is a 2-$(q^2, q, 1)$ design, that is, a finite affine plane of order $q$.
The symbols of each set $V_x$ correspond to $q$ identical blocks in $RD_U$.
On the other hand, for all $x, j$ we have $|V_x \cap C_j| = |C_j| - |R_x \cap C_j| = r - \lrc = 1$, so each set $V_x$ corresponds to a parallel class of $CD_U$.
It follows that $U$ is resolvable.
\end{proof}

\begin{remark}
Theorem~\ref{thm:TA->AG} implies that any $((q + 1) \times q^2, q(q + 1))$-unordered triple array can be obtained as described in Remark~\ref{rm:aff->TA}, that is, using the \nameref{constr:res} with the unique resolution of an affine plane of order $q$ as one of the components.
On the other hand, by Theorem~\ref{thm:Agrawal-univ}, it can also be obtained using \nameref{constr:Agrawal}.
In fact, both essentially describe the same process in two different ways.
For a finite projective plane $\mathcal{P}$, applying \nameref{constr:Agrawal} to its dual $\mathcal{P}^*$, which has lines of $\mathcal{P}$ as points, and points of $\mathcal{P}$ as blocks, corresponds to the classic construction of an affine plane by deleting a line from a projective plane.
The column-sets in both cases correspond to the points of the affine plane, and the row-sets correspond to the `points at infinity', or, equivalently, to parallel classes of the affine lines, in the sense that each row-set `misses' one parallel class.
\end{remark}

\begin{remark}
There is a bijection between isomorphism classes of $((q + 1) \times q^2, q(q + 1))$-unordered triple arrays and isomorphism classes of affine planes of order $q$.
Indeed, isomorphic unordered triple arrays have isomorphic column designs, so they do correspond to isomorphic affine planes.
On the other hand, the symmetric 2-design required for the \nameref{constr:res} is a trivial 2-$(q + 1, q, q - 1)$ design, thus the construction would always produce isomorphic unordered triple arrays, no matter how the parallel classes, blocks and symbols of the affine plane are labeled.
\end{remark}

\begin{definition}
For an affine plane $\mathcal{A}$ of order $q$, denote by $U_\mathcal{A}$ the unique (up to isomorphism) resolvable $((q + 1) \times q^2, q(q + 1))$-unordered triple array whose column design is isomorphic to $\mathcal{A}$.
\end{definition}

We see that for parameter sets $((q + 1) \times q^2, q(q + 1))$, the \nameref{pbm:exUTA} is equivalent to the question of existence of projective (or equivalently, affine) planes of order $q$.
In the remainder of this section, we present two equivalent reformulations of the \nameref{pbm:UTAtoTA} for these parameter sets that may prove to be easier to work with: one in terms of derangements, and one in terms of hypergraphs and partiteness.

\subsection{A derangement problem}

\begin{problem}\label{pbm:derang}
For an affine plane of order $q$ with point set $P$, label the parallel classes of its lines as $\mathcal{L}_1, \dots, \mathcal{L}_{q + 1}$.
Find a way to assign a permutation $\sigma_p \in \mathfrak{S}_{q + 1}$ to each point $p \in P$ such that:
\begin{enumerate}[(a)]
\item\label{sigma:derang} $\sigma_p$ is a \textit{derangement} for any point $p \in P$, that is, $\sigma_p(i) \neq i$ for all $i$;
\item\label{sigma:2points} for any pair of distinct points $p, u \in P$, if the unique line through $p$ and $u$ belongs to the parallel class $\mathcal{L}_i$, then $\sigma_p(i) \neq \sigma_u(i)$.
\end{enumerate}
\end{problem}

\begin{proposition}\label{prop:equiv-derang}
For any affine plane $\mathcal{A}$ of order $q$, Problem~\ref{pbm:derang} has a solution if and only if the \nameref{pbm:UTAtoTA} has a solution for the $((q + 1) \times q^2, q(q + 1))$-unordered triple array $U_\mathcal{A}$.
\end{proposition}
\begin{proof}
Label the points of $\mathcal{A}$ as $p_1, \dots, p_{q^2}$.
Recalling the \nameref{constr:res}, the symbols of $U_\mathcal{A}$ are the lines of the plane, and, without loss of generality, the column-set $C_j$ contains the lines incident to the point $p_j$, and the row-set $R_i$ contains the lines from all parallel classes but $\mathcal{L}_i$.
Denote by $\mathcal{L}_i(j)$ the unique line such that $p_j \in \mathcal{L}_i(j) \in \mathcal{L}_i$.

Let $T$ be a solution to the \nameref{pbm:UTAtoTA}, that is, a triple array with $U_T = U_\mathcal{A}$.
Then, for each point $p_j$ define a permutation $\sigma_{p_j} \in \mathfrak{S}_{q + 1}$ as follows: for each parallel class $\mathcal{L}_i$, let $\sigma_{p_j}(i) = s$ if $\mathcal{L}_i(j)$ belongs to the cell $(s, j)$ of $T$.
These permutations form a solution to Problem~\ref{pbm:derang}: indeed, condition~\ref{sigma:derang} follows from the fact that $R_i \cap \mathcal{L}_i = \emptyset$, and condition~\ref{sigma:2points} holds because $T$ is binary, so no line appears two times in the same row.

Conversely, for a collection of permutations $\sigma_{p_j} \in \mathfrak{S}_{q + 1}$ forming a solution to Problem~\ref{pbm:derang}, define an $r \times c$ row-column design $T$ by placing $\mathcal{L}_i(j)$ in the cell $(\sigma_{p_j}(i), j)$ of $T$.
Repeating the arguments above in reverse, $T$ is a triple array with $U_T = U_\mathcal{A}$.
\end{proof}

It is not hard to show by case analysis that there are no $(3 \times 4, 6)$-triple arrays.
This, however, gives no insights into why this parameter set is exceptional with regards to Conjecture~\ref{conj:Agrawal}.
In the next corollary, as an indication of the usefulness of the reformulation in Problem~\ref{pbm:derang}, we give an alternative proof of this non-existence fact using Proposition~\ref{prop:equiv-derang}, which highlights one reason why this parameter set is special.

\begin{corollary}\label{cor:no346}
There are no $(3 \times 4, 6)$-triple arrays.
\end{corollary}
\begin{proof}
Combining Theorem~\ref{thm:TA->AG} with Proposition~\ref{prop:equiv-derang}, any $(3 \times 4, 6)$-triple array corresponds to a solution to Problem~\ref{pbm:derang} for $\AG(2, 2)$, the affine plane of order 2.
Such a solution would include 4 pairwise different derangements from $\mathfrak{S}_3$, however, there are only 2 derangements on 3 elements.
\end{proof}

The next parameter set in the family $((q + 1) \times q^2, q(q + 1))$ is $(4 \times 9, 12)$, corresponding to $q = 3$.
In this case, solving Problem~\ref{pbm:derang} requires 9 pairwise different derangements from $\mathfrak{S}_4$, and there are precisely that many derangements on 4 elements.
This provides some explanation for why the $(4 \times 9, 12)$-triple array is unique up to isotopism, see Table~\ref{tbl:extr-UTA-counts} in Section~\ref{sec:enum}.

With $q$ growing further, the number of derangements grows exponentially, quickly outpacing what is required for Problem~\ref{pbm:derang}, that is, the number of columns $q^2$.
For example, there are 44 derangements vs. 16 columns for $q = 4$, and 265 derangements vs. 25 columns for $q = 5$.
This indicates that solving Problem~\ref{pbm:derang} should become easier as $q$ grows, which is further supported by the fact that there are many non-isotopic $(5 \times 16, 20)$-triple arrays (see Table~\ref{tbl:extr-UTA-counts} in Section~\ref{sec:enum}).
Computationally we have found many solutions for the larger $q$ we have tested, but we have not been able to find a direct description of a solution to Problem~\ref{pbm:derang} for general $q$.

\subsection{A hypergraph problem}

To state the \nameref{pbm:UTAtoTA} for parameter sets $((q + 1) \times q^2, q(q + 1))$ in yet another way, we need to introduce some terminology from hypergraph theory.
A \textit{$k$-graph} $(V, E)$ on a \textit{vertex} set $V$ is a collection of $k$-subsets of $V$ called \textit{edges}.
It is \textit{linear} if any two distinct edges share at most one vertex.
It is $k$-\textit{partite} if one can partition the vertex set $V$ into $k$ parts such that each edge contains one vertex from each part.

\begin{problem}\label{pbm:partite}
For an affine plane of order $q$ with points $p_1, \dots, p_{q^2}$, label the parallel classes of its lines as $\mathcal{L}_1, \dots, \mathcal{L}_{q + 1}$, and denote by $\mathcal{L}_i(j)$ the unique line such that $p_j \in \mathcal{L}_i(j) \in \mathcal{L}_i$.
Define a $(q + 1)$-graph $H = (V, E)$ as follows: let
\[
V := \{v_{ij}\,:\,1 \leq i \leq q + 1, 1 \leq j \leq q^2\} \cup \{w_i\,:\,1 \leq i \leq q + 1\},
\]
for each point $p_j$ add to $E$ the edge $e(p_j) := \{v_{ij}\,:\,1 \leq i \leq q + 1\}$, for each line $l \in \mathcal{L}_i$ add to $E$ the edge $e(l) := \{v_{ij}\,:\,l = \mathcal{L}_i(j)\} \cup \{w_i\}$, and, finally, add to $E$ one more edge $e_0 := \{w_i\,:\, 1 \leq i \leq q + 1\}$.

Find a partition of $V$ into $q+1$ parts $V_1, V_2, \ldots, V_{q+1}$ such that any edge meets each $V_i$ in exactly one vertex.
\end{problem}

We note that $H$ has $(q+1)q^2 + (q+1)$ vertices and $q^2 + (q^2+q) + 1$ edges, and is linear.
If the partition sought in Problem~\ref{pbm:partite} can be found, this in other words means that the hypergraph $H$ is $(q+1)$-partite. We prove in the following proposition that Problem~\ref{pbm:partite} has a solution if and only if Problem~\ref{pbm:derang} does, and consequently, by Proposition~\ref{prop:equiv-derang}, if and only if
the \nameref{pbm:UTAtoTA} does.

\begin{proposition}\label{equiv:partite}
For any affine plane of order $q$, Problem~\ref{pbm:partite} has a solution if and only if Problem~\ref{pbm:derang} has a solution.
\end{proposition}

\begin{proof}
Suppose that $H$ is $(q + 1)$-partite and let $V_1, \dots, V_{q + 1}$ be the corresponding vertex partition.
The edge $e_0$ has one vertex in each part, so, without loss of generality, we may assume by relabeling vertices that $w_i \in V_i$.
For each point $p_j$, define $\sigma_{p_j}(i) := s$ if $v_{ij} \in V_s$.
Each edge $e(p_j)$ has one vertex in each part, so $\sigma_{p_j}$ is a permutation.
Each edge $e(l)$ has one vertex in each part, and $w_i \in V_i$, which implies conditions~\ref{sigma:2points} and~\ref{sigma:derang} from Problem~\ref{pbm:derang}, respectively.

Conversely, for a collection of permutations $\sigma_{p_j} \in \mathfrak{S}_{q + 1}$ forming a solution to Problem~\ref{pbm:derang}, define the partition of the vertex set of $H$ as $V_s := \{v_{ij}\,:\,1 \leq j \leq q^2, \sigma_{p_j}(i) = s\} \cup \{w_s\}$, $1 \leq s \leq q + 1$.
Repeating the arguments above in reverse, each edge of $H$ contains one vertex from each of $V_s$.
\end{proof}

\section{Concluding remarks}
\label{sec:concl}

\subsection{No triple arrays via \texorpdfstring{$\alpha$}{alpha}-resolutions}

For a positive integer $\alpha$, the following notions generalize parallel classes, resolutions and resolvable 2-designs.
An $\alpha$-\textit{parallel class} in a 2-design is a collection of blocks such that each point occurs in exactly $\alpha$ of these blocks.
An $\alpha$-\textit{resolution} of a 2-design is a partition of the collection of blocks into $\alpha$-parallel classes.
A 2-design that admits an $\alpha$-resolution is called $\alpha$-\textit{resolvable}.

A potential generalization of the \nameref{constr:res} could use an $\alpha$-resolution $\mathcal{B}$ of a 2-$(c, v, r, e, \lcc)$ design for some $\alpha > 1$.
Note that the number of $\alpha$-parallel classes in such an $\alpha$-resolution is $\frac{ve}{\alpha c} = \frac{r}{\alpha}$.
For the construction to produce an unordered triple array, the second component would instead of a symmetric 2-$(r, e, \lrrc)$ design need to be a 2-design with $r$ points and $\frac{r}{\alpha}$ blocks.
Unfortunately, due to \nameref{thm:fisher}, this is only possible when $\alpha = 1$, so this approach does not lead to any new unordered triple arrays.

\subsection{Open questions}

First, we reiterate the two versions of Agrawal's conjecture on the existence of extremal triple arrays.
Note that our computational results (see Observation~\ref{obs:UTAgivesTA}) provide some new evidence towards these conjectures.

\theoremstyle{plain}
\newtheorem*{conjagr}{Conjecture~\ref{conj:Agrawal}}
\begin{conjagr}
If there is a symmetric 2-$(r + c, r, \lcc)$ design with $r - \lcc > 2$, then there is an $(r \times c, r + c - 1)$-triple array.
\end{conjagr}

\newtheorem*{conjagr2}{Conjecture~\ref{conj:Agrawal-gen}}
\begin{conjagr2}
For any $(r \times c, r + c - 1)$-unordered triple array $U$ with $e > 2$, there is an $(r \times c, r + c - 1)$-triple array $T$ with $U_T = U$.
\end{conjagr2}

Resolving these conjectures in the special case of parameters $((q + 1) \times q^2, q(q + 1))$ amounts to solving either of the Problems~\ref{pbm:derang} and~\ref{pbm:partite}.
In particular, solving these problems for Galois affine planes would result in a new infinite series of triple arrays.

\begin{conjecture}
Problems~\ref{pbm:derang} and~\ref{pbm:partite} have solutions for the Galois affine plane $\AG(2, q)$ for any prime power $q > 2$.
Consequently, there is a $((q + 1) \times q^2, q(q + 1))$-triple array for any prime power $q > 2$.
\end{conjecture}

Theorem~\ref{thm:PG3q} provides an infinite family of non-extremal unordered triple arrays.
We conjecture that, for each parameter set in the family, at least some of these unordered triple arrays can be ordered.

\begin{conjecture}
For any prime power $q$, there is a $(\frac{q^3 - 1}{q - 1}, \frac{q^4 -1}{q - 1}, \frac{(q^4 - 1)(q^3 - 1)}{(q^2 - 1)(q - 1)})$-triple array.
\end{conjecture}

The next step towards this conjecture would be resolving the case $q = 3$, which corresponds to a $(13 \times 40, 130)$-triple array. It may well be that the \nameref{constr:res} can succeed, but, as discussed in Section~\ref{subsec:enum-non-extr}, other ideas may be needed.

\begin{problem}\label{pbm:1340130}
Construct a $(13 \times 40, 130)$-triple array.
\end{problem}

Since not all extremal triple arrays are resolvable, it would be natural to assume that the same holds in the non-extremal case.
However, the \nameref{constr:res} which is the only construction method currently available in this case, by definition only produces resolvable arrays.
Table~\ref{tbl:non-extr-non-res} lists all non-extremal parameter sets with $r \leq 30$ or $c \leq 30$ which are admissible for triple arrays but not for resolvable triple arrays.
Note that some of these parameter sets are admissible for quad arrays.

\begin{table}[!ht]
\begin{center}
\begin{tabular}{|rr|rrrr|rr|}
\hline
$r \times c$ & $v$ & $e$ & $\lrc$ & $\lrr$ & $\lcc$ & $\lrrc$ & $k$ \\
\hline
$16 \times 21$ & $56$ & $6$ & $6$ & $4$ & $7$ & $2$ & $7/2$\\
$16 \times 25$ & $100$ & $4$ & $4$ & $2$ & $5$ & $4/5$ & $25/4$\\
$16 \times 81$ & $216$ & $6$ & $6$ & $1$ & $27$ & $2$ & $27/2$\\
$16 \times 145$ & $232$ & $10$ & $10$ & $1$ & $87$ & $6$ & $29/2$\\
$21 \times 16$ & $56$ & $6$ & $6$ & $7$ & $4$ & $3/2$ & $8/3$\\
$21 \times 36$ & $126$ & $6$ & $6$ & $3$ & $9$ & $3/2$ & $6$\\
$25 \times 16$ & $100$ & $4$ & $4$ & $5$ & $2$ & $1/2$ & $4$\\
$36 \times 21$ & $126$ & $6$ & $6$ & $9$ & $3$ & $6/7$ & $7/2$\\
$81 \times 16$ & $216$ & $6$ & $6$ & $27$ & $1$ & $3/8$ & $8/3$\\
\hline
\end{tabular}
\end{center}
\caption{All non-extremal parameter sets with $r \leq 30$ or $c \leq 30$ admissible for triple arrays but not for resolvable triple arrays.} 
\label{tbl:non-extr-non-res}
\end{table}

\begin{problem}
Construct a triple array for any of the parameters in Table~\ref{tbl:non-extr-non-res}.
\end{problem}

Further, as we have seen in the extremal case, see Example~\ref{ex:7814-res-or-not} and Table~\ref{tbl:extr-QTA-counts}, a triple array does not have to be resolvable even when its parameter set is admissible for resolvable triple arrays. Despite our full enumeration of resolvable $(7 \times 15,35)$-triple arrays, it is therefore still an open question whether there are further, non-resolvable triple arrays on these parameters.

\begin{problem}
Construct a non-resolvable $(7 \times 15, 35)$-triple array.
\end{problem}

In Lemma~\ref{lm:res-transpose}, we showed that a set of parameters admissible for a resolvable triple array cannot be admissible for the transpose of a resolvable triple array, so no triple array can be resolvable in both orientations.
We can ask the same question for quad arrays.

\begin{question}\label{q:quad_transpose}
Can the transpose of a quad array also be a quad array?
\end{question}

One way to answer Question~\ref{q:quad_transpose} in the negative is to show that there are no parameter sets that are admissible both for a quad array and the transpose of a quad array. We have checked computationally by exhaustive generation of admissible parameters that  there are no parameter sets admissible for quad arrays in both orientations for $e \leq 10^5$.

As noted in Observation~\ref{obs:QA-vs-res}, all (unordered) quad arrays we constructed computationally turned out to also be resolvable, which leads to the next question.

\begin{question}
Is there a non-resolvable quad array?    
\end{question}

\section*{Acknowledgements}

Alexey Gordeev was supported by Kempe foundation grant JCSMK23-0058.

\bibliographystyle{abbrv}
\bibliography{main}

@article{smithConstructionYoudenSquares1948,
  title = {The {{Construction}} of {{Youden Squares}}},
  author = {Smith, C. A. B. and Hartley, H. O.},
  year = 1948,
  journal = {Journal of the Royal Statistical Society: Series B (Methodological)},
  volume = {10},
  number = {2},
  pages = {262--263},
  doi = {10.1111/j.2517-6161.1948.tb00015.x},
  urldate = {2023-07-28}
}

@article {DennistonPackings,
    AUTHOR = {Denniston, Ralph H. F.},
     TITLE = {Some packings of projective spaces},
   JOURNAL = {Atti della Accademia Nazionale dei Lincei. Rendiconti. Classe
              di Scienze Fisiche, Matematiche e Naturali. Serie VIII},
    VOLUME = {52},
      YEAR = {1972},
     PAGES = {36--40},
      ISSN = {0392-7881},
   MRCLASS = {50D30},
  MRNUMBER = {331207},
MRREVIEWER = {Stanley\ E.\ Payne},
}

@article {BeutelspacherParallelisms,
    AUTHOR = {Beutelspacher, Albrecht},
     TITLE = {On parallelisms in finite projective spaces},
   JOURNAL = {Geometriae Dedicata},
  FJOURNAL = {Geometriae Dedicata},
    VOLUME = {3},
      YEAR = {1974},
     PAGES = {35--40},
   MRCLASS = {50D30},
  MRNUMBER = {341270},
       DOI = {10.1007/BF00181359},
       URL = {https://doi.org/10.1007/BF00181359},
}

@article {BettenPackings,
    AUTHOR = {Betten, Anton},
     TITLE = {The packings of {$\mathrm{PG}(3,3)$}},
   JOURNAL = {Designs, Codes and Cryptography},
    VOLUME = {79},
      YEAR = {2016},
    NUMBER = {3},
     PAGES = {583--595},
      ISSN = {0925-1022,1573-7586},
   MRCLASS = {51E23 (05E18)},
  MRNUMBER = {3489759},
       DOI = {10.1007/s10623-015-0074-6},
       URL = {https://doi.org/10.1007/s10623-015-0074-6},
}

@article {cole_kirkman_parades,
    AUTHOR = {Cole, F. N.},
     TITLE = {Kirkman parades},
   JOURNAL = {Bulletin of the American Mathematical Society},
    VOLUME = {28},
      YEAR = {1922},
    NUMBER = {9},
     PAGES = {435--437},
      ISSN = {0002-9904},
   MRCLASS = {99-04},
  MRNUMBER = {1560613},
       DOI = {10.1090/S0002-9904-1922-03599-9},
       URL = {https://doi.org/10.1090/S0002-9904-1922-03599-9},
}

@inproceedings {yucas_structure_7x15,
    AUTHOR = {Yucas, Joseph L.},
     TITLE = {The structure of a {$7\times 15$} triple array},
 BOOKTITLE = {Proceedings of the {T}hirty-third {S}outheastern
              {I}nternational {C}onference on {C}ombinatorics, {G}raph
              {T}heory and {C}omputing ({B}oca {R}aton, {FL}, 2002), Congressus Numerantium},
   JOURNAL = {Congr. Numer.},
  FJOURNAL = {Congressus Numerantium. A Conference Journal on Numerical
              Themes},
    VOLUME = {154},
      YEAR = {2002},
     PAGES = {43--47},
      ISSN = {0384-9864},
   MRCLASS = {05B15 (05B07)},
  MRNUMBER = {1980027},
}

@incollection{preeceNonorthogonalGraecolatinDesigns1976,
  title = {Non-orthogonal {G}raeco-{L}atin designs},
  booktitle = {Combinatorial {{Mathematics IV}}. Lecture Notes in Mathematics},
  author = {Preece, D. A.},
  editor = {Casse, Louis R. A. and Wallis, Walter D.},
  year = 1976,
  volume = {560},
  pages = {7--26},
  publisher = {Springer, Berlin, Heidelberg},
  doi = {10.1007/BFb0097364},
  urldate = {2025-09-12},
  copyright = {http://www.springer.com/tdm},
  isbn = {978-3-540-08053-4 978-3-540-37537-1}
}

@article{mcsorleyDoubleArraysTriple2005a,
  title = {Double Arrays, Triple Arrays and Balanced Grids},
  author = {McSorley, John P. and Phillips, N. C. K. and Wallis, W. D. and Yucas, J. L.},
  year = {2005},
  journal = {Designs, Codes and Cryptography},
  volume = {35},
  number = {1},
  pages = {21--45},
  doi = {10.1007/s10623-003-6149-9},
  urldate = {2023-07-27}
}

@article{jagerEnumerationConstructionRowColumn2025,
  title = {Enumeration and {{Construction}} of {{Row}}-{{Column Designs}}},
  author = {J{\"a}ger, Gerold and Markstr{\"o}m, Klas and Shcherbak, Denys and {\"O}hman, Lars-Daniel},
  year = {2025},
  journal = {Journal of Combinatorial Designs},
  volume = {33},
  number = {9},
  pages = {357--372},
  publisher = {Wiley},
  doi = {10.1002/jcd.21991},
  urldate = {2025-07-15},
  copyright = {http://creativecommons.org/licenses/by/4.0/}
}

@article{agrawalMethodsConstructionDesigns1966,
  title = {Some Methods of Construction of Designs for Two-Way Elimination of Heterogeneity, 1},
  author = {Agrawal, Hiralal},
  year = {1966},
  journal = {Journal of the American Statistical Association},
  volume = {61},
  number = {316},
  eprint = {2283205},
  eprinttype = {jstor},
  pages = {1153--1171},
  publisher = {[American Statistical Association, Taylor \& Francis, Ltd.]},
  doi = {10.2307/2283205},
  urldate = {2023-07-27}
}

@article{bagchiTwowayDesigns1998,
  title = {On Two-Way Designs},
  author = {Bagchi, Sunanda},
  year = {1998},
  journal = {Graphs and Combinatorics},
  volume = {14},
  number = {4},
  pages = {313--319},
  doi = {10.1007/PL00021181},
  urldate = {2024-06-17}
}

@article{bailey1994extremal,
  title={Extremal Row--Column Designs with Maximal Balance and Adjusted Orthogonality},
  author={Bailey, R. A. and Heidtmann, P.},
  journal={Preprint, Goldsmiths’ College, University of London},
  year={1994}
}

@manual{sagemath25,
  Key          = {SageMath},
  Author       = {{The Sage Developers}},
  Title        = {{S}ageMath, the {S}age {M}athematics {S}oftware {S}ystem ({V}ersion 10.7)},
  note         = {\url{https://www.sagemath.org}},
  Year         = {2025}
}

@manual{gordeevZenodoLibrary25,
    title={A library of triple arrays and unordered triple arrays},
    note={Zenodo. \url{https://doi.org/10.5281/zenodo.17854868}},
    author={A. Gordeev},
    year={2025}
}

@article{gordeevTripleArrays2026,
  title = {Near Triple Arrays},
  author = {Gordeev, Alexey and Markstr{\"o}m, Klas and {\"O}hman, Lars-Daniel},
  year = 2026,
  journal = {Journal of Combinatorial Theory, Series A},
  volume = {219},
  pages = {106121},
  doi = {10.1016/j.jcta.2025.106121},
  urldate = {2025-10-07}
}

@incollection{baileyRelationsPartitions2017,
  title = {Relations among Partitions},
  booktitle = {Surveys in {{Combinatorics}} 2017},
  author = {Bailey, R. A.},
  editor = {Claesson, Anders and Dukes, Mark and Kitaev, Sergey and Manlove, David and Meeks, Kitty},
  year = {2017},
  edition = {1},
  pages = {1--86},
  publisher = {Cambridge University Press},
  doi = {10.1017/9781108332699.002},
  urldate = {2025-01-24},
  isbn = {978-1-108-33269-9 978-1-108-41313-8}
}

@article{mckayPracticalGraphIsomorphism2014,
  title = {Practical Graph Isomorphism, {{II}}},
  author = {McKay, Brendan D. and Piperno, Adolfo},
  year = {2014},
  journal = {Journal of Symbolic Computation},
  volume = {60},
  pages = {94--112},
  doi = {10.1016/j.jsc.2013.09.003},
  urldate = {2024-02-12}
}

@article{preecePaleyTripleArrays2005,
  title = {Paley Triple Arrays},
  author = {Preece, Donald A. and Wallis, Walter D. and Yucas, Joseph L.},
  year = {2005},
  journal = {The Australasian Journal of Combinatorics},
  volume = {33},
  pages = {237--246},
  urldate = {2023-11-17}
}

@article{nilsonTripleArraysDifference2017,
  title = {Triple Arrays from Difference Sets},
  author = {Nilson, Tomas and Cameron, Peter J.},
  year = {2017},
  journal = {Journal of Combinatorial Designs},
  volume = {25},
  number = {11},
  pages = {494--506},
  doi = {10.1002/jcd.21569},
  urldate = {2024-06-29},
  copyright = {{\copyright} 2017 Wiley Periodicals, Inc.}
}

@article{nilsonTripleArraysYouden2015,
  title = {Triple Arrays and {{Youden}} Squares},
  author = {Nilson, Tomas and {\"O}hman, Lars-Daniel},
  year = {2015},
  journal = {Designs, Codes and Cryptography},
  volume = {75},
  number = {3},
  pages = {429--451},
  doi = {10.1007/s10623-014-9926-8},
  urldate = {2023-08-10}
}

@phdthesis{mulder1917kirkman,
  type = {Thesis},
  title = {Kirkman-systemen},
  author = {Mulder, Pieter},
  year = {1917},
  school = {Rijksuniversiteit Groningen, Leiden},
}

@book{bethDesignTheory1999,
  title = {Design {{Theory}}},
  author = {Beth, Thomas and Jungnickel, Deiter and Lenz, Hanfried},
  year = {1999},
  edition = {2nd},
  publisher = {Cambridge University Press},
  doi = {10.1017/CBO9780511549533},
  urldate = {2024-10-30},
  copyright = {https://www.cambridge.org/core/terms},
  isbn = {978-0-521-44432-3 978-0-511-54953-3}
}

@article{fon-der-flaassArraysDistinctRepresentatives1997,
  title = {Arrays of Distinct Representatives --- a Very Simple {{NP-complete}} Problem},
  author = {{Fon-Der-Flaass}, Dmitri G.},
  year = {1997},
  journal = {Discrete Mathematics},
  volume = {171},
  number = {1},
  pages = {295--298},
  doi = {10.1016/S0012-365X(97)89167-4},
  urldate = {2023-08-14}
}

@book{knuth2025fascicle7,
  author    = {Donald E. Knuth},
  title     = {The Art of Computer Programming, Volume 4, Fascicle 7: Constraint Satisfaction},
  publisher = {Addison--Wesley / Pearson},
  year      = {2025},
}

@article{knuthDancingLinks2000,
  title = {Dancing Links},
  author = {Knuth, Donald E.},
  year = {2000},
  publisher = {arXiv},
  journal = {arXiv:cs/0011047},
  doi = {10.48550/ARXIV.CS/0011047},
  urldate = {2025-03-20},
  copyright = {Assumed arXiv.org perpetual, non-exclusive license to distribute this article for submissions made before January 2004}
}

@manual{heinlein2023,
    author={Heinlein, Daniel and Ivanov, Andrei and McKay, Brendan and {\"O}sterg\r{a}rd, Patric R. J.},
    title={A library of combinatorial 2-designs},
    year={2023},
    note={Zenodo. \url{https://doi.org/10.5281/zenodo.8262681}}
}

@article{aschbacherCollineationGroupsSymmetric1971,
  title = {On Collineation Groups of Symmetric Block Designs},
  author = {Aschbacher, Michael},
  year = 1971,
  journal = {Journal of Combinatorial Theory, Series A},
  volume = {11},
  number = {3},
  pages = {272--281},
  doi = {10.1016/0097-3165(71)90054-9},
  urldate = {2025-11-24},
  copyright = {https://www.elsevier.com/tdm/userlicense/1.0/}
}

@book{seberryHadamardMatricesConstructions2020,
  title = {Hadamard Matrices: Constructions Using Number Theory and Algebra},
  shorttitle = {Hadamard Matrices},
  author = {Seberry, Jennifer R. and Yamada, Mieko},
  year = 2020,
  publisher = {Wiley},
  address = {Hoboken (N.J.)},
  isbn = {978-1-119-52013-9 978-1-119-52027-6 978-1-119-52025-2},
  lccn = {511.6}
}

@article{paleyOrthogonalMatrices1933,
  title = {On {{Orthogonal Matrices}}},
  author = {Paley, R. E. A. C.},
  year = 1933,
  journal = {Journal of Mathematics and Physics},
  volume = {12},
  number = {1-4},
  pages = {311--320},
  doi = {10.1002/sapm1933121311},
  urldate = {2025-11-26},
  copyright = {http://onlinelibrary.wiley.com/termsAndConditions\#vor}
}

@incollection{mathon15FamilyBIBDs1989,
  title = {On the $(15,5,\lambda)$-{{Family}} of {{BIBDs}}},
  booktitle = {Annals of {{Discrete Mathematics}}},
  author = {Mathon, R. and Rosa, A.},
  editor = {Hartman, A.},
  year = 1989,
  series = {Combinatorial {{Designs}}---{{A Tribute}} to {{Haim Hanani}}},
  volume = {42},
  pages = {205--216},
  publisher = {Elsevier},
  doi = {10.1016/S0167-5060(08)70107-9},
  urldate = {2023-11-10}
}

@article{pavoneSevenNonisomorphicSolutions2023,
  title = {On the Seven Non-Isomorphic Solutions of the Fifteen Schoolgirl Problem},
  author = {Pavone, Marco},
  year = 2023,
  journal = {Discrete Mathematics},
  volume = {346},
  number = {6},
  pages = {113316},
  doi = {10.1016/j.disc.2023.113316},
  urldate = {2023-11-02}
}
\newpage

\appendix
\counterwithin{figure}{section}
\counterwithin{table}{section}

\clearpage
\section{Counts of resolvable \texorpdfstring{$(7 \times 15, 35)$}{(7 x 15, 35)}-triple arrays}\label{ap:7-15-35}

\begin{table}[!ht]
\begin{center}
\begin{tabular}{|rr|rrrrrrrr|r|}
\multicolumn{11}{c}{\multirow{2}{*}{Kirkman parade 1a}} \\
\multicolumn{11}{c}{} \\
\hline
\multicolumn{2}{|c|}{UTA} & $U_{0}$ & $U_{1}$ & $U_{2}$ & $U_{3}$ &  &  &  &  & All \\
\multicolumn{2}{|c|}{$|\Aut U_i|$} & $12$ & $21$ & $24$ & $168$ &  &  &  &  & \\
\hline
\multicolumn{2}{|c|}{Total \#} & $0$ & $0$ & $0$ & $0$ &  &  &  &  & 0 \\
\hline
\multicolumn{11}{c}{\multirow{2}{*}{Kirkman parade 1b}} \\
\multicolumn{11}{c}{} \\
\hline
\multicolumn{2}{|c|}{UTA} & $U_{4}$ & $U_{5}$ & $U_{6}$ & $U_{7}$ &  &  &  &  & All \\
\multicolumn{2}{|c|}{$|\Aut U_i|$} & $24$ & $12$ & $21$ & $168$ &  &  &  &  & \\
\hline
\multicolumn{2}{|c|}{Total \#} & $0$ & $0$ & $2$ & $1$ &  &  &  &  & 3 \\
\hline
$|\Aut|$  & $1$ &  &  & $2$ & $1$ &  &  &  &  & $3$ \\
\hline
\multicolumn{11}{c}{\multirow{2}{*}{Kirkman parade 7a}} \\
\multicolumn{11}{c}{} \\
\hline
\multicolumn{2}{|c|}{UTA} & $U_{8}$ & $U_{9}$ & $U_{10}$ & $U_{11}$ & $U_{12}$ & $U_{13}$ &  &  & All \\
\multicolumn{2}{|c|}{$|\Aut U_i|$} & $4$ & $4$ & $3$ & $3$ & $24$ & $24$ &  &  & \\
\hline
\multicolumn{2}{|c|}{Total \#} & $2$ & $2$ & $10$ & $5$ & $4$ & $1$ &  &  & 24 \\
\hline
$|\Aut|$  & $1$ & $2$ & $2$ & $4$ & $4$ &  &  &  &  & $12$ \\
 & $3$ &  &  & $6$ & $1$ & $4$ & $1$ &  &  & $12$ \\
\hline
\multicolumn{11}{c}{\multirow{2}{*}{Kirkman parade 7b}} \\
\multicolumn{11}{c}{} \\
\hline
\multicolumn{2}{|c|}{UTA} & $U_{14}$ & $U_{15}$ & $U_{16}$ & $U_{17}$ & $U_{18}$ & $U_{19}$ &  &  & All \\
\multicolumn{2}{|c|}{$|\Aut U_i|$} & $4$ & $4$ & $3$ & $3$ & $24$ & $24$ &  &  & \\
\hline
\multicolumn{2}{|c|}{Total \#} & $1$ & $0$ & $1$ & $2$ & $0$ & $0$ &  &  & 4 \\
\hline
$|\Aut|$  & $1$ & $1$ &  & $1$ & $2$ &  &  &  &  & $4$ \\
\hline
\multicolumn{11}{c}{\multirow{2}{*}{Kirkman parade 19a}} \\
\multicolumn{11}{c}{} \\
\hline
\multicolumn{2}{|c|}{UTA} & $U_{20}$ & $U_{21}$ & $U_{22}$ & $U_{23}$ & $U_{24}$ & $U_{25}$ & $U_{26}$ & $U_{27}$ & All \\
\multicolumn{2}{|c|}{$|\Aut U_i|$} & $12$ & $4$ & $12$ & $12$ & $3$ & $1$ & $3$ & $3$ & \\
\hline
\multicolumn{2}{|c|}{Total \#} & $2$ & $1$ & $2$ & $0$ & $3$ & $6$ & $3$ & $4$ & 21 \\
\hline
$|\Aut|$  & $1$ & $2$ & $1$ & $2$ &  & $3$ & $6$ & $3$ & $4$ & $21$ \\
\hline
\multicolumn{11}{c}{\multirow{2}{*}{Kirkman parade 19b}} \\
\multicolumn{11}{c}{} \\
\hline
\multicolumn{2}{|c|}{UTA} & $U_{28}$ & $U_{29}$ & $U_{30}$ & $U_{31}$ & $U_{32}$ & $U_{33}$ & $U_{34}$ & $U_{35}$ & All \\
\multicolumn{2}{|c|}{$|\Aut U_i|$} & $4$ & $12$ & $12$ & $12$ & $1$ & $3$ & $3$ & $3$ & \\
\hline
\multicolumn{2}{|c|}{Total \#} & $1$ & $1$ & $2$ & $0$ & $6$ & $3$ & $5$ & $3$ & 21 \\
\hline
$|\Aut|$  & $1$ & $1$ & $1$ & $2$ &  & $6$ & $3$ & $5$ & $3$ & $21$ \\
\hline
\multicolumn{11}{c}{\multirow{2}{*}{Kirkman parade 61}} \\
\multicolumn{11}{c}{} \\
\hline
\multicolumn{2}{|c|}{UTA} & $U_{36}$ & $U_{37}$ & $U_{38}$ & $U_{39}$ & $U_{40}$ & $U_{41}$ &  &  & All \\
\multicolumn{2}{|c|}{$|\Aut U_i|$} & $3$ & $3$ & $3$ & $3$ & $21$ & $21$ &  &  & \\
\hline
\multicolumn{2}{|c|}{Total \#} & $2$ & $3$ & $2$ & $4$ & $1$ & $0$ &  &  & 12 \\
\hline
$|\Aut|$  & $1$ & $2$ & $3$ & $2$ & $4$ & $1$ &  &  &  & $12$ \\
\hline
\end{tabular}
\end{center}
\caption{The number of non-isotopic resolvable $(7 \times 15, 35)$-triple arrays sorted by the corresponding Kirkman parade, the underlying unordered triple array, and the autotopism group order.
The unordered triple arrays are numbered in the order that they are given in~\cite{gordeevZenodoLibrary25}.
The naming convention for the Kirkman parades is borrowed from~\cite[Table 1]{pavoneSevenNonisomorphicSolutions2023}.}
\label{tbl:extended-7-15-35}
\end{table}

\begin{figure}[!ht]
\begin{center}
\begin{tabular}{c}
\begin{tabular}{|*{15}{c}|}
\hline
$ 1$ & $17$ & $18$ & $22$ & $16$ & $19$ & $25$ & $ 3$ & $23$ & $20$ & $24$ & $ 5$ & $21$ & $ 2$ & $ 4$ \\
$11$ & $ 1$ & $28$ & $14$ & $13$ & $ 4$ & $ 5$ & $27$ & $15$ & $ 2$ & $30$ & $29$ & $ 3$ & $26$ & $12$ \\
$ 6$ & $ 7$ & $ 1$ & $ 2$ & $ 3$ & $10$ & $31$ & $33$ & $ 4$ & $35$ & $ 5$ & $32$ & $34$ & $ 8$ & $ 9$ \\
$21$ & $12$ & $13$ & $ 9$ & $ 7$ & $22$ & $15$ & $14$ & $ 6$ & $23$ & $11$ & $10$ & $ 8$ & $24$ & $25$ \\
$16$ & $27$ & $ 8$ & $28$ & $29$ & $30$ & $ 7$ & $ 6$ & $17$ & $10$ & $ 9$ & $18$ & $20$ & $19$ & $26$ \\
$31$ & $32$ & $33$ & $16$ & $35$ & $13$ & $20$ & $19$ & $34$ & $11$ & $17$ & $14$ & $12$ & $15$ & $18$ \\
$26$ & $22$ & $23$ & $34$ & $24$ & $31$ & $28$ & $25$ & $29$ & $27$ & $33$ & $21$ & $30$ & $32$ & $35$ \\
\hline
\end{tabular}\\
\\
\begin{tabular}{|*{15}{c}|}
\hline
$16$ & $17$ & $ 1$ & $22$ & $24$ & $19$ & $20$ & $25$ & $ 4$ & $23$ & $ 5$ & $21$ & $ 3$ & $ 2$ & $18$ \\
$ 1$ & $27$ & $28$ & $14$ & $13$ & $30$ & $ 5$ & $ 3$ & $15$ & $ 2$ & $11$ & $29$ & $12$ & $26$ & $ 4$ \\
$31$ & $ 1$ & $33$ & $ 2$ & $ 3$ & $ 4$ & $ 7$ & $ 6$ & $34$ & $10$ & $ 9$ & $ 5$ & $ 8$ & $32$ & $35$ \\
$ 6$ & $12$ & $13$ & $ 9$ & $ 7$ & $22$ & $15$ & $14$ & $23$ & $11$ & $24$ & $10$ & $21$ & $ 8$ & $25$ \\
$26$ & $ 7$ & $ 8$ & $16$ & $29$ & $10$ & $28$ & $27$ & $ 6$ & $20$ & $17$ & $18$ & $30$ & $19$ & $ 9$ \\
$11$ & $32$ & $18$ & $34$ & $16$ & $13$ & $31$ & $19$ & $17$ & $35$ & $33$ & $14$ & $20$ & $15$ & $12$ \\
$21$ & $22$ & $23$ & $28$ & $35$ & $31$ & $25$ & $33$ & $29$ & $27$ & $30$ & $32$ & $34$ & $24$ & $26$ \\
\hline
\end{tabular}\\
\\
\begin{tabular}{|*{15}{c}|}
\hline
$16$ & $17$ & $ 1$ & $ 2$ & $19$ & $22$ & $18$ & $ 3$ & $ 4$ & $23$ & $ 5$ & $21$ & $20$ & $24$ & $25$ \\
$26$ & $ 1$ & $13$ & $29$ & $11$ & $ 4$ & $15$ & $14$ & $12$ & $ 2$ & $28$ & $ 5$ & $ 3$ & $27$ & $30$ \\
$ 1$ & $ 7$ & $ 8$ & $ 9$ & $ 3$ & $33$ & $ 5$ & $ 6$ & $35$ & $10$ & $31$ & $34$ & $32$ & $ 2$ & $ 4$ \\
$11$ & $12$ & $23$ & $22$ & $24$ & $14$ & $ 7$ & $25$ & $ 6$ & $15$ & $ 9$ & $10$ & $21$ & $ 8$ & $13$ \\
$ 6$ & $27$ & $28$ & $18$ & $ 7$ & $10$ & $26$ & $17$ & $29$ & $30$ & $20$ & $19$ & $ 8$ & $16$ & $ 9$ \\
$31$ & $32$ & $18$ & $11$ & $33$ & $20$ & $35$ & $34$ & $19$ & $17$ & $12$ & $13$ & $15$ & $14$ & $16$ \\
$21$ & $22$ & $33$ & $34$ & $30$ & $26$ & $25$ & $28$ & $23$ & $31$ & $24$ & $27$ & $29$ & $35$ & $32$ \\
\hline
\end{tabular}
\end{tabular}
\end{center}
\caption{Three pairwise non-isotopic resolvable $(7 \times 15, 35)$-triple arrays corresponding to the Kirkman parade `1b', one of two resolutions of $\PG(3, 2)$.
The top array is isotopic to the triple array in Figure~\ref{fig:MPWY-71535}, originally found by McSorley, Phillips, Wallis and Yucas~\cite{mcsorleyDoubleArraysTriple2005a}.
The top and middle triple arrays share the same underlying unordered triple array, labeled as $U_6$ in Table~\ref{tbl:extended-7-15-35}, while the bottom triple array is the unique (up to isotopism) ordering of $U_7$.}
\label{fig:PG-71535}
\end{figure}

\clearpage
\section{A \texorpdfstring{$(21 \times 15, 63)$}{(21 x 15, 63)}-triple array}\label{ap:21-15-63}

\begin{figure}[!ht]
\begin{center}
\begin{tabular}{c}
\begin{tabular}{c|*{15}{c}|}
\multicolumn{1}{c}{} & \multicolumn{1}{c}{\scriptsize $C_{1}$} & \multicolumn{1}{c}{\scriptsize $C_{2}$} & \multicolumn{1}{c}{\scriptsize $C_{3}$} & \multicolumn{1}{c}{\scriptsize $C_{4}$} & \multicolumn{1}{c}{\scriptsize $C_{5}$} & \multicolumn{1}{c}{\scriptsize $C_{6}$} & \multicolumn{1}{c}{\scriptsize $C_{7}$} & \multicolumn{1}{c}{\scriptsize $C_{8}$} & \multicolumn{1}{c}{\scriptsize $C_{9}$} & \multicolumn{1}{c}{\scriptsize $C_{10}$} & \multicolumn{1}{c}{\scriptsize $C_{11}$} & \multicolumn{1}{c}{\scriptsize $C_{12}$} & \multicolumn{1}{c}{\scriptsize $C_{13}$} & \multicolumn{1}{c}{\scriptsize $C_{14}$} & \multicolumn{1}{c}{\scriptsize $C_{15}$} \\
\hhline{~---------------}
{\scriptsize $R_{ 1}$} & $39$ & $ 2$ & $30$ & $ 1$ & $55$ & $38$ & $29$ & $56$ & $28$ & $ 3$ & $24$ & $37$ & $57$ & $22$ & $23$ \\
{\scriptsize $R_{ 2}$} & $49$ & $16$ & $ 1$ & $50$ & $51$ & $45$ & $17$ & $35$ & $ 2$ & $34$ & $43$ & $44$ & $18$ & $ 3$ & $36$ \\
{\scriptsize $R_{ 3}$} & $54$ & $52$ & $31$ & $46$ & $53$ & $47$ & $ 2$ & $ 1$ & $32$ & $13$ & $15$ & $48$ & $33$ & $14$ & $ 3$ \\
{\scriptsize $R_{ 4}$} & $ 1$ & $20$ & $40$ & $42$ & $41$ & $27$ & $58$ & $26$ & $60$ & $59$ & $ 3$ & $25$ & $ 2$ & $19$ & $21$ \\
{\scriptsize $R_{ 5}$} & $61$ & $ 7$ & $ 6$ & $10$ & $ 3$ & $ 2$ & $ 5$ & $ 9$ & $ 8$ & $11$ & $ 4$ & $ 1$ & $62$ & $63$ & $12$ \\
{\scriptsize $R_{ 6}$} & $26$ & $44$ & $54$ & $ 5$ & $45$ & $23$ & $24$ & $43$ & $25$ & $22$ & $53$ & $ 6$ & $27$ & $ 4$ & $52$ \\
{\scriptsize $R_{ 7}$} & $ 6$ & $35$ & $14$ & $15$ & $ 5$ & $42$ & $36$ & $ 4$ & $34$ & $56$ & $55$ & $13$ & $41$ & $57$ & $40$ \\
{\scriptsize $R_{ 8}$} & $28$ & $ 4$ & $46$ & $21$ & $19$ & $50$ & $49$ & $51$ & $ 6$ & $20$ & $47$ & $29$ & $30$ & $48$ & $ 5$ \\
{\scriptsize $R_{ 9}$} & $16$ & $31$ & $37$ & $59$ & $38$ & $ 4$ & $32$ & $60$ & $17$ & $ 6$ & $39$ & $33$ & $ 5$ & $58$ & $18$ \\
{\scriptsize $R_{10}$} & $41$ & $51$ & $50$ & $32$ & $22$ & $ 7$ & $ 8$ & $23$ & $40$ & $42$ & $31$ & $24$ & $49$ & $ 9$ & $33$ \\
{\scriptsize $R_{11}$} & $ 7$ & $37$ & $ 9$ & $34$ & $36$ & $53$ & $21$ & $52$ & $39$ & $54$ & $ 8$ & $19$ & $20$ & $35$ & $38$ \\
{\scriptsize $R_{12}$} & $47$ & $48$ & $25$ & $27$ & $16$ & $56$ & $55$ & $17$ & $57$ & $ 8$ & $18$ & $ 9$ & $46$ & $26$ & $ 7$ \\
{\scriptsize $R_{13}$} & $58$ & $29$ & $59$ & $ 7$ & $ 8$ & $13$ & $45$ & $15$ & $44$ & $43$ & $30$ & $60$ & $ 9$ & $28$ & $14$ \\
{\scriptsize $R_{14}$} & $15$ & $10$ & $11$ & $37$ & $26$ & $12$ & $25$ & $38$ & $14$ & $27$ & $51$ & $50$ & $13$ & $39$ & $49$ \\
{\scriptsize $R_{15}$} & $12$ & $42$ & $16$ & $53$ & $28$ & $17$ & $41$ & $11$ & $54$ & $29$ & $40$ & $18$ & $10$ & $52$ & $30$ \\
{\scriptsize $R_{16}$} & $21$ & $55$ & $20$ & $57$ & $10$ & $31$ & $11$ & $33$ & $19$ & $32$ & $12$ & $56$ & $43$ & $45$ & $44$ \\
{\scriptsize $R_{17}$} & $35$ & $24$ & $22$ & $23$ & $60$ & $36$ & $47$ & $46$ & $10$ & $48$ & $34$ & $12$ & $58$ & $11$ & $59$ \\
{\scriptsize $R_{18}$} & $22$ & $14$ & $62$ & $16$ & $13$ & $19$ & $15$ & $21$ & $23$ & $17$ & $20$ & $63$ & $24$ & $18$ & $61$ \\
{\scriptsize $R_{19}$} & $31$ & $27$ & $36$ & $28$ & $33$ & $30$ & $63$ & $29$ & $62$ & $61$ & $26$ & $34$ & $35$ & $32$ & $25$ \\
{\scriptsize $R_{20}$} & $44$ & $61$ & $45$ & $43$ & $46$ & $62$ & $37$ & $40$ & $47$ & $38$ & $63$ & $41$ & $39$ & $42$ & $48$ \\
{\scriptsize $R_{21}$} & $56$ & $60$ & $55$ & $63$ & $61$ & $58$ & $52$ & $62$ & $51$ & $49$ & $59$ & $53$ & $54$ & $50$ & $57$ \\
\hhline{~---------------}
\end{tabular}\\
\\
\renewcommand{\arraystretch}{0.4}
\setlength\tabcolsep{0.5pt}
\begin{tabular}{c|*{21}{c}|*{15}{c}|}
\multicolumn{1}{c}{} & \multicolumn{1}{c}{\tiny $R_{1}$} & \multicolumn{1}{c}{\tiny $R_{2}$} & \multicolumn{1}{c}{\tiny $R_{3}$} & \multicolumn{1}{c}{\tiny $R_{4}$} & \multicolumn{1}{c}{\tiny $R_{5}$} & \multicolumn{1}{c}{\tiny $R_{6}$} & \multicolumn{1}{c}{\tiny $R_{7}$} & \multicolumn{1}{c}{\tiny $R_{8}$} & \multicolumn{1}{c}{\tiny $R_{9}$} & \multicolumn{1}{c}{\tiny $R_{10}$} & \multicolumn{1}{c}{\tiny $R_{11}$} & \multicolumn{1}{c}{\tiny $R_{12}$} & \multicolumn{1}{c}{\tiny $R_{13}$} & \multicolumn{1}{c}{\tiny $R_{14}$} & \multicolumn{1}{c}{\tiny $R_{15}$} & \multicolumn{1}{c}{\tiny $R_{16}$} & \multicolumn{1}{c}{\tiny $R_{17}$} & \multicolumn{1}{c}{\tiny $R_{18}$} & \multicolumn{1}{c}{\tiny $R_{19}$} & \multicolumn{1}{c}{\tiny $R_{20}$} & \multicolumn{1}{c}{\tiny $R_{21}$} & \multicolumn{1}{c}{\tiny $C_{1}$} & \multicolumn{1}{c}{\tiny $C_{2}$} & \multicolumn{1}{c}{\tiny $C_{3}$} & \multicolumn{1}{c}{\tiny $C_{4}$} & \multicolumn{1}{c}{\tiny $C_{5}$} & \multicolumn{1}{c}{\tiny $C_{6}$} & \multicolumn{1}{c}{\tiny $C_{7}$} & \multicolumn{1}{c}{\tiny $C_{8}$} & \multicolumn{1}{c}{\tiny $C_{9}$} & \multicolumn{1}{c}{\tiny $C_{10}$} & \multicolumn{1}{c}{\tiny $C_{11}$} & \multicolumn{1}{c}{\tiny $C_{12}$} & \multicolumn{1}{c}{\tiny $C_{13}$} & \multicolumn{1}{c}{\tiny $C_{14}$} & \multicolumn{1}{c}{\tiny $C_{15}$} \\
\hhline{~------------------------------------}
{\tiny  1} & \cfl & \cfl & \cfl & \cfl & \cfl & & & & & & & & & & & & & & & & & \cfl & & \cfl & \cfl & & & & \cfl & & & & \cfl & & & \\
{\tiny  2} & \cfl & \cfl & \cfl & \cfl & \cfl & & & & & & & & & & & & & & & & & & \cfl & & & & \cfl & \cfl & & \cfl & & & & \cfl & & \\
{\tiny  3} & \cfl & \cfl & \cfl & \cfl & \cfl & & & & & & & & & & & & & & & & & & & & & \cfl & & & & & \cfl & \cfl & & & \cfl & \cfl \\
\hhline{~------------------------------------}
{\tiny  4} & & & & & \cfl & \cfl & \cfl & \cfl & \cfl & & & & & & & & & & & & & & \cfl & & & & \cfl & & \cfl & & & \cfl & & & \cfl & \\
{\tiny  5} & & & & & \cfl & \cfl & \cfl & \cfl & \cfl & & & & & & & & & & & & & & & & \cfl & \cfl & & \cfl & & & & & & \cfl & & \cfl \\
{\tiny  6} & & & & & \cfl & \cfl & \cfl & \cfl & \cfl & & & & & & & & & & & & & \cfl & & \cfl & & & & & & \cfl & \cfl & & \cfl & & & \\
\hhline{~------------------------------------}
{\tiny  7} & & & & & \cfl & & & & & \cfl & \cfl & \cfl & \cfl & & & & & & & & & \cfl & \cfl & & \cfl & & \cfl & & & & & & & & & \cfl \\
{\tiny  8} & & & & & \cfl & & & & & \cfl & \cfl & \cfl & \cfl & & & & & & & & & & & & & \cfl & & \cfl & & \cfl & \cfl & \cfl & & & & \\
{\tiny  9} & & & & & \cfl & & & & & \cfl & \cfl & \cfl & \cfl & & & & & & & & & & & \cfl & & & & & \cfl & & & & \cfl & \cfl & \cfl & \\
\hhline{~------------------------------------}
{\tiny 10} & & & & & \cfl & & & & & & & & & \cfl & \cfl & \cfl & \cfl & & & & & & \cfl & & \cfl & \cfl & & & & \cfl & & & & \cfl & & \\
{\tiny 11} & & & & & \cfl & & & & & & & & & \cfl & \cfl & \cfl & \cfl & & & & & & & \cfl & & & & \cfl & \cfl & & \cfl & & & & \cfl & \\
{\tiny 12} & & & & & \cfl & & & & & & & & & \cfl & \cfl & \cfl & \cfl & & & & & \cfl & & & & & \cfl & & & & & \cfl & \cfl & & & \cfl \\
\hhline{~------------------------------------}
{\tiny 13} & & & \cfl & & & & \cfl & & & & & & \cfl & \cfl & & & & \cfl & & & & & & & & \cfl & \cfl & & & & \cfl & & \cfl & \cfl & & \\
{\tiny 14} & & & \cfl & & & & \cfl & & & & & & \cfl & \cfl & & & & \cfl & & & & & \cfl & \cfl & & & & & & \cfl & & & & & \cfl & \cfl \\
{\tiny 15} & & & \cfl & & & & \cfl & & & & & & \cfl & \cfl & & & & \cfl & & & & \cfl & & & \cfl & & & \cfl & \cfl & & & \cfl & & & & \\
\hhline{~------------------------------------}
{\tiny 16} & & \cfl & & & & & & & \cfl & & & \cfl & & & \cfl & & & \cfl & & & & \cfl & \cfl & \cfl & \cfl & \cfl & & & & & & & & & & \\
{\tiny 17} & & \cfl & & & & & & & \cfl & & & \cfl & & & \cfl & & & \cfl & & & & & & & & & \cfl & \cfl & \cfl & \cfl & \cfl & & & & & \\
{\tiny 18} & & \cfl & & & & & & & \cfl & & & \cfl & & & \cfl & & & \cfl & & & & & & & & & & & & & & \cfl & \cfl & \cfl & \cfl & \cfl \\
\hhline{~------------------------------------}
{\tiny 19} & & & & \cfl & & & & \cfl & & & \cfl & & & & & \cfl & & \cfl & & & & & & & & \cfl & \cfl & & & \cfl & & & \cfl & & \cfl & \\
{\tiny 20} & & & & \cfl & & & & \cfl & & & \cfl & & & & & \cfl & & \cfl & & & & & \cfl & \cfl & & & & & & & \cfl & \cfl & & \cfl & & \\
{\tiny 21} & & & & \cfl & & & & \cfl & & & \cfl & & & & & \cfl & & \cfl & & & & \cfl & & & \cfl & & & \cfl & \cfl & & & & & & & \cfl \\
\hhline{~------------------------------------}
{\tiny 22} & \cfl & & & & & \cfl & & & & \cfl & & & & & & & \cfl & \cfl & & & & \cfl & & \cfl & & \cfl & & & & & \cfl & & & & \cfl & \\
{\tiny 23} & \cfl & & & & & \cfl & & & & \cfl & & & & & & & \cfl & \cfl & & & & & & & \cfl & & \cfl & & \cfl & \cfl & & & & & & \cfl \\
{\tiny 24} & \cfl & & & & & \cfl & & & & \cfl & & & & & & & \cfl & \cfl & & & & & \cfl & & & & & \cfl & & & & \cfl & \cfl & \cfl & & \\
\hhline{~------------------------------------}
{\tiny 25} & & & & \cfl & & \cfl & & & & & & \cfl & & \cfl & & & & & \cfl & & & & & \cfl & & & & \cfl & & \cfl & & & \cfl & & & \cfl \\
{\tiny 26} & & & & \cfl & & \cfl & & & & & & \cfl & & \cfl & & & & & \cfl & & & \cfl & & & & \cfl & & & \cfl & & & \cfl & & & \cfl & \\
{\tiny 27} & & & & \cfl & & \cfl & & & & & & \cfl & & \cfl & & & & & \cfl & & & & \cfl & & \cfl & & \cfl & & & & \cfl & & & \cfl & & \\
\hhline{~------------------------------------}
{\tiny 28} & \cfl & & & & & & & \cfl & & & & & \cfl & & \cfl & & & & \cfl & & & \cfl & & & \cfl & \cfl & & & & \cfl & & & & & \cfl & \\
{\tiny 29} & \cfl & & & & & & & \cfl & & & & & \cfl & & \cfl & & & & \cfl & & & & \cfl & & & & & \cfl & \cfl & & \cfl & & \cfl & & & \\
{\tiny 30} & \cfl & & & & & & & \cfl & & & & & \cfl & & \cfl & & & & \cfl & & & & & \cfl & & & \cfl & & & & & \cfl & & \cfl & & \cfl \\
\hhline{~------------------------------------}
{\tiny 31} & & & \cfl & & & & & & \cfl & \cfl & & & & & & \cfl & & & \cfl & & & \cfl & \cfl & \cfl & & & \cfl & & & & & \cfl & & & & \\
{\tiny 32} & & & \cfl & & & & & & \cfl & \cfl & & & & & & \cfl & & & \cfl & & & & & & \cfl & & & \cfl & & \cfl & \cfl & & & & \cfl & \\
{\tiny 33} & & & \cfl & & & & & & \cfl & \cfl & & & & & & \cfl & & & \cfl & & & & & & & \cfl & & & \cfl & & & & \cfl & \cfl & & \cfl \\
\hhline{~------------------------------------}
{\tiny 34} & & \cfl & & & & & \cfl & & & & \cfl & & & & & & \cfl & & \cfl & & & & & & \cfl & & & & & \cfl & \cfl & \cfl & \cfl & & & \\
{\tiny 35} & & \cfl & & & & & \cfl & & & & \cfl & & & & & & \cfl & & \cfl & & & \cfl & \cfl & & & & & & \cfl & & & & & \cfl & \cfl & \\
{\tiny 36} & & \cfl & & & & & \cfl & & & & \cfl & & & & & & \cfl & & \cfl & & & & & \cfl & & \cfl & \cfl & \cfl & & & & & & & & \cfl \\
\hhline{~------------------------------------}
{\tiny 37} & \cfl & & & & & & & & \cfl & & \cfl & & & \cfl & & & & & & \cfl & & & \cfl & \cfl & \cfl & & & \cfl & & & & & \cfl & & & \\
{\tiny 38} & \cfl & & & & & & & & \cfl & & \cfl & & & \cfl & & & & & & \cfl & & & & & & \cfl & \cfl & & \cfl & & \cfl & & & & & \cfl \\
{\tiny 39} & \cfl & & & & & & & & \cfl & & \cfl & & & \cfl & & & & & & \cfl & & \cfl & & & & & & & & \cfl & & \cfl & & \cfl & \cfl & \\
\hhline{~------------------------------------}
{\tiny 40} & & & & \cfl & & & \cfl & & & \cfl & & & & & \cfl & & & & & \cfl & & & & \cfl & & & & & \cfl & \cfl & & \cfl & & & & \cfl \\
{\tiny 41} & & & & \cfl & & & \cfl & & & \cfl & & & & & \cfl & & & & & \cfl & & \cfl & & & & \cfl & & \cfl & & & & & \cfl & \cfl & & \\
{\tiny 42} & & & & \cfl & & & \cfl & & & \cfl & & & & & \cfl & & & & & \cfl & & & \cfl & & \cfl & & \cfl & & & & \cfl & & & & \cfl & \\
\hhline{~------------------------------------}
{\tiny 43} & & \cfl & & & & \cfl & & & & & & & \cfl & & & \cfl & & & & \cfl & & & & & \cfl & & & & \cfl & & \cfl & \cfl & & \cfl & & \\
{\tiny 44} & & \cfl & & & & \cfl & & & & & & & \cfl & & & \cfl & & & & \cfl & & \cfl & \cfl & & & & & & & \cfl & & & \cfl & & & \cfl \\
{\tiny 45} & & \cfl & & & & \cfl & & & & & & & \cfl & & & \cfl & & & & \cfl & & & & \cfl & & \cfl & \cfl & \cfl & & & & & & & \cfl & \\
\hhline{~------------------------------------}
{\tiny 46} & & & \cfl & & & & & \cfl & & & & \cfl & & & & & \cfl & & & \cfl & & & & \cfl & \cfl & \cfl & & & \cfl & & & & & \cfl & & \\
{\tiny 47} & & & \cfl & & & & & \cfl & & & & \cfl & & & & & \cfl & & & \cfl & & \cfl & & & & & \cfl & \cfl & & \cfl & & \cfl & & & & \\
{\tiny 48} & & & \cfl & & & & & \cfl & & & & \cfl & & & & & \cfl & & & \cfl & & & \cfl & & & & & & & & \cfl & & \cfl & & \cfl & \cfl \\
\hhline{~------------------------------------}
{\tiny 49} & & \cfl & & & & & & \cfl & & \cfl & & & & \cfl & & & & & & & \cfl & \cfl & & & & & & \cfl & & & \cfl & & & \cfl & & \cfl \\
{\tiny 50} & & \cfl & & & & & & \cfl & & \cfl & & & & \cfl & & & & & & & \cfl & & & \cfl & \cfl & & \cfl & & & & & & \cfl & & \cfl & \\
{\tiny 51} & & \cfl & & & & & & \cfl & & \cfl & & & & \cfl & & & & & & & \cfl & & \cfl & & & \cfl & & & \cfl & \cfl & & \cfl & & & & \\
\hhline{~------------------------------------}
{\tiny 52} & & & \cfl & & & \cfl & & & & & \cfl & & & & \cfl & & & & & & \cfl & & \cfl & & & & & \cfl & \cfl & & & & & & \cfl & \cfl \\
{\tiny 53} & & & \cfl & & & \cfl & & & & & \cfl & & & & \cfl & & & & & & \cfl & & & & \cfl & \cfl & \cfl & & & & & \cfl & \cfl & & & \\
{\tiny 54} & & & \cfl & & & \cfl & & & & & \cfl & & & & \cfl & & & & & & \cfl & \cfl & & \cfl & & & & & & \cfl & \cfl & & & \cfl & & \\
\hhline{~------------------------------------}
{\tiny 55} & \cfl & & & & & & \cfl & & & & & \cfl & & & & \cfl & & & & & \cfl & & \cfl & \cfl & & \cfl & & \cfl & & & & \cfl & & & & \\
{\tiny 56} & \cfl & & & & & & \cfl & & & & & \cfl & & & & \cfl & & & & & \cfl & \cfl & & & & & \cfl & & \cfl & & \cfl & & \cfl & & & \\
{\tiny 57} & \cfl & & & & & & \cfl & & & & & \cfl & & & & \cfl & & & & & \cfl & & & & \cfl & & & & & \cfl & & & & \cfl & \cfl & \cfl \\
\hhline{~------------------------------------}
{\tiny 58} & & & & \cfl & & & & & \cfl & & & & \cfl & & & & \cfl & & & & \cfl & \cfl & & & & & \cfl & \cfl & & & & & & \cfl & \cfl & \\
{\tiny 59} & & & & \cfl & & & & & \cfl & & & & \cfl & & & & \cfl & & & & \cfl & & & \cfl & \cfl & & & & & & \cfl & \cfl & & & & \cfl \\
{\tiny 60} & & & & \cfl & & & & & \cfl & & & & \cfl & & & & \cfl & & & & \cfl & & \cfl & & & \cfl & & & \cfl & \cfl & & & \cfl & & & \\
\hhline{~------------------------------------}
{\tiny 61} & & & & & \cfl & & & & & & & & & & & & & \cfl & \cfl & \cfl & \cfl & \cfl & \cfl & & & \cfl & & & & & \cfl & & & & & \cfl \\
{\tiny 62} & & & & & \cfl & & & & & & & & & & & & & \cfl & \cfl & \cfl & \cfl & & & \cfl & & & \cfl & & \cfl & \cfl & & & & \cfl & & \\
{\tiny 63} & & & & & \cfl & & & & & & & & & & & & & \cfl & \cfl & \cfl & \cfl & & & & \cfl & & & \cfl & & & & \cfl & \cfl & & \cfl & \\
\hhline{~------------------------------------}
\end{tabular}
\end{tabular}
\end{center}
\caption{An example of a $(21 \times 15, 63)$-triple array and a representation of its underlying unordered triple array.}
\label{fig:ex21-15-63}
\end{figure}

\clearpage
\section{Counts of extremal (unordered) triple arrays}\label{ap:extr}

\begin{table}[!ht]
\begin{center}
\begin{tabular}{|rr|r|r|r|r|r|r|}
\hline
\multicolumn{2}{|r|}{$r \times c$} & $3 \times 4$ & $4 \times 9$ & $5 \times 6$ & $5 \times 16$ & $6 \times 10$ & $7 \times 8$ \\
\multicolumn{2}{|r|}{$v$} & $6$ & $12$ & $10$ & $20$ & $15$ & $14$ \\
\hline
\multicolumn{8}{c}{\multirow{2}{*}{Unordered triple arrays}} \\
\multicolumn{8}{c}{} \\
\hline
\multicolumn{2}{|c|}{Total \#} & $1$ & $1$ & $1$ & $1$ & $3$ & $10$ \\
\hline
$|\Aut|$ & $12$ &  &  &  &  &  & $2$ \\
 & $16$ &  &  &  &  &  & $1$ \\
 & $21$ &  &  &  &  &  & $1$ \\
 & $24$ & $1$ &  &  &  & $1$ & $1$ \\
 & $48$ &  &  &  &  & $1$ & $1$ \\
 & $60$ &  &  & $1$ &  &  &  \\
 & $96$ &  &  &  &  &  & $1$ \\
 & $168$ &  &  &  &  &  & $1$ \\
 & $192$ &  &  &  &  &  & $1$ \\
 & $432$ &  & $1$ &  &  &  &  \\
 & $720$ &  &  &  &  & $1$ &  \\
 & $1344$ &  &  &  &  &  & $1$ \\
 & $5760$ &  &  &  & $1$ &  &  \\
\hline
\multicolumn{8}{c}{\multirow{2}{*}{Triple arrays}} \\
\multicolumn{8}{c}{} \\
\hline
\multicolumn{2}{|c|}{Total \#} & $0$ & $1$ & $7$ & $26804$ & $270119$ & $684782$ \\
\hline
$|\Aut|$ & $1$ &  &  &  & $26714$ & $263790$ & $682054$ \\
 & $2$ &  &  &  &  & $5280$ & $1266$ \\
 & $3$ &  & $1$ & $2$ & $90$ & $260$ & $1277$ \\
 & $4$ &  &  & $1$ &  & $579$ & $98$ \\
 & $5$ &  &  &  &  & $1$ &  \\
 & $6$ &  &  & $1$ &  & $69$ & $48$ \\
 & $7$ &  &  &  &  &  & $2$ \\
 & $8$ &  &  &  &  & $88$ & $12$ \\
 & $10$ &  &  &  &  & $2$ &  \\
 & $12$ &  &  & $2$ &  & $17$ & $9$ \\
 & $16$ &  &  &  &  & $11$ &  \\
 & $18$ &  &  &  &  & $1$ &  \\
 & $20$ &  &  &  &  & $4$ &  \\
 & $21$ &  &  &  &  &  & $8$ \\
 & $24$ &  &  &  &  & $9$ & $7$ \\
 & $36$ &  &  &  &  & $2$ &  \\
 & $48$ &  &  &  &  & $4$ &  \\
 & $60$ &  &  & $1$ &  &  &  \\
 & $120$ &  &  &  &  & $1$ &  \\
 & $168$ &  &  &  &  &  & $1$ \\
 & $720$ &  &  &  &  & $1$ &  \\
\hline
\end{tabular}
\end{center}
\caption{The number of non-isomorphic unordered triple arrays sorted by automorphism group order and non-isotopic triple arrays sorted by autotopism group order.}
\label{tbl:extr-UTA+TA}
\end{table}

\begin{table}[!ht]
\begin{center}
\begin{tabular}{c}
\begin{tabular}{|rr|r|r|r|r|r|r|r|}
\hline
\multicolumn{2}{|r|}{$r \times c$} & $6 \times 25$ & $8 \times 49$ & $9 \times 10$ & $9 \times 16$ & $9 \times 28$ & $9 \times 64$ & $10 \times 21$ \\
\multicolumn{2}{|r|}{$v$} & $30$ & $56$ & $18$ & $24$ & $36$ & $72$ & $30$ \\
\hline
\multicolumn{2}{|c|}{Total \#} & $1$ & $1$ & $22$ & $1382$ & $8$ & $1$ & $3809$ \\
\hline
$|\Aut|$ & $1$ &  &  & $4$ & $1172$ &  &  & $3648$ \\
 & $2$ &  &  & $7$ & $136$ & $1$ &  & $21$ \\
 & $3$ &  &  &  & $22$ &  &  & $125$ \\
 & $4$ &  &  &  & $36$ &  &  &  \\
 & $6$ &  &  & $2$ & $8$ & $1$ &  & $4$ \\
 & $7$ &  &  &  &  &  &  & $1$ \\
 & $8$ &  &  & $6$ & $2$ &  &  &  \\
 & $9$ &  &  & $2$ &  & $1$ &  & $5$ \\
 & $12$ &  &  &  & $4$ &  &  &  \\
 & $14$ &  &  &  &  &  &  & $1$ \\
 & $21$ &  &  &  &  &  &  & $4$ \\
 & $24$ &  &  &  & $2$ &  &  &  \\
 & $42$ &  &  &  &  & $1$ &  &  \\
 & $54$ &  &  &  &  & $2$ &  &  \\
 & $72$ &  &  & $1$ &  &  &  &  \\
 & $168$ &  &  &  &  & $1$ &  &  \\
 & $1512$ &  &  &  &  & $1$ &  &  \\
 & $12000$ & $1$ &  &  &  &  &  &  \\
 & $98784$ &  & $1$ &  &  &  &  &  \\
 & $677376$ &  &  &  &  &  & $1$ &  \\
\hline
\end{tabular}\\
\\
\begin{tabular}{|rr|r|r|r|r|r|}
\hline
\multicolumn{2}{|r|}{$r \times c$} & $10 \times 81$ & $11 \times 12$ & $11 \times 45$ & $13 \times 14$ & $13 \times 66$ \\
\multicolumn{2}{|r|}{$v$} & $90$ & $22$ & $55$ & $26$ & $78$ \\
\hline
\multicolumn{2}{|c|}{Total \#} & $7$ & $23360$ & $16$ & $5606594$ & $\geq 8$ \\
\hline
$|\Aut|$ & $1$ &  & $22095$ & $1$ & $5596116$ &  \\
 & $2$ &  & $920$ & $2$ & $7982$ & $\geq 2$ \\
 & $3$ &  & $199$ & $1$ & $2316$ &  \\
 & $4$ &  & $57$ & $4$ & $96$ &  \\
 & $5$ &  & $34$ &  &  & $\geq 1$ \\
 & $6$ &  & $34$ &  & $77$ &  \\
 & $8$ &  &  & $2$ &  &  \\
 & $10$ &  & $2$ &  &  & $\geq 2$ \\
 & $11$ &  & $2$ &  &  &  \\
 & $12$ &  & $8$ &  &  &  \\
 & $13$ &  &  &  & $2$ &  \\
 & $16$ &  &  & $3$ &  &  \\
 & $39$ &  &  &  & $4$ &  \\
 & $55$ &  & $2$ &  &  & $\geq 1$ \\
 & $60$ &  & $6$ &  &  &  \\
 & $72$ &  &  & $1$ &  &  \\
 & $78$ &  &  &  & $1$ &  \\
 & $110$ &  &  &  &  & $\geq 2$ \\
 & $144$ &  &  & $1$ &  &  \\
 & $432$ & $1$ &  &  &  &  \\
 & $660$ &  & $1$ &  &  &  \\
 & $1440$ &  &  & $1$ &  &  \\
 & $2592$ & $1$ &  &  &  &  \\
 & $3456$ & $1$ &  &  &  &  \\
 & $3840$ & $1$ &  &  &  &  \\
 & $31104$ & $1$ &  &  &  &  \\
 & $311040$ & $1$ &  &  &  &  \\
 & $933120$ & $1$ &  &  &  &  \\
\hline
\end{tabular}
\end{tabular}
\end{center}
\caption{The number of non-isomorphic unordered triple arrays sorted by automorphism group order.}
\label{tbl:extr-UTA-only}
\end{table}

\clearpage

\begin{table}[!ht]
\begin{center}
\begin{tabular}{|rr|rrr|r|}
\hline
\multicolumn{2}{|c|}{UTA} & $U_{0}$ & $U_{1}$ & $U_{2}$ & All \\
\multicolumn{2}{|c|}{$|\Aut U_i|$} & $24$ & $48$ & $720$ & \\
\hline
\multicolumn{2}{|c|}{Total \#} & $162202$ & $96890$ & $11027$ & $270119$ \\
\hline
$|\Aut|$  & $1$ & $161827$ & $93977$ & $7986$ & $263790$ \\
 & $2$ & $158$ & $2650$ & $2472$ & $5280$ \\
 & $3$ & $211$ & $17$ & $32$ & $260$ \\
 & $4$ &  & $212$ & $367$ & $579$ \\
 & $5$ &  &  & $1$ & $1$ \\
 & $6$ & $6$ & $14$ & $49$ & $69$ \\
 & $8$ &  & $18$ & $70$ & $88$ \\
 & $10$ &  &  & $2$ & $2$ \\
 & $12$ &  &  & $17$ & $17$ \\
 & $16$ &  &  & $11$ & $11$ \\
 & $18$ &  &  & $1$ & $1$ \\
 & $20$ &  &  & $4$ & $4$ \\
 & $24$ &  & $2$ & $7$ & $9$ \\
 & $36$ &  &  & $2$ & $2$ \\
 & $48$ &  &  & $4$ & $4$ \\
 & $120$ &  &  & $1$ & $1$ \\
 & $720$ &  &  & $1$ & $1$ \\
\hline
\end{tabular}
\end{center}
\caption{The number of non-isotopic $(6 \times 10, 15)$-triple arrays sorted by the underlying unordered triple array and autotopism group order.
The unordered triple arrays are numbered in the order that they are given in~\cite{gordeevZenodoLibrary25}.}
\label{tbl:61015}
\end{table}

\begin{table}[!ht]
\begin{center}
\setlength\tabcolsep{3pt}
\begin{tabular}{|rr|rrrrrrrrrr|r|r|}
\hline
\multicolumn{2}{|c|}{UTA} & $U_{0a}$ & $U_{0b}$ & $U_{1a}$ & $U_{1b}$ & $U_{1c}$ & $U_{2a}$ & $U_{2b}$ & $U_{3}$ & $U_{4a}$ & $U_{4b}$ & RTA & All \\
\multicolumn{2}{|c|}{Resolvable?} & $+$ & $-$ & $+$ & $-$ & $-$ & $-$ & $-$ & $+$ & $+$ & $-$ & & \\
\multicolumn{2}{|c|}{$|\Aut U_{ij}|$} & $168$ & $12$ & $96$ & $16$ & $12$ & $24$ & $21$ & $1344$ & $192$ & $48$ & & \\
\hline
\multicolumn{2}{|c|}{Total \#} & $9968$ & $148676$ & $18574$ & $118760$ & $165804$ & $71694$ & $86826$ & $3096$ & $11750$ & $49634$ & $43388$ & $684782$ \\
\hline
$|\Aut|$  & $1$ & $9861$ & $148468$ & $18433$ & $118760$ & $165514$ & $71563$ & $86588$ & $2248$ & $11054$ & $49565$ & $41596$ & $682054$ \\
 & $2$ &  &  &  &  &  &  &  & $659$ & $607$ &  & $1266$ & $1266$ \\
 & $3$ & $102$ & $208$ & $141$ &  & $290$ & $131$ & $235$ & $50$ & $51$ & $69$ & $344$ & $1277$ \\
 & $4$ &  &  &  &  &  &  &  & $64$ & $34$ &  & $98$ & $98$ \\
 & $6$ &  &  &  &  &  &  &  & $48$ &  &  & $48$ & $48$ \\
 & $7$ & $2$ &  &  &  &  &  &  &  &  &  & $2$ & $2$ \\
 & $8$ &  &  &  &  &  &  &  & $12$ &  &  & $12$ & $12$ \\
 & $12$ &  &  &  &  &  &  &  & $5$ & $4$ &  & $9$ & $9$ \\
 & $21$ & $3$ &  &  &  &  &  & $3$ & $2$ &  &  & $5$ & $8$ \\
 & $24$ &  &  &  &  &  &  &  & $7$ &  &  & $7$ & $7$ \\
 & $168$ &  &  &  &  &  &  &  & $1$ &  &  & $1$ & $1$ \\
\hline
\end{tabular}
\end{center}
\caption{The number of non-isotopic $(7 \times 8, 14)$-triple arrays sorted by the underlying unordered triple array and autotopism group order.
Column RTA lists the number of non-isotopic resolvable $(7 \times 8, 14)$-triple arrays sorted by autotopism group order.
The unordered triple arrays are given, from left to right, in the order they appear in~\cite{gordeevZenodoLibrary25}.
They have the same first index, for example $U_{0a}$ and $U_{0b}$, if they come from the same symmetric 2-design in \nameref{constr:Agrawal}.}
\label{tbl:7814}
\end{table}

\end{document}